\documentclass[12pt]{amsbook} 
\usepackage
{amsmath}
\usepackage{amscd, amsthm, amsfonts, amssymb}
\usepackage{diagrams}
\usepackage{mathptmx}
\newtheorem{thm}{Theorem}[section] 
\newtheorem{cor}[thm]{Corollary}
\newtheorem{lem}[thm]{Lemma}
\newtheorem{prop}[thm]{Proposition}
 
\newtheorem{defn}[thm]{Definition}

\newtheorem{preremark}[thm]{Remark}
\newenvironment{remark}%
  {\begin{preremark}\upshape}{\end{preremark}}
\newtheorem{preexample}[thm]{Example}
\newenvironment{example}%
  {\begin{preexample}\upshape}{\end{preexample}}

\numberwithin{thm}{section}
\numberwithin{equation}{section}
\newenvironment{Myenumerate}
{

\begin{enumerate}}{\end{enumerate}}

\newcommand{\ad}{\operatorname{ad}}

\newcommand{\CToda}{\mathrm{CToda}}
\newcommand{\Cg}{C{\mathfrak g}}
\newcommand{\Cz}{\mathbb C_{\mathbb Z}} 
\newcommand{\Czpol}{\mathbb C_{\mathbb Z}^{\mathrm{pol}}} 
\newcommand{\Cztor}{\mathbb C_{\mathbb Z}^{\mathrm{tor}}} 
\newcommand{\Cztf}{\mathbb C_{\mathbb Z}^{\mathrm{tf}}} 
\newcommand{\hatCzpol}{\hat{\mathbb C}_{\mathbb Z}^{\mathrm{pol}}} 
\newcommand{\Der}{\operatorname{\mathrm{Der}}}
\newcommand{\Derev}{\operatorname{\mathrm{Der}^{\mathrm{ev}}}}

\newcommand{\End}{\operatorname{End}}

\newcommand{\Hol}{\mathrm{Hol}}
\newcommand{\Hom}{\operatorname{\mathrm{Hom}}} 
\newcommand{\inv}{^{-1}}
\newcommand{\K}{\operatorname{\mathrm{K}}}
\newcommand{\KZ}{\operatorname{\mathrm{K}_{\mathbb Z}}}

\newcommand{\barKsing}{\overline \K_{\Sing}}
\newcommand{\Ksing}{\operatorname{K}_{\mathrm{Sing}}}
\newcommand{\LToda}{\mathcal{L}\mathrm{Toda}}
\newcommand{\LTodaann}{\mathcal{L}\mathrm{Toda}_{\mathrm{ann}}}
\newcommand{\LTodacr}{\mathcal{L}\mathrm{Toda}_{\mathrm{cr}}}
\newcommand{\Lg}{\mathcal{L}\mathfrak g}
\newcommand{\Lgcr}{\mathcal{L}\mathfrak g_{\mathrm{cr}}}
\newcommand{\Lgann}{\mathcal{L}\mathfrak g_{\mathrm{ann}}}

\newcommand{\md}{\operatorname{m}^*}

\newcommand{\mStwod}{\operatorname{m}_{1\otimes S}^*}
\newcommand{\poissontensor}{\mathop{\overset{\otimes},}}

\newcommand{\otimesd}{\dot\otimes}
\newcommand{\Res}{\operatorname{\mathrm{Res}}}
\newcommand{\Sing}{{\mathrm{Sing}}}
\newcommand{\Sp}{\operatorname{\mathrm{Sp}}}

\newcommand{\Specm}{\operatorname{\mathrm{Specm}}}
\newcommand{\Tr}{\operatorname{\mathrm{Tr}}}

\newcommand{\VToda}{\mathrm{VToda}}
\newcommand{\Yhol}{{Y_{\Hol}}}
\newcommand{\YholS}{{Y^S_{\Hol}}}
\newcommand{\Ysing}{{Y_{\Sing}}}
\newcommand{\YsingLC}{{Y^{LC}_{\Sing}}}
\newcommand{\YsingC}{{Y^{C}_{\Sing}}}
\newcommand{\YsingS}{{Y_{\Sing}^S}}
\begin{document}
\thispagestyle{empty}
\vspace*{\fill}
\begin{center}
\Huge\textsf{$H_T$- Vertex Algebras }\par
\end{center}
\begin{center}
\Large\textsf{and}\par
\end{center}
\begin{center}
\Huge\textsf{ the Infinite Toda  Lattice}\par
\end{center}
                                                                                                                                            
\begin{center}
\Large\textsf{M.J. Bergvelt}\par
\end{center}
\vspace*{\fill}
\begin{center}
\textsf{Department of Mathematics} \\
\textsf{University of Illinois at Urbana-Champaign}\\
\textsf{email: bergv@uiuc.edu}
\end{center}
\clearpage

         \begin{abstract} Let $H_T=\mathbb C[T, T\inv]$ be the Hopf
           algebra of symmetries of a lattice of rank 1, or
           equivalently $H_T$ is the group algebra of a free Abelian
           group with one generator $T$. We construct conformal
           algebras, vertex Poisson algebras and vertex algebras with
           $H_T$ as symmetry. For example, the Hamiltonian structure
           for the infinite Toda lattice gives rise to an $H_T$-vertex
           Poisson structure on a free difference algebra. Examples of
           $H_T$-vertex algebras are constructed from representations
           of a class of infinite dimensional Lie algebras related to
           $H_T$ in the same way loop algebras are related to the Hopf
           algebra $H_D=\mathbb C[D]$ of infinitesimal translations
           used in the usual vertex algebras.
          \end{abstract}

\tableofcontents

\chapter{Introduction}
\label{chap:Introduction}

\section{Classical Vertex Algebras and the Hopf Algebra $H_D$ }

\label{sec:vertalg}

In this section we recall what a vertex algebra in the usual sense
is and indicate how it is related to a Hopf algebra $H_D$.  See e.g.,
Frenkel-Lepowsky-Meurman (\cite{MR90h:17026}) or Kac
(\cite{MR99f:17033} for more details on classical vertex algebras and
Borcherds (\cite{MR1653021}) on the idea of emphasizing the role of
the Hopf algebra, cf., also Snydal (\cite{math.QA/9904104}).

Briefly, a vertex algebra structure on a vector
space $V$ is a map $Y\colon V\otimes V\to V((z))$ associating to a pair of
vectors $a,b$ in $V$ a formal Laurent series $Y(a\otimes b,z)$ with  
values in $V$.  As usual we write $Y(a\otimes b,z)=Y(a,z)b$ to obtain
the \emph{vertex operator} $Y(a,z)=\sum a_{(n)}z^{-n-1}$ with the
$a_{(n)}$ linear maps $V\to V$.  Also an essential ingredient is the
\emph{infinitesimal translation operator} $D\colon V\to V$, such that
$[D,Y(a,z)]=\partial_z Y(a,z)$.  These ingredients satisfy a number of
axioms that we will not recall at this point.  One of the properties
of a vertex algebra is that two vertex operators satisfy a commutator
formula of the form
\begin{equation}
  [Y(a,z),Y(b,w)]=\sum_{n=0}^N Y(a_{(n)}b, w)\partial^{(n)}_w\delta(z,w),
\label{eq:vertexcommut}
\end{equation}
where $\delta(z,w)=\sum_{k\in\mathbb Z} z^kw^{-k-1}$ is the formal
Dirac delta function, and for any linear operator $D$ we write
$D^{(n)}=\frac{D^n}{n!}$ .

The simplest nontrivial example is given by the vertex algebra of the Heisenberg
algebra, the Lie algebra spanned by a central element $c$ and elements
$b_{(n)}$, $n\in \mathbb Z$, with relations
\[
[b_{(m)}, b_{(n)}]=m\delta_{m+n}c.
\]
The vector space $V=\mathbb C[b_{(-n-1)}]_{n\ge 0}$ is then a vertex
algebra, with basic vertex operator the free boson
$b(z)=Y(b_{(-1)},z)=\sum_{n\in \mathbb Z} b_{(n)}z^{-n-1}$. All other
vertex operators in $V$ are obtained by taking linear combinations of
\emph{normal ordered products} of the free boson. The commutator of the free
boson with itself is
\begin{equation}
  \label{eq:freeboson}
  [b(z),b(w)]=\partial_w\delta(z,w).
\end{equation}
Let $H_D=\mathbb C[D]$ be the algebra of polynomials in the linear
operator $D$.  So by definition every vertex algebra is a module over
$H_{D}$.  We can think of $H_D$ as the universal enveloping algebra of
the 1-dimensional Abelian Lie algebra with basis element $D$.  As such
$H_D$ is a Hopf-algebra, with coproduct
 $\Delta_{H}\colon  D\mapsto D\otimes
1+1\otimes D$, antipode $S\colon D\mapsto -D$ and counit $\epsilon\colon  D\mapsto 0$.
Now the dual Hopf algebra of $H_D$ is the Hopf algebra $H_D^*=\mathbb
C[[t]]$ 
of formal power series in a variable $t$, the functions on the
formal disk.  Note that the map $Y$ in a vertex algebra structure
takes values in formal Laurent series, which are the \emph{singular}
functions on the formal disk.  In the same spirit note that the
commutator is \emph{singular}, in the sense that we cannot put $z=w$
in this formula.  The delta function can be obtained from the
difference of two expansions of the basic singularity $\frac 1{z-w}$.

The occurence of $H_D$ and the singularities in its dual is no
coincidence but the tip of an iceberg.  Borcherds, \cite{MR1653021}
has a general construction of what he calls $G$-vertex algebras. The
construction is based on the notion of a vertex group $G$, which is,
roughly speaking,
\begin{itemize}
\item the choice of a (cocommutative) Hopf algebra $H$,
\item the choice of singularities; this is in examples some
  localization $K$ of the dual (commutative) algebra $H^{*}$.
\end{itemize}
The classical notion of a vertex algebra corresponds to the choice of
$H=H_D$ and $K=K_D=\mathbb C[[t]][t\inv]$. Borcherds' work and examples
are very interesting, but lacks the level of detail we know about the
usual vertex algebras, as developed for instance in Frenkel, Lepowsky
and Meurman (\cite{MR90h:17026}) or Kac (\cite{MR99f:17033}).

In this paper we study in detail a new class of vertex algebras.
These are non-trivial examples of $G$-vertex algebras of the sort
envisioned by Borcherds.  We choose as a Hopf algebra not, as in the
classical case, the universal enveloping algebra $H_D=\mathcal U(L)$ of an
1-dimensional Abelian Lie algebra $L=\mathbb C D$, but the group
algebra 
\begin{equation}
  \label{eq:defH_T}
H_T=\mathbb CA=\mathbb C[T,T^{-1}]  
\end{equation}
of a free rank one Abelian group $A=\langle T\rangle$ generated by an
element $T$.  The dual algebra is the algebra $\Cz$ of functions on
the integers, and we have found a natural way to introduce
singularities for $\Cz$ by localizing to get an algebra $\K$.  This
leads to a class of what we call $H_T$-vertex algebras.  These occurs
naturally in the theory of the infinite Toda Lattice.  To explain this
we first recall some facts about Gelfand-Dickey Hamiltonian
structures, and how they are related to classical vertex algebras.

\section{Gelfand-Dickey Structures and the Hopf Algebra $H_D$}
\label{sec:GDHamstructH_D}

In the simplest case one wants to describe an infinite hierarchy of
differential equations for a single unknown $v$ depending on spatial
variable $x$ and infinitely many time variables $t_{i}$, $i>0$. 
For instance, the Modified Korteweg-de Vries (mKdV) hierarchy is a system of
differential equations of the form
\begin{equation}
  \label{eq:mkdv}
\partial_{t_i}v=f_i(v,v^\prime,\dots), \quad v^{\prime}=\partial_{x}v,
\quad i>0,
\end{equation}
where $f_i$ is a differential polynomial, i.e., an element of the
polynomial algebra $V=\mathbb C[v^{(n)}]_{n\ge0}$. This algebra is a
module over $H_D$ via $Dv^{(n)}=v^{(n+1)}$. $V$ is a commutative
algebra, and its (maximal) spectrum
  $\operatorname{Specm}(V)=\operatorname{Hom}_{\mathbb
C-\operatorname{alg}}(V,\mathbb C)$
 can be identified with $\mathbb
C[[t]]$. Indeed, a formal power series $s\in \mathbb C[[t]]$ acts on
an element $f=f(v,v^{(1)},\dots)\in V$ by
\[
s(f)=\Res \left( \frac{ f(s,s^\prime,\dots)}t\right ),\quad
s^{\prime}=\partial_{t}s,
\]
and we have $s(fg)=s(f)s(g)$.  Of course, usually we will think of $f$
as a function on the phase space $\Specm(V)$, and we will write
$f[s]=s(f)$. More generally, one can define functionals
\begin{equation}\label{eq:functionals}
f_{(n)}[s]=\Res\left( f(s,s^\prime,\dots) t^n\right)
\end{equation}
so that $f[s]=f_{(-1)}[s]$. For $n\ge0$ the functionals $f_{(n)}$ are
identically zero on $s\in\mathbb C[[t]]$, so we will from now on allow
$s\in\mathbb C((t))$. The reader will be excited to note that here the
same space of singular functions on the formal disk appears as in the
previous section in the context of vertex algebras. But there is more
\dots

Define the \emph{Classical Field} $C(f,z)$ of $f\in V$ as the
generating series of functionals
\[
C(f,z)=\sum_{n\in \mathbb Z}f_{(n)} z^{-n-1}.
\]
A Gelfand-Dickey Hamiltonian structure on $V$ is then a prescription
to define the Poisson brackets of functionals $f_{(n)}$, or,
equivalently, for the corresponding classical fields.  For a precise
definition (of the Poisson bracket $\{f_{(0)},g_{(0)}\}$ of the
residues of the classical fields) see Dickey's book
(\cite{MR1964513}).  The simplest nontrivial example is, perhaps, the
Hamiltonian structure for the modified Korteweg-de Vries hierarchy.
In this case the Poisson bracket for two functionals is given by
\[
\{ f_{(m)},g_{(n)}\}=\Res\left( D\big(\frac{\delta f_{(m)}}{\delta
    v}\big )\frac{\delta g_{(n)}}{\delta v}\right),
\]
where the variational derivative of a functional is given by
\[
\frac{\delta f_{(m)}}{\delta v}=\sum_{i=0}^{\infty} (-D)^{i}(\frac
{\partial f}{\partial v^{(i)}}\otimes t^{m}),\quad
D=D_{V}+\partial_{t}.
\]
Then one checks that the Poisson bracket of the basic classical field
$C(v,z)$ is given by
\[
\{C(v,z),C(v,w)\}=\partial_w \delta(z,w).
\]
More generally, the Poisson bracket of two classical fields will be a
linear combination of derivatives of delta functions, with classical
fields as coefficients.  Of course, this is very similar to the
commutator formula (\ref{eq:freeboson}) for free bosons.

The basic reason $H_D$ appears in the theory of equations of the form
(\ref{eq:mkdv}) is that if $v(x)$ solves this system, also the
infinitesimal translate $v(x+\delta)$ will solve it.  In the theory of
Toda lattices one studies equations of the following form (or the
obvious multicomponent generalization)
\begin{equation}
  \label{eq:infTodaprototype}\partial_{t_j} x_{n}= f_j(\dots, x_{n-1},x_n,x_{n+1},\dots),
\end{equation}
where $x=(x_n(t))_{n\in\mathbb Z}$ is a sequence (of functions of the
times $t_j$, $j>0$) and $f_j$ is a polynomial in the entries of the sequence
$x$.  In this case we get as basic symmetry invariance under shifts:
if $x(t)$ solves (\ref{eq:infTodaprototype}) then also the shifted
sequence $Tx:=(x_{n+1}(t))_{n\in \mathbb Z}$ is a solution.  In this
case the symmetry algebra is the Hopf algebra $H_T$ of
\eqref{eq:defH_T}. 

First we will develop some of the machinery of the
Hopf algebra $H_{T}$ in the next chapter. Then we use this to define
in chapter \ref{chap:SingHamStruc} the notion of a Singular Hamiltonian
structure and the associated Poisson brackets. As an example we discuss
the Hamiltonian structure for the infinite Toda lattice in chapter
\ref{chap:HamstrucInfToda}. Then we define a generalization of singular
Hamiltonian structures called $H_T$-conformal algebras in
\ref{chap:H_TConformalAlg}. Finally we define a class of vertex algebras
associated to $H_{T}$ in chapter \ref{chap:VertexAlg} and discuss some
simple examples.

\chapter{The Hopf Algebra $H_{T}$ and Sequences}
\label{chap:HopfAlgebraSeq}
\section{Introduction}
\label{sec:introHopf}

The goal of this chapter is to indicate how one can construct for the
Hopf algebra $H_{T}$ analogs of well known structures related to
formal series used in classical vertex algebras. Recall that the dual
of the symmetry algebra $H_{D}=\mathbb C[D]$ is just the power series
ring $\mathbb C[[t]]$. Among the structures we have in mind here are

\begin{itemize}
\item The localization of the dual. This gives the formal Laurent
  series algebra $K_{D}=\mathbb C((t))$ from $H_{D}^{*}=\mathbb
  C[[t]]$ by inverting $t$.
  
\item The residue map $\Res_0\colon  K_{D}\to \mathbb C$, characterized by
  $\operatorname{Res}_{0}(Df)=0$ and the normalization
  $\operatorname{Res}_{0}(\frac 1t)=1$.  Note that the counit $\epsilon$
  on $H_{D}$ is given by $\epsilon(D)=0$, so that we can write
  $\operatorname{Res}_{0}(Df)=\epsilon(D)\operatorname{Res}_{0}(f)$.
  
\item The expansion maps (in regions $|x|>|y|$, resp. $|y|>|x|$)
  \begin{equation}
    \label{eq:expansionH_D}
    \frac1{x-y}\mapsto e^{-y\partial_{x}}\frac1x,\quad\frac1{x-y}\mapsto 
e^{-x\partial_{y}}\frac1{-y}.
  \end{equation}
  
\item The formal delta function, the difference of the two expansions:
\[
\delta(x,y)=e^{-y\partial_{x}}\frac1x + e^{-x\partial_{y}}\frac1y.
\]  
\end{itemize}
For $H_{T}=\mathbb C[T,T\inv]$ the dual Hopf algebra is the algebra
$\Cz$ of (arbitrary) functions from the integers to $\mathbb C$ as
will be explained below. In other words the elements of $\Cz$ are
sequences $s=(s_{n})_{n\in \mathbb Z}$, with $s_n\in \mathbb C$.  Let
$\delta_{n}$ be the Kronecker sequence with value 1 at $n$ and zero
elsewhere, so that $s=\sum_{n\in \mathbb Z}s_{n}\delta_{n}$.  Note
that Kronecker sequences are zero divisors in $\Cz$, as are all non
invertible elements of $\Cz$..

Because of the zero divisors in $\Cz$ localizing is a somewhat
delicate matter. In general, we want to choose a multiplicative set
$M\subset \Cz$ of elements that can be inverted, but we must be
careful not to have $f, g\in M$ with $f.g=0$ , because
in that case $1/f=g/f.g$.  In particular inverting more than one of
the Kronecker sequences $\delta_{n}$ would not be a useful strategy.
It turns out that we need to find appropriate analogs in $\Cz$ of the
powers of $t$ occuring in the dual $\mathbb C[[t]]$ of $H_D$. These
are the sequences $\tau(\ell)$ introduced in section
\ref{sec:HopfDual}.

\section{The Hopf Algebra $H_T$}
\label{sec:TodaSymHopf}
Let $H_T=\mathbb C[T, T\inv]$, as before, so that $H_T$ is a
commutative and cocommutative Hopf algebra and we have the following structures on $H_{T}$:
\begin{itemize}
\item The unit $i\colon \mathbb C\to H_T$, $1\mapsto T^0$ and
  multiplication map $m\colon H_{T}\otimes H_{T}\to H_{T}$, given by
  $T^{m}\otimes T^{n}\mapsto T^{{m+n}}$, $m,n\in \mathbb Z$.
    
\item The coproduct $\pi\colon  H_{T}\to H_{T}\otimes H_{T}$, given by
  $T^{m}\mapsto T^{m}\otimes T^{m}$ and the counit $\epsilon\colon  H_{T}\to
  \mathbb C$, $T^{m}\mapsto 1$.
  
\item The antipode $S\colon H_{T}\to H_{T}$, $T^{m}\mapsto T^{-m}$.
\end{itemize}
These maps satisfy a large number of compatibility conditions that we
will not recall here, see e.g., \cite{MR96h:17014}.  

The powers of the translation $T$ form a basis $\{T^m\}_{m\in \mathbb
  Z}$ for $H_T$.  Another basis will also be useful for us.  Instead
of translations we can use differences as basis: let $\Delta=T-1$ and
$\bar \Delta=1-T\inv$.  The (divided) powers of $\Delta$ and $\bar
\Delta$ form the \emph{difference operator basis}
$\{\Delta[m]\}_{m\in\mathbb Z}$, where
\begin{equation}
\Delta[m]=
\begin{cases}
  \Delta^m/m! &m\ge 0,\\
  \bar\Delta^{k}/k!&m=-k<0.
\end{cases}\label{eq:3}
\end{equation}
Note that $\bar \Delta=-S\Delta$. Then we get for the antipode acting
on this basis:
\begin{equation}
S\Delta[m]=
\left .\begin{cases}
 (-\bar\Delta)^m/m! &\text{if $m\ge 0$}\\
(-\Delta)^{k}/k!&\text{if $m=-k<0$}
\end{cases}\right\}=(-1)^m \Delta[-m].
\label{eq:4}
\end{equation}
The coproduct of these basis elements becomes
\begin{equation}
  \label{eq:coproddiffbasis}
  \pi(\Delta[n])=\sum_{s=0}^n \Delta[s]\otimes
  T^s\Delta[n-s]=\sum_{s=0}^n \Delta[s]T^{n-s}
\otimes \Delta[n-s].
\end{equation}

\section{The Dual of $H_{T}$ and Sequences}\label{sec:dualseq}
$H_T$ is the group algebra (over the complex numbers) of the free
Abelian group of rank 1, i.e., of the additive group $\mathbb Z$.
Therefore the linear dual $H_T^*$ is the space of complex-valued maps
on $\mathbb Z$. So let $\mathbb C_{\mathbb Z}$ be the vector space of
two sided infinite sequences $s=(s_{n})_{n\in \mathbb Z}$, $s_{n}\in
\mathbb C$. Then the sequence $s$ defines a function on $\mathbb Z$
with value $s_n$ at $n\in\mathbb Z$. Let $\delta_n$ be the Kronecker
sequence with a 1 on the $n$th position and for the rest zeroes, so
that any sequence $s$ has an expansion 
\[
s=\sum_{n\in\mathbb Z} s_n\delta_n.
\]
Define an action of $H_T$ on $\Cz$ by $T\delta_k=\delta_{k-1}$ so that
if $s,\hat s$ are sequences with $\hat s=Ts$ then $\hat s_n=s_{n+1}$.
Define a pairing
\begin{equation}
\langle \, , \rangle\colon   H_T\otimes \mathbb C_{\mathbb Z}\to
\mathbb C,\quad P(T)\otimes s\mapsto
(P(T)s)_{n=0},\label{eq:H_TCZpairing}
\end{equation}
so that we have $\langle T^k, \delta_l\rangle=\delta_{k,l}$. Now
$\{T^m\}$ and $\{\delta_{m}\}$ are dual \lq bases\rq \ of $H_T$ and
$\Cz$ respectively (we need infinite sums of the $\delta_m$ in $\Cz$
so they don't form a basis in the algebraic sense).  We use this to
identify $\Cz$ with the dual $H_{T}^{*}$ of $H_{T}$.  The pairing
satisfies
\[
\langle T P(T),s\rangle=\langle P(T), Ts\rangle, \quad P(T)\in H_T,
s\in\Cz.
\]
The coproduct on $H_{T}$ induces via the pairing
(\ref{eq:H_TCZpairing}) a multiplication on
$\Cz$: we have
\[
\langle \pi(T^{m}),\delta_{n}\otimes\delta_{p}\rangle=\langle T^{m},
\delta_{n}\cdot \delta_p\rangle=\delta_{m,n}\delta_{m,p},
\]
so that the multiplication of basis elements of $\Cz$ is
$\delta_{n}\cdot\delta_{p}=\delta_{n}\delta_{n,p}$. The multiplication
is therefore just the usual termwise multiplication of sequences,
\[
s.\tilde s=\sum_{n\in\mathbb Z}s_n\tilde s_n\delta_n.
\]
The antipode $S$ of $H_T$ induces an antipode on $\Cz$ via
\[\langle P(T), Sf\rangle=\langle S(P(T)),f\rangle,\quad P(T)\in H_T,
f\in \Cz,
\]
so that in particular
\begin{equation}
  \label{eq:antipodeonCfin}
  S\delta_n=\delta_{-n}.
\end{equation}
The product on $H_T$ induces a map
\[
m^*\colon \Cz\to (H_T\otimes H_T)^*.
\]
In particular 
\[
m^*(\delta_n)=\sum_{k\in\mathbb Z}\delta_k\otimes \delta_{n-k}.
\]
We see that the image of $\delta_n$ under $m^*$ does not belong to
$\Cz\otimes\Cz$. In a sense $\Cz$ is too large.

\section{The Hopf Dual of $H_T$ and Difference Equations}
\label{sec:HopfDual}

Recall that the Hopf dual of a Hopf algebra is the subspace $H^\circ$
of the full linear dual $H^*$ consisting of $\phi$ such that 
\[
m^*(\phi)\in H^*\otimes H^*\subset (H\otimes H)^*.
\]
One proves that in fact in this case
\[
m^*(\phi)\in H^\circ \otimes H^\circ,
\]
and that $H^\circ$ is alternatively characterized as the subspace of
elements $\phi$ such that the $H$-submodule generated by $\phi$ is
finite dimensional. $H^\circ$ is then a Hopf algebra in its own right,
see \cite{MR83a:16010} for details. 

In case of $H=H_T$ the Hopf dual $H_T^\circ$ is the space of sequences
$s\in\Cz$ such that 
\[
H_T.s\simeq H_T/Ann_s
\]
is finite dimensional, where $Ann_s$ is the annihilator of $s$, the
ideal of elements of $H_T$ killing $s$. Since $H_T$ has Krull
dimension 1 the submodule $H_T.s$ is finite dimensional whenever
$Ann_s$ is not zero. So $H_T^\circ$ is the space of sequences that
solve at least one non trivial homogeneous constant coefficient difference
equation
\begin{equation}
P(T)s=0,\quad P(T)\ne 0\in H_T.\label{eq:HomConstantCoeff}
\end{equation}
Solutions of these equations are easy to describe explicitly. First note that
without loss of generality we can assume that 
\[
P(T)=\sum_{n=0}^Np_nT^n,\quad p_0\ne0,p_N\ne0.
\]
The space of solutions of \eqref{eq:HomConstantCoeff} is then $N$ dimensional. 
Now we can factor $P(T)=\prod_i^n(T-\lambda_i)^{d_i}$ and solving
\eqref{eq:HomConstantCoeff} reduces to solving $(T-\lambda)^ds=0$.

Define for $\lambda\ne0$ an exponential sequence
\begin{equation}
E_\lambda=\sum_{n\in\mathbb
  Z}\lambda^n\delta_n,\label{eq:defexpsequences}
\end{equation}
so that we have
\[
(T-\lambda)E_\lambda=0.
\]
For any sequence $f$ we have
\[
(T-\lambda)fE_\lambda=\lambda(\Delta F)E_\lambda,\quad \Delta=T-1,
\]
and so $fE_\lambda$ satisfies $(T-\lambda)^d fE_\lambda=0$ if and only
if $\Delta^df=0$. To find such $f$ we introduce polynomial sequences
$\tau(\ell)\in H_T^\circ,\ell\ge0$.  

First put $\tau(0)=\sum_{n\in\mathbb Z}\delta_n$, the sequence with 1
in all positions; this is the identity $\mathbf 1_{\Cz}$ in the
$\mathbb C$-algebra $\Cz$.  Then define recursively sequences
$\tau(\ell)$, as the unique solutions of the difference equations
\begin{equation}
\Delta  \tau(\ell)= \ell\tau(\ell-1),\label{eq:deftau}
\end{equation}
satisfying the initial conditions $\tau(\ell)_0=0$. We have an
explicit formula for the $\tau(\ell)$. First of all
\[
\tau(1)=\sum_{n\in\mathbb Z}n\delta_n.
\]
We abbreviate $\tau=\tau(1)$. Then
\begin{equation}\label{eq:explicitformtau}
\tau(\ell)
=\prod_{j=0}^{\ell-1}T^{-j}\tau=\sum_{n\in\mathbb Z}(n)_\ell\delta_n,
\end{equation}
where $(n)_\ell$ is the \emph{Jordan factorial}:
$(n)_\ell=n(n-1)\dots(n-\ell+1)$, for $n,\ell\in\mathbb Z,\ell\ge0$.
In particular $\tau(\ell)$ is a sequence with $0,1,\dots,\ell-1$ as
its only zeroes. The formula \eqref{eq:explicitformtau} will follow
from \eqref{eq:factortau} below.

The space 
\[
S_{\lambda,d}=\bigoplus_{\ell=0}^{d-1}\mathbb C \tau(\ell) E_\lambda
\]
is the $d$-dimensional solution space of the equation $(T-\lambda)^d
s=0$. Define
\begin{equation}
  \Czpol=\bigoplus_{\ell\ge0}\mathbb C\tau(\ell).\label{eq:defCzpol}
  \end{equation}
Then
\begin{equation}
H^\circ_T=\bigoplus_{\lambda\in\mathbb C^\times}\Czpol
E_\lambda.\label{eq:defHopfdualHT}
\end{equation}
An element $fE_\lambda$ of $H^\circ_T$ is a function on $\mathbb Z$
obtained by restricting the function $f(x)\lambda^x$, with
$f(x)\in\mathbb C[x]$, from $\mathbb C$ to $\mathbb Z$, the sequence
$\tau(1)=\tau$ corresponding to the polynomial $x$.

Note that the sequences $\tau(\ell)$ are analogs for $H_T$ (or
$\mathbb C[\Delta]$) of the powers $t^\ell\in H_D^*=\mathbb C[[t]]$,
satisfying the differential recursion
\[
Dt^\ell=\ell t^{\ell-1},\quad t^\ell\mid_0=0.
\]
Of course, the multiplicative properties of the $\tau(\ell)$'s are
more complicated. We have
\begin{equation}
\tau(\ell)T^{-\ell}\tau(m)=\tau(\ell+m),\label{eq:factortau}
\end{equation}
and hence we can calculate products in $\Czpol$:
\begin{equation}
\label{eq:prodtauellm}
  \begin{split}
    \tau(\ell)\tau(m)&=\tau(\ell)T^{-\ell}T^\ell \tau(m)
    =\tau(\ell)T^{-\ell}(1+\Delta)^\ell \tau(m)\\
    &=\sum_{k=0}^\ell {\binom{\ell}{k}} (m)_k \tau(\ell + m-k).
  \end{split}
\end{equation}
We see that the product of two basis elements of $\Czpol$ is a finite
(integral) linear combination of basis elements and $\Czpol$ is a
subalgebra of $H_T^\circ$.

To derive \eqref{eq:factortau} we use the Leibniz rule for the
difference of a product: by the coproduct formula~\eqref{eq:coproddiffbasis} 
\begin{equation}
    \label{eq:Leibniz}
    \Delta[n](fg)=\sum_{s=0}^n\Delta[s](f)T^{s}\Delta[n-s]g.
\end{equation}
We have
  \[
  \Delta(\tau T\inv\tau)=\Delta(\tau)\tau+\tau(T\inv\Delta
  \tau)=2\tau, \quad (\tau T\inv \tau)_0=0
  \]
so that by uniqueness $\tau(2)=\tau T\inv\tau$. Assuming that
(\ref{eq:factortau}) is true for all $\ell,m$ such that $\ell+m=n$
then
  \begin{align*}
    \Delta\bigg(\tau(\ell+1)T^{-\ell-1}\tau(m)\bigg)&=
    (\ell+1)\tau(\ell)T^{-\ell}\tau(m)+\\
    &\qquad\qquad\qquad+m\tau(\ell+1)T^{-\ell-1}\tau(m-1)\\
    &=(\ell+m+1)\tau(\ell+m),
  \end{align*}
by induction and so by uniqueness of solutions of (\ref{eq:deftau})
(\ref{eq:factortau}) is true for $\ell+1,m$. A similar calculation
shows the same for $\ell, m+1$.

\section{The antipodal sequences $\overline{\tau(\ell)}$}
\label{sec:Antipodaltaubar}

We defined the sequences $\tau(\ell)$ as solutions of the difference
equation \eqref{eq:deftau}, using $\Delta$. We can instead use
$\overline\Delta$ to define $\overline{\tau(\ell)}$ by
$\overline{\tau(0)}=1$ and 
\[
\overline{\Delta}\overline{\tau(\ell)}=\ell\overline{\tau(\ell-1)},\quad
\overline{\tau(\ell)}_0=0.
\]
We see that
\[
\overline{\tau(1)}=\tau(1)=\tau,
\]
and
\begin{equation}
\overline{\tau(\ell)}=\prod_{j=0}^{\ell-1}T^j\tau=\sum_{n\in\mathbb
  Z}(n)_{-\ell}\delta_n,\label{eq:explicitformtaubar}
\end{equation}
where $(n)_{-\ell}=n(n+1)\dots n+\ell-1$.
Introduce notation
\begin{equation}
  \label{eq:deftausquare}
  \tau[m]=
  \begin{cases}
    \tau(m)&\quad m\ge 0, \\
{}&\\
    \overline{\tau(k)},&\quad m=-k<0,
  \end{cases}
\end{equation}
so that
\begin{equation}
  \label{eq:tau[k]termskronecker}
  \tau[m]=\sum_{n\in\mathbb Z}(n)_m\delta_n.
\end{equation}
Using
\[
S\tau=S\sum_n n\delta_n=\sum_n n\delta_{-n}=-\tau,
\]
we find
\begin{equation}
  \label{eq:Antipodetau}
  S\tau[m]=(-1)^m \tau[-m],\quad m\in\mathbb Z.
\end{equation}
In case either $0\ge s, \ell$ or $0\le s, \ell$ we have the
factorization
\begin{equation}
  \label{eq:factorizationtau[k]}
  \tau[s+\ell]=\tau[s] T^{-s}\tau[\ell].
\end{equation}

\section{The Dual Difference Basis}
\label{sec:DualDiff}

The difference operator basis $\{\Delta[n]\}_{n\in \mathbb Z}$ of
\eqref{eq:3} has as dual the set of sequences $\Delta^*[n]$ defined by
$\langle \Delta[n],\Delta^*[m]\rangle=\delta_{nm}$. We refer to
$\{\Delta^*[n]\}_{n\in \mathbb Z}$ as the dual difference
basis. Introduce notation
\begin{equation}
  \label{eq:defthetak}
\theta(k)=
\begin{cases}
  1,\quad k\ge 0, \\
  -1,\quad k<0.
\end{cases}
\end{equation}
\begin{lem}
  For all $\ell\in\mathbb Z$
  \begin{equation}
    \label{eq:DualDiffinKronecker}
    \Delta^*[\ell]=\sum_{\theta(\ell)[n-\ell]\ge0}(n)_{\ell}\delta_n.
  \end{equation}
\end{lem}
\begin{proof}
  For any $\phi\in \Cz$ we have $\phi=\sum_n \langle
  T^n,\phi\rangle \delta_n$. Now $T=1+\Delta$, $T\inv=1-\bar\Delta$,
  so that if $n,k>0$
\[
\langle T^n, \Delta^*[-k]\rangle=0=\langle T^{-n},\Delta^*[k]\rangle.
\]
By the binomial formula we then obtain
\begin{align*}
  \Delta^*[\ell]&=\sum_{n\ge \ell}\langle
  T^n,\Delta^*[\ell]\rangle\delta_n=\sum_{n\ge 0}(n)_\ell \delta_n,\\
\intertext{and}
  \Delta^*[-\ell]&=\sum_{n\le- \ell}\langle
  T^n,\Delta^*[-\ell]\rangle\delta_n=\sum_{n\le 0}(n)_{-\ell} \delta_n.
\end{align*}
\end{proof}
Comparing \eqref{eq:DualDiffinKronecker} with
\eqref{eq:tau[k]termskronecker} we see that for $m\ne 0$ $\Delta^*[m]$
is a projection of $\tau[m]$:
\begin{equation}
\Delta^*[m]=\Pi_{\theta(m)}\tau[m]=
\begin{cases}
  \Pi_+\tau(m)& m>0, \\
  \Pi_-\overline{\tau(-m)}&m<0,
\end{cases}\label{eq:Delta*tau}
\end{equation}
where $\Pi_+$ is the restriction of functions on $\mathbb Z$ to
$\mathbb Z_{\ge0}$ and $\Pi_-$ is the restriction to $\mathbb Z_{<0}$.

We can conversely express the polynomials $\tau[k]$ in terms of the
dual difference basis. We have by calculating $\tau[k]=\sum_n \langle
\Delta[n],\tau[k]\rangle \Delta^*[k]$ in general
\begin{equation}
\tau[\pm \ell]=\Delta^*[\pm\ell] + \sum_{s<0}c_{\ell s}\Delta^*[\pm s],\label{eq:tauindualdifferencebasis}
\end{equation}
for certain integers $c_{\ell s}$, see section
\ref{sec:LocCoproductsDeltaDistributions} for an explicit formula.
In particular
\begin{equation}
  \label{eq:tauindualdiff}
  \tau[\pm 1]=\Delta^*[1]+\Delta^*[-1].
\end{equation}
\section{Localization for $\Czpol$}
\label{sec:localizationHTcirc}

Let $M\subset\Czpol$ be the multiplicative set generated by 1 and the
translates $T^k\tau$ of $\tau$, $k\in\mathbb Z$. Elements of $M$
correspond to polynomial functions on $\mathbb C$ with only zeroes at
the integers. Define 
\[
\K=M\inv\Czpol.
\]  
So an element of $\K$ is of the form $q=\frac fg$, with $f\in
\Czpol$, $g\in M$.  

\begin{remark}
  We could also consider the localization 
\[
K^\circ=M\inv H_T^\circ.
\]
This does not seem to give an essential different theory, so we will
restrict ourselves to considering $K$.
\end{remark}
Since $M$ is stable under the action of the group generated by $T$ we
can extend the action of $H_T$ on $\Czpol$ to $\K$ via
\[
T^k\frac fg= \frac{T^k f}{T^k g}, \quad f\in \Czpol,g\in M.
\]
Using the Leibniz rule (\ref{eq:Leibniz}) 
we find the quotient rule
\begin{equation}
  \Delta[\pm1](\frac fg)=\frac{\Delta[\pm1](f)g-f\Delta[\pm1](g)}{gT^{\pm1}g}.
\label{eq:Deltafrac}
\end{equation}
It follows from this that
\begin{equation}
\frac{(S\Delta)^\ell}{\ell!} \frac 1\tau=  \frac 1{  \tau(\ell+1)},
\quad
\frac{(S\overline{\Delta})^\ell}{\ell!} \frac
  1\tau=  \frac 1{\overline{ \tau(\ell+1)}}.
\label{eq:Delta1taubartau}
\end{equation}
Then we can summarize \eqref{eq:Delta1taubartau}, using the notations
\eqref{eq:defthetak}, \eqref{eq:deftausquare},
\begin{equation}
S\Delta[k]\frac1{\tau[\theta(k)]}=\frac1{\tau[\theta(k)+k]},\quad
k\in\mathbb Z.\label{eq:quotientruletau}
\end{equation}
Recall that elements of $\Czpol$ correspond to polynomial functions on
$\mathbb C$, via the algebra homomorphism
\begin{equation}
\Phi\colon \mathbb C[x]\to \Czpol, \quad f(x)\mapsto \sum_{n\in\mathbb Z}
f(n)\delta_n.\label{eq:defphi}
\end{equation}
Also recall that by the \emph{division algorithm} in $\mathbb C[x]$ we
can write, uniquely, for $f,g\in \mathbb C[x], g\ne 0$,
\[
\frac fg=q + \frac rg, \quad q, r\in\mathbb C[x],
\]
where in case $r\ne0$ we have $\deg r<\deg g$. We can and will assume
that $r$ and $g$ have no common zeroes. Furthermore we have a
\emph{partial fraction expansion} for $\frac rg$: if $g$ has zeroes
$\{n_i\}_{i\in I}$ with $n_i\in\mathbb Z$ and with multiplicities $d_i$ then
\[
\frac rg= \sum_{i\in I} \sum_{k=0}^{d_i-1}
\frac{a_{i,k}}{(x-n_i)^{k+1}},\quad a_{i,k}\in \mathbb C.
\]
Translating this result from $\mathbb C[x]$ to $\Czpol$, using the map
$\Phi$ of (\ref{eq:defphi}), we see that any $y\in \K$ can be
uniquely written as
\begin{align}
  y&=y^{\text{pol}}+y_-, \label{eq:decompK}\\
\intertext{with $y^{\text{pol}}\in \Czpol$ and} 
y_-&= \sum_{i\in I}\sum_{k=0}^{d_i-1}
  \frac{a_{i,k}}{T^{n_i}\tau^{k+1}}\label{eq:minuspartdecomp}.
\end{align}
Here the $n_i$ are integers (as $y_-$ is supposed to come, via $\Phi$,
from a fraction $\frac rg$ where $r,g\in\mathbb C[x]$ and $g$ has only
integral zeroes).  Hence
\begin{equation}
\K=\Czpol\bigoplus \Ksing, \quad \Ksing=\bigoplus_{n\in\mathbb
  Z}\bigoplus_{k\ge 0} \mathbb C
T^n\frac1{\tau^{k+1}}.\label{eq:decomK}
\end{equation}
We have then a projection
\begin{equation}
  \label{eq:sinprojK}
  \Sing\colon \K\to \Ksing.
\end{equation}
\section{$\hat H_T$ and $\hat H_T^*$}
\label{sec:hatH_Tanddual}
From \eqref{eq:decomK} it is clear that $\Ksing$ is freely generated
by $\frac1\tau$ as a module over the Hopf algebra
\[
\hat H_T=H_T[\partial_\tau]=\mathbb
C[\Delta,\bar\Delta,\partial_\tau],
\]
where $\partial_\tau$ commutes with $\Delta,\overline \Delta$, has
antipode $S(\partial_\tau)=-\partial_\tau$, coproduct
$\pi(\partial_\tau)=\partial_\tau\otimes 1+1\otimes\partial_\tau$ and
counit $\epsilon(\partial_\tau)=0$.  We have an action of $\hat H_T$
on $\Czpol$ by putting
\begin{equation}
  \label{eq:defpartialtau}
  \partial_\tau=-\log(T\inv)=\sum_{k=1}^\infty\bar\Delta^k/k,
\end{equation}
and this gives $\K$ an $\hat H_T$-module structure.

A basis for $\hat H_T$ is given by
\begin{equation}
  \label{eq:basishatH_T}
  e_{k,\ell}=\Delta[k]\partial_\tau^{(\ell)},\quad k,\ell\in\mathbb
  Z,\quad \ell\ge0.
\end{equation}
Consider the algebra $\Cz[[\sigma]]$ of formal power series in a
variable $\sigma$ with coefficients in $\Cz$, the dual of $H_T$. We define an
action of $\hat H_T$ on $\Cz[[\sigma]]$ by extending the action of
$H_T$ on $\Cz$ by putting:
\[
\partial_\tau \Delta^*[k]=0,\quad
\Delta[k]\sigma^\ell=\delta_{0,k}\sigma^\ell,\quad
\partial_\tau \sigma^{\ell}=\ell\sigma^{\ell-1}.  
\]
An element $f\in \Cz[[\sigma]]$ has an expansion
\[
f=\sum f_{n,k}\Delta^*[n]\sigma^k,
\]
and we extend the counit of $\Cz$ to a map $\epsilon:\Cz[[\sigma]]\to
\mathbb C$ by putting
\[
\epsilon(f)=f_{0,0}.
\]
We use this to define a pairing
\begin{equation}\label{eq:pairinghatHTCzsigma}
  \begin{split}
    \langle\,,\rangle\colon \hat H_T&\otimes \Cz[[\sigma]]\to\mathbb C\\
    P(\Delta,\partial_\tau)&\otimes f\quad\mapsto
    \epsilon(P(\Delta,\partial_\sigma)f).
\end{split}
\end{equation}
This identifies $\Cz[[\sigma]]$ with the dual $\hat H_T^*$, and we get
functions dual to $e_{k,\ell}$ by defining
\[
e^*_{k,\ell}  = \Delta^*[k]\sigma^\ell.
\]
\section{Action of $\Cz[[\sigma]]$ on $\hat H_T$}
\label{sec:actionCzsigmaHatHT}

The duality \eqref{eq:pairinghatHTCzsigma} allows us to define an
action
\[
\Cz[[\sigma]]\otimes \hat H_T\to \hat H_T,
\]
by putting, for $e^*\in\Cz[[\sigma]]$, $e\in\hat H_T$
\begin{equation}
e^*.e=\sum \langle e, e^*.e_{n,\ell}^*\rangle
e_{n,\ell}.\label{eq:CzsigmaactionH_T}
\end{equation}
One easily checks the following Lemma
\begin{lem}
  \label{lem:ActionCzsigmaHatHT}
  \begin{Myenumerate}
      \item 
        \[
        \sigma.e_{n,\ell}=
        \begin{cases}
          e_{n,\ell-1}& \ell>0\\
          0 & \ell=0
        \end{cases}
        \]
      \item 
        \[
        \Delta^*[1].e_{n,\ell}=
        \begin{cases}
          0& n\le0\\
          Te_{n-1,\ell} & n>0
        \end{cases}
        \]
      \item 
        \[
        \Delta^*[-1].e_{n,\ell}=
        \begin{cases}
          T\inv e_{n+1,\ell-1}& n<0\\
          0 & n\ge0
        \end{cases}
        \]
      \item 
        \begin{equation}
          \label{eq:oppositevansish}
          \Delta^*[k].\Delta[-\ell]=\Delta^*[-k].\Delta[\ell]=0,\quad k,\ell>0.
        \end{equation}
  \end{Myenumerate}
\end{lem}

\section{Action of $\Cz[[\sigma]]$ on $\Ksing$}
\label{sec:ActionCzsigmaKsing}

Let $\barKsing=\K/\Czpol$. As a vector space and $\hat H_T$-module
this is the same as $\Ksing$ but $\barKsing$ is a $\Czpol$-module,
which $\Ksing$ is not. 

Now $\Czpol$ is a subalgebra of $\Cz[[\sigma]]$, via the action of
$\hat H_T$ on $\Czpol$ (see \eqref{eq:defpartialtau}): we identify $f\in \Czpol$ with $\sum e_{n,\ell}.f|_0
e^*_{n,\ell}$. In particular 
\begin{equation}
\tau\mapsto\Delta^*[1]+\Delta^*[-1]+\sigma=e_{1,0}^*+e_{-1,0}^*+e_{0,1}.\label{eq:tauinhatHTdual}
\end{equation}
We have an $\hat H_T$-module isomorphism
\begin{equation}\label{eq:isomorphalpha}
  \begin{split}
    \alpha\colon \hat H_T &\to \barKsing\\
    P(\Delta,\partial_\tau)&\mapsto
    S(P(\Delta,\partial_\tau))\frac1\tau
\end{split}.
\end{equation}
Domain and range of $\alpha$ are $\Czpol$-modules.

\begin{lem}
  \label{lem:tauhatHTbasisaction}
  \begin{Myenumerate}

  \item \label{item:HT}In $\hat H_T$ we have
\[
\tau.e_{k,\ell}=
\begin{cases}
  T e_{k-1,\ell}+e_{k,\ell-1}&k>0, \\
  e_{0,\ell-1}&k=0, \\
  T\inv e_{k+1,\ell}+e_{k,\ell-1}&k<0 .
\end{cases}
\]

\item \label{item:barKsing}  In $\barKsing$ we have
\[
\tau.S(e_{k,\ell})\frac1\tau=
\begin{cases}
  S(T e_{k-1,\ell}+e_{k,\ell-1})\frac1\tau&k>0 \\
  S(e_{0,\ell-1})\frac1\tau&k=0 \\
  S(T\inv e_{k+1,\ell}+e_{k,\ell-1})\frac1\tau&k<0 .
\end{cases}
\]

\item \label{item:alphaiso}
Hence the map $\alpha$ of \eqref{eq:isomorphalpha} is also an
  isomorphism of $\Czpol$-modules.
\end{Myenumerate}
\end{lem}
\begin{proof}
  Part \eqref{item:HT} follows from Lemma \eqref{lem:ActionCzsigmaHatHT}
  and \eqref{eq:tauinhatHTdual}. Part \eqref{item:barKsing} follows from
  the observation that in any $\hat H_T$- and $\Czpol$-module we have
\[
 [\Delta[k],\tau]=\Delta\left[k-\theta(k)\right]T^{\theta(k)},\quad
[\partial^{(\ell)}_\tau,\tau]=\partial^{(\ell-1)}_\tau,
\]
recalling   the notation \eqref{eq:defthetak}. Then part
\eqref{item:alphaiso} follows from the other parts.
\end{proof}
We use $\alpha$ to extend the $\Czpol$ action on $\barKsing$ to an action
of $\Cz[[\sigma]]$ on $\barKsing$. In particular the product
$e_{n,k}^*.F$ is well defined for all $F\in \barKsing$, and for
instance we have, by \eqref{eq:oppositevansish},
\begin{equation}
  \label{eq:oppositeactionKsing}
  \begin{split}
    \Delta^*[k].\frac 1{\overline {\tau[\ell+1]}} &= \Delta^*[-k].\frac
    1 {\tau[\ell+1]}=0, \quad k,\ell>0.\\
  \Delta^*[\ell].
  \frac1{{\tau(m+1)}}&=\tau(\ell)\frac1{\tau(m+1)}, \quad\ell, m\ge0\\
\Delta^*[-\ell].
  \frac1{{\overline{\tau(m+1)}}}&=\overline{\tau(\ell)}
  \frac1{\overline{\tau(m+1)}},\quad \ell, m\ge0.
\end{split}
\end{equation}
\begin{remark}\label{rem:actionCzKsing}
  The action of $\Czpol$ on $\barKsing$ can be written as
  \begin{equation}
f.(\frac gh)=\frac{fg}{h}-(\frac{fg}{h})_\Hol, \quad f\in\Czpol, \frac
gh\in \Ksing.\label{eq:CzpolactionKsing}
\end{equation}
Here we use that $\K$ is the localization $M\inv\Czpol$ (with $M$ the
multiplicative set of translates of $\tau$), and that $\K$
has a decomposition in singular and holomorphic elements. Now we can
localize also the full dual $\Cz=H_T^*$: define 
\[
\KZ=M\inv \Cz.
\]
Now, maybe, one expects that similarly the $\Cz$-action on $\barKsing$
can be written in the form (\ref{eq:CzpolactionKsing}), or
equivalently that the $\Cz$-action on $\barKsing$ comes from a
projection of the natural action of $\Cz$ on $\KZ$. 
However, this is \emph{not} the case. Indeed, let $f=\delta_k\in\Cz$,
the Kronecker sequence with support at $k>0$, for example. Then we
have 
\[
\delta_k\tau(k+1)=0, \quad\text{see (\ref{eq:explicitformtau})},
\]
so that for all $\frac gh\in \K$ (or $\K_{\mathbb Z}$) we have
\[
\frac{\delta_kg}{h}=\frac{\delta_k\tau(k+1)g}{\tau(k+1)h}=0\in\KZ.
\]
On the other hand we have by (\ref{eq:CzsigmaactionH_T})
\[
\delta_k.\frac1{\tau-\ell}=\delta_{k,ell}\frac1{\tau-\ell}\ne 0.
\]
In fact, inside $\KZ$ there is no distinction between
holomorphic and singular elements. To see this note that we have an
exact sequence
\[
0\to \Cztor\to\Cz\overset\beta\to \KZ
\]
where $\beta\colon \Cz\to \KZ$ is the canonical map
$f\mapsto \frac f1$, and 
\[
\Cztor=\oplus_{k\in\mathbb Z}\mathbb C\delta_k.
\]
Define $\Cztf=\Cz/\Cztor$. Then we see that also
\[
\KZ=M\inv\Cztf,
\]
where we identify $M$ with its image in $\Cztf$, as
$M\cap\Czpol=\emptyset$.  Now if $m\in M$ then $m$ has at most a
finite number of zeroes (as function on the integers), so modulo
$\Cztor$ $m$ is invertible. Hence as subset of $\Cztf$ $M$ consists of
units and so we find $\KZ\simeq \Cztf$. There are no singular elements
in $\Cztf$.
\end{remark}

\section{Trace for $\K$}
\label{sec:TraceK}

Now $\hat H_T$, being a Hopf algebra, has a counit, a multiplicative
map $\epsilon\colon \hat H_T\to\mathbb C$ given on the basis
\eqref{eq:basishatH_T} by 
\[
\epsilon(e_{k,\ell})=\epsilon(\Delta[k]\partial_\tau^{(\ell)})=\delta_{k,0}\delta_{\ell,0}.
\]
The counit induces, via the isomorphism $\alpha$, a map
$\barKsing\to\mathbb C$, which we extend (by zero) to all of
$\K$ to get a map, called the \emph{Trace}
\begin{equation}
  \label{eq:deftraceK}
  \Tr\colon \K\to\mathbb C,\quad \frac fg\mapsto
  \epsilon(\alpha\inv(\Sing(\frac fg))),
\end{equation}
where $\Sing$ is the projection \eqref{eq:sinprojK}
\begin{lem}
  \begin{Myenumerate}
  \item If $F\in \K$, $h\in \hat H_T$, then
\[ \Tr(h F)=\epsilon(h)\Tr(F).
\]
\item In particular we have
\[
\Tr(T^k F)=\Tr(F),\quad \Tr(\Delta F)=\Tr(\overline\Delta
F)=\Tr(\partial_\tau F)=0.
\]
\item (Partial Summation/Integration) For all $F, G\in \K$ and
  $h\in\hat H_T$ we have
  \begin{equation}
    \label{eq:PartialSum}
\Tr(h\langle F\rangle G)=\Tr(F S(h)\langle G\rangle).    
  \end{equation}
  \end{Myenumerate}
\end{lem}

Recall that an element $F\in\K$ corresponds to a rational function $F(x)$ on
$\mathbb C$ with at most a finite number of finite order poles at the
integers. Let $F(x)dx$ be the corresponding rational 1-form on
$\mathbb C$. Then the trace of $F$ is the sum of residues of this 1-form.

\begin{lem}
  If $F\in \K$ then
\begin{equation}
\Tr(F)= \sum_{n\in\mathbb Z}\Res_{x=n}
(F(x)dx).\label{eq:TraceResidue}
\end{equation}
\end{lem}

Note that if $f\in \Czpol$ the value of $f$ at 0 can be expressed
as a trace
\begin{equation}
  \label{eq:tracevaluezero}
  f|_0=\Tr(\frac f\tau).
\end{equation}
The terminology Trace for the map (\ref{eq:deftraceK}) is used to
indicate the analogy with the trace of finite matrices. The
Hamiltonians of the finite Toda lattice are traces of Lax
matrices. Similarly the Hamiltonians for the infinite Toda lattice
are traces in the sense of (\ref{eq:deftraceK}). This is
discussed in more detail at the end of section \ref{sec:LaxOp}.

\begin{remark}
  The trace on $\K$ can not reasonably extended to a trace on the
  localization $\KZ$ introduced in Remark~\ref{rem:actionCzKsing}. As
  explained there, in $\KZ$ there is no distinction between singular
  and holomorphic elements. 
  
  However, what we can and will do is define a pairing, by abuse of
  notation also called trace, between $\Cz$ and $\Ksing$:
  \[
  \Cz\otimes \Ksing\to \mathbb C,\quad f\otimes \frac gh\mapsto
  \Tr(f.(\frac gh))
  \]
  We warn the reader that in general, for $\frac gh\in \Ksing$
  \[
  \Tr(f.(\frac gh))\ne\Tr( (mf).(\frac g{mh}) ),\quad m\in M,
  \]
  unless $f\in\Czpol$, see Lemma \ref{lem:tauhatHTbasisaction}, part
  ~\ref{item:alphaiso}.

\end{remark}
\section{Orthogonality Relations}
\label{sec:orth-relat}

We introduce notation:  for $\ell\ge0$ we put
\begin{equation}
\tau(-\ell-1)=\frac1{\tau(\ell+1)},
\quad\overline{\tau(-\ell-1)}=\frac1{\overline{\tau(\ell+1)}}.\label{eq:deftau-n-1}
\end{equation}
\begin{lem}\label{lem:orthtautau} For all $m,n\in\mathbb Z$
  \[
\Tr\left(\tau(m)\tau(n)\right)=\delta_{m,-n-1},\quad 
\Tr\left(\overline{\tau(m)}\overline{\tau(n)}\right)=\delta_{m,-n-1}.
\]
\end{lem}
\begin{proof} We consider the left hand identity first. 
  In case $m,n\ge 0$ the trace is zero because $\tau(m)\tau(n)$ is
  nonsingular. In case $m,n\le -1$ the one form $\frac1 f dx$
  corresponding to $\frac1{\tau(-m)\tau(-n)}$ is regular at $\infty$
  and has only poles at the integers, so we can, using
  \eqref{eq:TraceResidue}, replace the trace by the sum over the
  residues of $\frac1 f dx$ at all points of $\mathbb P^1$, which is
  zero by the residue theorem.  So assume $m\ge0$ and $n=-k-1$ for
  $k\ge0$. Then we get using the quotient rule
  (\ref{eq:Delta1taubartau}), the partial summation rule
  \eqref{eq:PartialSum}, \eqref{eq:deftau}, and
  \eqref{eq:factortau} and \eqref{eq:tracevaluezero}
\[
   \Tr\left(\tau(m){\tau(n)}\right)=
    \Tr\left(\tau(m)\frac{(-\overline{\Delta}^{\,k})}{k!}\frac
      1{\tau}\right)=\Tr\left(\frac{\Delta^k}{k!}(\tau(m))\frac
      1\tau\right)=\delta_{m,k}.
\]
The rest of the proof is similar.
\end{proof}

Recall the action of $\Cz[[\sigma]]$ on $\barKsing$ described by
combining Lemma \ref{lem:ActionCzsigmaHatHT} with the isomorphism $\alpha$. 
\begin{lem}\label{lemma:orthogonalityhatHT}
  For all $n,m,\ell, k\in\mathbb Z$, $\ell,k\ge 0$ we have
\[
\Tr(e_{n,\ell}^*.S(e_{m,k})\frac1\tau)=\delta_{n,m}\delta_{\ell,k}.
\]
\end{lem}

 \section{Expansions in $\K$}
\label{sec:ExpansionsK}

A basis for $\K$ is given by
\begin{equation}
\{\tau(k),k\ge0\}\cup\{S(e_{n,\ell})\frac1\tau,n,\ell\in\mathbb
Z,\ell\ge0\}.\label{eq:basisK}
\end{equation}
So any $F\in K$ has a (finite) expansion
\begin{equation}
  \label{eq:ExpansionK}
  F=\sum F_{\{ \frac1{\tau(k+1)}\}} \tau(k)+ F_{\{e^*_{n,\ell}\}} S(e_{n,\ell})\frac1\tau.
\end{equation}
The coefficients of $F$ are traces:
\[
 F_{\{ \frac1{\tau(k+1)}\}}=\Tr(\frac1{\tau(k+1)}F),\quad
 F_{\{e^*_{n,\ell}\}}=\Tr(e^*_{n,\ell}.F),
\]
using the action of $\Cz[[\sigma]]$ on $\Ksing$, see section
\ref{sec:ActionCzsigmaKsing}.

\section{(Twisted) Coproduct and (Twisted) Exponential Operators}
\label{sec:exp_op} 

Let $H$ be a commutative and cocommutative Hopf algebra, for
simplicity, and let $H^*$ be its dual. Fix a basis $\{e_i\}_{i\in I}$
for $H$ and let $e_j^*$ be linear functions on $H$ such that
\[
e_j^*(e_i)=\delta_{ij}.
\]
Any $\phi\in H^*$ can then be expanded as
\begin{equation}
\phi=\sum \phi(e_i)e_i^*,
\label{eq:expamsionalphaHTdual}
\end{equation}
and similarly we have for $\omega\in(H\otimes H)^*$
\[
\omega=\sum \omega(e_i\otimes e_j)e_i^*\otimes e_j^*.
\]
Let now $A,B$ be two linear maps $H\to H$. (Later on $A,B$ will be
either the identity or the antipode $S$.) Define, if $m$ is the usual
multiplication of $H$, the $A\otimes B$-twisted multiplication on $H$
by
\[
m_{A\otimes B}\colon H\otimes H\to H,\quad x\otimes y\mapsto m(Ax\otimes
By).
\]
Dual to the twisted multiplication we get a map
\[
m_{A\otimes B}^*\colon H^*\to (H\otimes H)^*.
\]

\begin{lem}\label{lem:ABtwistedcoproductExp}
  We have for all $\phi\in H^*$
\[
m_{A\otimes B}^*(\phi)=\mathcal L_{A\otimes B}(\phi)=\mathcal
R_{A\otimes B}(\phi),
\]
where
\begin{align*}
  \mathcal L_{A\otimes B}(\phi)&=\sum_{i\in I}e_i^*\otimes
  B^*\left((Ae_i).\phi\right),\\
  \mathcal R_{A\otimes B}(\phi)&=\sum_{i\in I}A^*((Be_i).\phi)\otimes e_i^*.
\end{align*}
\end{lem}

\begin{proof}
  If, for instance, $B:H\to H$ is a linear map, then the dual is given
  by $B^*\phi(e)=\phi(Be)$ for $\phi\in H^*$, $e\in H$. Then we have
\begin{align*}
  m_{A\otimes B}^*(\phi)&= \sum m_{A\otimes B}^*(\phi)(e_i\otimes e_j) e_i^*\otimes e_j^* = 
  \sum \phi(A(e_i)B( e_j)) e_i^*\otimes e_j^*\\
  &= \sum (Ae_i).\phi (B e_j) e_i^*\otimes e_j^*=\mathcal L_{A\otimes B}(\phi)
\end{align*}
By a similar calculation we find the expression of the $A\otimes B$
twisted coproduct in terms of $\mathcal R_{A\otimes B}$.
\end{proof}
Because of Example \ref{ex:HDexponentials} below we call $\mathcal
L_{A\otimes B},\mathcal R_{A\otimes B}$ 
twisted exponential operators. In case $A=B=1$ we write 
$\mathcal L_{1\otimes 1}=\mathcal L, \mathcal
R_{1\otimes1}=\mathcal R$
and refer to $\mathcal L, \mathcal R$ as untwisted exponential operators.
We will mainly interested in the case $A=1,B=S$, and we will write for
brevity
\[
\mathcal{L}^S=\mathcal{L}_{1\otimes S},\quad \mathcal{R}^S=\mathcal{R}_{1\otimes S},
\]

\begin{example}\label{ex:HDexponentials}
  Take $H=H_D$, with as basis the divided powers $e_i=D^{(i)}=D^i/i!$.
  The dual is $H_D^*=\mathbb C[[t]]$, with $e_i^*=t^i$ . The action of
  $H_D$ on $H^*_D$ is given by $D=\partial_t$. If $f(t)\in\mathbb
  C[[t]]$ then the untwisted comultiplication is $m^*(f)=f(t_1+t_2)$.
  Here and below we will write $t_1, t_2$ for $t\otimes1$, $1\otimes
  t$.  We have two expansions:
  \[
  f(t_1+t_2)=\exp\left(t_2\partial_{t_1}\right)f(t_1)=\exp\left(t_1\partial_{t_2}\right)f(t_2).
  \]
  In this case the untwisted exponential operators are
  \[
\mathcal L=\sum_i e_i^*\otimes e_i=\exp(t_1\partial_{t_2})  , \quad
\mathcal R=\sum_i e_i\otimes e_i^*=\exp(t_2\partial_{t_1}). 
  \]
  The $1\otimes S$-twisted comultiplication is
\[
m_{1\otimes S}^*(f(t))=f(t_1-t_2),
\]
and
\[
\mathcal L^{S}=S_2\exp(t_1\partial_2),\quad \mathcal
R^{S}=\exp(-t_2\partial_1).
\]
\qed
\end{example}

\begin{example}\label{ex:HTexponentials}
  Now take $H=H_T$, with as basis the difference operators
  $\{e_n=\Delta[n]\}$ of \eqref{eq:3} and as dual basis the sequences
  $\{e_n^*=\Delta^*[n]\}$ introduced in section
  \ref{sec:DualDiff}. If $\phi\in \Cz$ then $\phi$ is a
  function on the integers. The coproduct $m^*(\phi)$ and the
  twisted coproduct $\mStwod(\phi)$ are functions on the product
  $\mathbb Z\times\mathbb Z$ given by
  \begin{equation}
m^*(\phi)(m,n)=\phi(m+n),\quad
\mStwod(\phi)(m,n)=\phi(m-n).\label{eq:cproductsasfunction}
\end{equation}
The untwisted and twisted exponential operators for $H_T$ using the
difference operator basis are
\begin{equation}
\begin{split} 
\mathcal L&=\sum_{n\in\mathbb Z}\Delta^*[n]\otimes \Delta[n],\quad
\mathcal L^{S}(-)=\sum_{n\in\mathbb
  Z}\Delta^*[n]\otimes S\left(\Delta[n](-)\right),
\\
\mathcal R&=\sum_{n\in\mathbb Z}\Delta[n]\otimes
\Delta^*[n],\quad
\mathcal R^{S}(-)=\sum_{n\in\mathbb
 Z}S\left(\Delta[n]\right)(-)\otimes \Delta^*[n],
\end{split}\label{eq:HTexponentials}
\end{equation}
so that $m^*(\phi)=\mathcal L(\phi)=\mathcal R(\phi)$ and
$m_{1\otimes S}^*(\phi)=\mathcal L^{S}(\phi)=\mathcal
R^{S}(\phi)$.\qed
\end{example}

\section{$H$-Covariance of (Twisted) Exponentials}
\label{sec:HinvarianceTwistedExp}

The domain $H^*$ and the range $(H\otimes H)^*$ of the twisted
comultiplications admit actions of $H$ and $H\otimes H$, and one can
expect that for suitable $A,B$ the twisted exponentials relate these
actions.

We will write $h_1=h\otimes 1$ and $h_2=1\otimes h$ and similar for
other operators acting on tensor products.

\begin{lem}\label{lem:HpropertiesExponentials}
  For all $\phi\in H^*$, $h\in H$
  \begin{Myenumerate}
  \item \label{Part:1}
$\mathcal L(h\phi)=(h)_1\mathcal L(\phi)=h_2 \mathcal L(\phi)$.
\item \label{part:2}
$\mathcal R(h\phi)=(h)_1\mathcal R(\phi)=h_2 \mathcal R(\phi)$.
  \item \label{Part:3}
$\mathcal L^{S}(h \phi)=h_1\mathcal L^{S}(\phi)=(Sh)_2\mathcal L^{S}(\phi)$.
\item \label{Part:4}
$\mathcal R^{S}(h \phi)=h_1\mathcal R^{S}(\phi)=(Sh)_2\mathcal R^{S}(\phi)$.
  \end{Myenumerate}
\end{lem}
\begin{proof}
We have by the commutativity of $H$
\[
  \mathcal L(h\phi)=\sum e_i^*\otimes e_i h\phi=1\otimes h \sum
  e_i^*\otimes e_i\phi=h_2\mathcal L(\phi).
\]
On the other hand, introducing constants $h^j_i\in\mathbb C$
describing the action of $h\in H$ on $H^*$ and $H$:
\[
 h e_i=\sum_j h^j_i e_j, \quad h e_j^*=\sum_i h^j_i e_i^*,
\]
we see that
\[
  \mathcal L(h\phi)=\sum_i e_i^*\otimes  e_ih\phi
=\sum_i e_i^*\otimes \sum_j h^j_ie_j\phi=\sum_jhe_j^*\otimes e_j\phi =h_1 \mathcal L(\phi).
\]
The proof of part \eqref{part:2} is the same. 

Similarly we find
\[
  \mathcal L^{S}(h \phi)=\sum_i e_i^*\otimes S(e_ih\phi)
  =1\otimes S(h)\sum_i e_i^*\otimes  S(e_i\phi)=(S h)_2\mathcal L^{S}(\phi).
\]
On the other hand
\[
  \mathcal L^{S}(h \phi)= \sum_{i,j} e_i^*\otimes S(h_i^j
  e_j\phi)=  \sum_j he_j^*\otimes S(e_j\phi)=(h)_1\mathcal L^{S}(\phi).
\]
Part \eqref{Part:4} is proved similarly.
\end{proof}

\begin{remark}\label{rem:ExponentialsMHcovariance}
  Let $M$ be an $H$-module. Then we can define for $m\in M$ the action of,
  say, twisted exponentials by
\begin{equation}
\begin{split} 
\mathcal L_M^{S}(m)&=\sum e^*_i\otimes S_M\left(e_im\right),
\\
\mathcal R_M^{S}(m)&=\sum S_{H}\left(e_i\right)(m)\otimes e_i^*,
\end{split}\label{eq:HTexponentialsM}
\end{equation}
where $S_M:M\to M$ is supposed to be a linear map such that
\[
S_M(hm)=S_H(h)S_M(m),\quad h\in H_T,\quad m\in M.
\]
In this situation Lemma \ref{lem:HpropertiesExponentials} will still
be true, with the same proof. We will frequently omit the subscript
${}_M$ on the exponentials in case the module they are acting on
should be clear.

\end{remark}
\section{Multiplicativity of (Twisted) Exponentials}
\label{sec:MultTwistedExp}

The domain $H^*$ and range $(H\otimes H)^*$ of twisted exponentials
are algebras so one can expect, for suitable twistings $A,B$, that
these algebra structures are related. We will need the following cases.

\begin{lem}\label{lem:multiplicativityexps}
  The exponentials $\mathcal L, \mathcal L^{S}, \mathcal R, \mathcal
  R^{S}$ are multiplicative: we have for all $\phi,\beta\in H^*$
\begin{Myenumerate}
\item $\mathcal L(\phi\beta)=\mathcal L(\phi)\mathcal L(\beta)$,
  $\mathcal R(\phi\beta)=\mathcal R(\phi)\mathcal R(\beta)$.
  \item $\mathcal L^{S}(\phi\beta)=\mathcal L^{S}(\phi)\mathcal
    L^{S}(\beta)$, $\mathcal R^{S}(\phi\beta)=\mathcal R^{S}(\phi)\mathcal R^{S}(\beta)$.
  \end{Myenumerate}
\end{lem}

\begin{proof}
  Introduce structure constants $C^k_{ij}\in\mathbb C$ for the
  multiplication in $H^*$ and the comultiplication in $H$:
\[
e_i^*e_j^*=\sum_k C^k_{ij}e_k^*,\quad \pi(e_k)=\sum_{i,j}
C^k_{ij}e_i\otimes e_j.
\]
Then, for instance,
\begin{align*}
  \mathcal L^{S}(\phi\beta)&=\sum_k e_k^*\otimes
  S(e_k(\phi\beta))=\sum_{k,i,j} e_k^*\otimes
  S(C^k_{ij}e_i(\phi)e_j(\beta))\\
  &=\sum_{k,i,j} C^k_{ij}e_k^*\otimes S(e_i(\phi)e_j(\beta))
  =\sum_{i,j} e_i^*e_j^*\otimes S(e_i(\phi))S(e_j(\beta))\\
  &=\mathcal L^{S}(\phi)\mathcal L^{S}(\beta).
\end{align*}
Similar computations gives the multiplicativity of the other exponentials.
\end{proof}

\section{$H_T$-Leibniz Algebras and Exponentials}
\label{sec:HTLeibnizExpo}

If $H$ is a Hopf algebra, then an \emph{$H$-Leibniz algebra} is a
associative algebra in the category of $H$-modules.  So $V$ is both an
$H$-module and an algebra, such that the unital algebra structure maps
\begin{align*}
i\colon \mathbb C  &\to V                &   m\colon V\otimes V &\to V    \\
 \lambda&\mapsto \lambda\mathbf{1}_V       &       a\otimes b &\mapsto ab
\end{align*}
are $H$-module morphisms (giving $\mathbb C$ the trivial
$H$-module structure): if $\epsilon$ is the counit and $\pi(h)=\sum h^\prime\otimes
h^{\prime\prime}$ is the coproduct, then
\begin{equation}\label{eq:HLeibnizalgebra}
  hi(\lambda)=i(\epsilon(h)\lambda),\quad hm(a\otimes b)=\sum
  m(h^\prime a\otimes h^{\prime\prime}b).
\end{equation}
In case $H=H_D=\mathbb C[D]$ a (commutative) $H_D$-Leibniz algebra is
just a (commutative) differential algebra, and the last equation of
\eqref{eq:HLeibnizalgebra} is the Leibniz rule of differentiation of
a product. For any cocommutative Hopf algebra the dual $H^*$ is a
commutative $H$-Leibniz algebra. We will assume cocommutativity from
now on, being interested in the case $H=H_T$.

Examples of $H_T$-Leibniz algebras $V$ are $\Cz$, $\Czpol$,
$H_T^\circ$ and the localizations $\K$, $\K^\circ$.
\begin{remark}
We noted in Remark \ref{rem:ExponentialsMHcovariance}
that (twisted) exponentials act on any $H$-module and that they will
have then the $H$-covariance properties of Lemma
\ref{lem:HpropertiesExponentials}. In case $V$ is a
$H_T$-Leibnitz algebra (twisted) exponentials will act on $V$, and
in this situation furthermore the multiplicativity properties of
Lemma~\ref{lem:multiplicativityexps} will still hold, with the same
proof.
\end{remark}

\section{Inverses of exponentials.}
\label{sec:InverseExp}
The twisted exponential operators $\mathcal{L}^{S},
\mathcal{R}^{S}$ are, in a sense, inverses of the untwisted
ones. Indeed, recall that for $\phi\in\Cz$, i.e., if $\phi$ is a
function on the integers, then the copoducts $\md(\phi)$ and $\md
_{1\otimes S}$ are functions on $\mathbb Z\times \mathbb Z$ given by
\eqref{eq:cproductsasfunction}. We can define functions on $\mathbb
Z\times \mathbb Z\times \mathbb Z$ by
\[
\md(\md_{1\otimes S}(\phi))(p,q,r)=\phi(p-q+r),
\]
and restricting to the diagonal $q=r$ we get $\phi$ back. This holds
not only in $\Cz$ but is an identity for exponential operators acting
on any $H_T$-module. More precisely, we have the following Lemma.

\begin{lem}\label{lem:InverseExp}
  Let $V$ be an $H_T$-module. Then
  \begin{align*}
    \mathcal R\circ\mathcal R^{S}&=\mathcal R_{1\otimes
      S}\circ \mathcal R=1_V\otimes 1_{H^*_T},\\
    \mathcal L\circ\mathcal L^{S}&=\mathcal L_{1\otimes
      S}\circ \mathcal L=1_{H^*_T}\otimes 1_V.
  \end{align*}
\end{lem}
 \begin{proof}
   For instance
   \begin{align*}
     \mathcal R^{S}\circ\mathcal R(v)&
     =\sum_{i,j}S(e_i)e_j.(v)\otimes e_j^*e_i^*=
     \sum_{i,j,k}\gamma_{ij}^{S,k} e_k(v)\otimes e_j^*e_i^*.
 \intertext{Now $ \gamma_{ij}^{S,k}=\langle
       S(e_i)e_j,e_k^*\rangle=\langle S(e_i)\otimes
       e_j,m^*(e_k^*)\rangle$, where $m^*(e^*_k)=\sum
       e_k^{*\prime}\otimes e_k^{*\prime\prime}$ is the coproduct on
       $H^\circ$ dual to the product on $H$. Hence}
     \mathcal R^{S}\circ\mathcal R(v)    &=\sum_{i,j,k}e_k(v)\otimes
     \gamma_{ij}^{S,k}e_j^*e_i^*
=\sum_{k}e_k(v)\otimes \mu_{S\otimes 1}(m^*(e_k)),
 \intertext{where $\mu_{S\otimes 1}$ is the twisted multiplication
       on $H^\circ$, so}
     \mathcal R^{S}\circ\mathcal R(v)    &=\sum_{k}e_k(v)
       \otimes \sum S(e_k^{*\prime})e_k^{*\prime\prime}=\sum_{k}e_k(v)\otimes i^*(e_k^*)I_{H^*}, 
 \intertext{using the connection $i^*(\phi)1=\sum
   S(\phi^\prime)\phi^{\prime\prime}$ between counit, antipode
       and coproduct on a Hopf algebra, hence}
     \mathcal R^{S}\circ\mathcal R(v)    &=\sum_{k}i^*(e_k^*)e_k(v)\otimes I_{H^*}=v\otimes 1_{H^*},
 \end{align*}
 since $i^*(e_k^*)e_k=1_H$. Similarly one proves the other parts.
 \end{proof}
\section{Adjoint action of exponentials}
\label{sec:AdjointExp}

Define for $h\in H$ and $X\in \operatorname{End}(V)$, for some
$H$-module $V$, the adjoint action of $h$ by
\[
\ad_h(X)=\sum h^\prime\circ X\circ S(h^{\prime\prime}).
\]
The twisted coproduct $\pi_{1\otimes S}\colon h\mapsto h^\prime\otimes
S(h^{\prime\prime})$ is dual to the twisted multiplication $\mu_{1\otimes
S}$ on $H^*$. 
Extend the adjoint action termwise to exponential operators:
\[
\ad_{\mathcal R}(X)=\sum \ad_{e_i}(X)\otimes e_i^*,\quad \ad_{\mathcal
  R^{S}}(X)=\sum \ad_{Se_i}(X)\otimes e_i^*.
\]

\begin{lem}\label{lem:Adjointactionexponentials}
  For all $X\in \End(V)$ we have
\[
\ad_{\mathcal R}(X)=\mathcal R\circ X\circ \mathcal R^{S},
\quad 
\ad_{\mathcal R^{S}}(X)=\mathcal R^{S}\circ X\circ \mathcal R.
\]
\end{lem}

\begin{proof}
We have
\[
\pi_{1\otimes S}(\mathcal R)=\sum_i \pi_{1\otimes S} (e_i)\otimes e_i^*=\sum
\gamma^i_{jk} e_j\otimes e_k\otimes e_i^*,
\]
where 
\[
\gamma^i_{jk}=\langle \pi_{1\otimes S} (e_i), e_j^*\otimes
e_k^*\rangle= \langle e_i, e_j^*S(e_k^*)\rangle.
\]
So
\[
\pi_{1\otimes S}(\mathcal R)= \sum e_j\otimes e_k\otimes
e_j^*S(e_k^*)=\mathcal R\dot\otimes \mathcal R^{S},
\]
where we use
\begin{equation}
  \label{eq:twoformsR1otimesS}
  \mathcal{R}^{S}(f)=\sum S(e_i)(f)\otimes e_i^*=\sum
  e_j(f)\otimes S(e_j^*).
\end{equation}
Hence
\[
\ad_{\mathcal R }(X)=\mathcal R\circ X\circ\mathcal R^{S}.
\]
The other part is proved similarly.
\end{proof}

\begin{lem}
  If $V$ is an $H$-module, the map $H\to \End(\End(V))$,
\[
h\mapsto \ad_h: \End(V)\to \End(V),\quad X\mapsto \sum h^\prime \circ
X\circ S(h^{\prime\prime})
\]
gives $\End(V)$ the structure of $H$-Leibniz algebra. In particular
\[
\ad_h(X\circ Y)=\sum \ad_{h^\prime}(X)\circ\ad_{h^{\prime\prime}}(Y).
\]
\end{lem}
\begin{cor}
  \[
  \ad_{\mathcal{R}}(X\circ Y)=\ad_{\mathcal{R}}(X)\circ
  \ad_{\mathcal{R}}(Y).
\]
\end{cor}

\section{Multivariable expansions}
\label{sec:MultivariableExpansions}

Let $H$ be a Hopf algebra, commutative and cocommutative for
simplicity.  If $M$ is an $H$-module, and $m\in M$, then we get a
distribution $\mathcal{R}_M(m)$ on $H$ and the twisted version
$\mathcal{R}^S_M(m)$. We think of these distributions as expandable in a
basis $e_i^*$ of $H^*$, and to emphasize this we write
$\mathcal{R}_M(e^*)(m)$, $\mathcal{R}_M^S(e^*)(m)$. We will also consider
distributions on $H\otimes H$ and $H\otimes H\otimes H$.  For instance
if $\mathcal{D}$ is a distribution on $H^{\otimes 3}$ it will be
expandable in $e_{i,1}^*=e_{i}^*\otimes 1\otimes 1$, $e_{i,2}=1\otimes
e_i^*\otimes 1$, $e_{i,3}=1\otimes 1\otimes e_i^*$, and we write
$\mathcal{D}(e^*_1,e^*_2,e^*_3)$. In this section we discuss how the
invertibility properties of exponential operators, see section
\ref{sec:InverseExp}, leads to identities between distributions in
several variables. 

First we need some notation. If $e^*\in H^*$ we
write $e^*_{12}=\mStwod(e^*)=\sum
e^{*\prime}\otimes e^{*\prime\prime}$, $e_{13}^*=\sum
e^{*\prime}\otimes 1\otimes e^{*\prime\prime}$ and 
$e_{23}^*=1\otimes e^{*,\prime}\otimes e^{*,\prime\prime}$.
We can consider exponential operators in these variables, for instance
\[
\mathcal{R}^S_M\left(e^*_{23}\right)m=\sum S(e_i)_M m \otimes e^*_{i, 23}.
\]
\begin{lem}\label{lem:productExpandcoproduct}
  For all $m\in M$ we have
\[
\mathcal{R}_M(e_1^*)\mathcal{R}^S_M(e_2^*)(m)=\mathcal{R}_M(e^*_{12})m.
\]
\end{lem}

\begin{proof}
  We have
  \begin{align*}
    \mathcal{R}_M(e^*_1)\mathcal{R}^S_M(e^*_2)(m)
    &=\sum e_iS(e_j)(m)\otimes e_i^*\otimes e_j^*\\
    &=\sum c_{ij}^k e_k(m) \otimes e_i^*\otimes e_j^*
\intertext{where $\sum e_iS(e_j)=\sum c_{ij}^k e_k$. We also have
  $\mStwod(e_k^*)=\sum c^k_{ij}e_i^*\otimes e^*_j$, so that}
    &=\sum  e_k(m) \otimes \mStwod(e_k^*)=\mathcal{R}_M(e^*_{12})(m).
  \end{align*}
\end{proof}

Then, by section \ref{sec:InverseExp} and Lemma \ref{lem:productExpandcoproduct},
we get identities like
\[
\mathcal{R}_M(e^*_1)(m)=\mathcal{R}_M(e^*_2)\mathcal{R}^S_M(e^*_2)\mathcal{R}_M(e^*_1)(m)=
\mathcal{R}_M(e^*_2)\mathcal{R}_M(e^*_{12})(m).
\]
We now specialize to the case  $H=\mathbb C[\Delta]$ and $M=\Czpol$.
\begin{lem}\label{lem:multivariabletauexp}
 Then for all $\ell\ge0$ we have
\[
\tau_{12}[\ell]=\sum_{s=0}^\ell \binom{\ell}{s}
\tau_{13}[\ell-s]\tau_{32}[s].
\]
Here $\tau_{32}[s]=\mathcal{R}^S(\tau_2)\tau_3[s]$.
\end{lem}
\begin{proof}
 We have
 \begin{align*}
   \tau_{12}[\ell]&=\mathcal{R}^S(\tau_2)(\tau_1[\ell])&&\\
   &=\mathcal{R}^S(\tau_3)\mathcal{R}(\tau_3)\mathcal{R}^S(\tau_2)(\tau_1[\ell])&&\text{Lemma
     \ref{lem:InverseExp}}\\
   &=\mathcal{R}^S(\tau_{23})\mathcal{R}^S(\tau_3)(\tau_1[\ell])&&\text{Lemma
     \ref{lem:productExpandcoproduct}}
   \\
   &= \mathcal{R}^S(\tau_{23})(\tau_{13}[\ell])=\sum_{s=0}^\ell
   \binom{\ell}{s} \tau_{13}[\ell-s]\tau_{32}[s].&&
 \end{align*}
\end{proof}

\section{Extension of Exponential Distributions}
\label{sec:extExpDistr}

Let $W, H$ be vector spaces. A \emph{$W$-valued distribution on $H$} is just
a linear map $H\to W$. In this section we describe the following
phenomenon.  Assume that $H$ is an algebra and let $H\subset H^\prime$
be an inclusion of $H$ in a bigger algebra $H^\prime$. Let $\{e_i\}$
be a basis for $H$, with $e_i^*\in H^*$ dual functions. Let $M$ be an
$H^\prime$-module and consider for $m\in M$ the action by the
exponential for $H$ on $m$:
\begin{equation}
\mathcal{R}^H(m)=\sum e_i(m)e_i^*.\label{eq:ExpactionM}
\end{equation}
So $\mathcal{R}^H(M): H\to M$ is the $M$-valued distribution given by
\[
\mathcal{R}^H(m)(h)=h.m,\quad h\in H.
\]
Now it can happen that the infinite series \eqref{eq:ExpactionM} makes
sense, for all $m\in M$, as a distribution on the bigger algebra
$H^\prime$, in such a way that
\[
\mathcal{R}^H(m)(h^\prime)=h^\prime.m,\quad h^\prime\in H^\prime.
\]
In this case we say that the distributions $\mathcal{R}^H(m)$ (on $H$)
can be extended to $H^\prime$.  This means that we can replace on $M$
the exponential operator $\mathcal{R}^{H^\prime}$ by the
simpler $\mathcal{R}^H$.

The basic example we are interested in is where $H=\mathbb C[\Delta]$
and $H^\prime$ is either $H_T$ or $\hat H_T$. Here $\mathbb C[\Delta]$
is the semigroup algebra of the semigroup generated by $T$, with
basis $\Delta[\ell]$ and dual functions $\tau[\ell]$, $\ell\ge0$. So
the exponential operators of $\mathbb C[\Delta]$ on $M$ are defined by
\begin{equation}
  \begin{aligned}
  \mathcal L^{\Delta}(f)&=\sum_{\ell\ge0}\tau[\ell]\otimes
  \Delta[\ell](f), &  \mathcal R^{\Delta}(f)&=\sum_{\ell\ge0}\Delta[\ell](f)\otimes
  \tau[\ell],
\\
  \mathcal L^{\Delta}_{1\otimes S}(f)&=\sum_{\ell\ge0}\tau[\ell]\otimes
  S\left(\Delta[\ell](f)\right), &
  \mathcal R^{\Delta}_{1\otimes S}(f)
  &=\sum_{\ell\ge0}S\left(\Delta[\ell]\right)(f)\otimes
  \tau[\ell],
\end{aligned}\label{eq:DeltaExponentials}
\end{equation}
although $S(\Delta[\ell])$ does not belong to $\mathbb C[\Delta]$, but
exists as an operator on $M$.

\begin{lem}\label{lem:Extensiondistributions}
  Let $H^\prime$ be either $H_T$ or $\hat H_T$, with basis $\{e_i\}$
  and dual functions $\{e_i^*\}$ and let $M$ be an $H^\prime$-module
  such that the $H^\prime$ generators acting on $M$ can be expanded as
  \begin{equation}
    \label{eq:expansioneiDelta[n]}
    e_{i,M}=\sum_{\ell\ge0} c_{i,\ell}\Delta_M[\ell],\quad
    c_{i,\ell}=\Tr(\frac{e_i.\tau(\ell)}{\tau}).
  \end{equation}
Then, for all $m\in M$, we have
\begin{align*}
  \mathcal{L}^\Delta(m)&=\mathcal{L}^{H^\prime}(m), &
  \mathcal{R}^\Delta(m)&=\mathcal{R}^{H^\prime}(m),\\
  \mathcal{L}^\Delta_{1\otimes S}(m)&=\mathcal{L}^{H^\prime}_{1\otimes S}(m),
  &\mathcal{R}^\Delta_{1\otimes S}(m)&=\mathcal{R}^{H^\prime}_{1\otimes S}(m).
\end{align*}
\end{lem}

\begin{proof}
  We have an inclusion of $\Czpol$ in $(H^\prime)^*$ given by
\[
\tau(\ell)=\sum_i c_{i,\ell}e_i^*,\quad c_{i,\ell}=
\Tr\left(\frac{e_i.\tau(\ell)}{\tau}\right),
\]
so that $c_{i,\ell}$ are the same coefficients that appear in
\eqref{eq:expansioneiDelta[n]}. Then the proof of the lemma is a
simple computation. For instance,
\begin{align*}
  \mathcal{R}^\Delta_{1\otimes S}(m)&=
  \sum_{\ell\ge0}S(\Delta[\ell])(m)\tau(\ell) = \sum_{\ell\ge0}S(\Delta[\ell])(m)\sum_i c_{i,\ell}e_i^*\\
  &= \sum_{i}S\left(\sum_{\ell\ge0} c_{i,\ell}\Delta[\ell]\right)(m)e_i^*= \sum_{i}S(e_i)(m)e_i^*\\
  &=\mathcal{R}^{H^\prime}_{1\otimes S}(m).
\end{align*}
The rest is proved similarly.
\end{proof}

\begin{example}\label{example:CZpolExtension}
  Consider the $H_T$-module $M=\Czpol$. This satisfies the
  conditions of Lemma \ref{lem:Extensiondistributions}. Indeed, on
  $\Czpol$ the action of $\Delta$ is locally nilpotent and we have on
  $\Czpol$
  \begin{align*}
    \overline \Delta&=
    (1-T\inv)=(1-\frac1{1+\Delta})=-\sum_{k=1}^\infty(-\Delta)^k,\\
    \partial_\tau&=\log(T)=\log(1+\Delta)=\sum_{k=1}^\infty (-\Delta)^k/k.
  \end{align*}
  From this one easily checks the condition
  \eqref{eq:expansioneiDelta[n]}. So on $\Czpol$ we can replace the
  exponential operators \eqref{eq:HTexponentials} for $H_T$ by the
  ones for $\mathbb C[\Delta]$ in \eqref{eq:DeltaExponentials}.  In
  particular we can calculate the coproduct on $\Czpol$ use these
  simpler exponentials. For instance
\[
m^*_{1\otimes S}(\tau)=\mathcal{R}^{S}(\tau)=\tau\otimes
1-1\otimes \Delta^*[1]-1\otimes \Delta^*[-1]
\]
can be written as
\begin{equation}
m^*_{1\otimes S}(\tau)=\mathcal{R}^{\mathbb C[\Delta]}_{1\otimes S}(\tau)=\tau\otimes
1-1\otimes \tau.\label{eq:twcoproductau}
\end{equation}
More generally, on $\Czpol$ we have an action of $\hat H_T$, so we have
a twisted exponential for $\hat H_T$ on $\Czpol$:
\[
\hat {\mathcal{R}}^{S}(f)=\sum S(e_{n,k})(f)e_{n,k}^*,
\]
and we calculate the twisted coproduct on $\Czpol$ by thinking of
it as a subspace of the dual of $\hat H_T$. For instance, in this way
we get
\[
m_{1\otimes S}^*(\tau)=\tau\otimes 1-1\otimes e_{1,0}^*-1\otimes
e_{-1,0}^*-1\otimes e_{0,1}^*,
\]
which reduces to the prevous result \eqref{eq:twcoproductau} by
\eqref{eq:tauinhatHTdual}.\qed
\end{example}

\begin{example}\label{example:KsingExtension}
  Now take $M=\Ksing$. We think of it as an $H=\mathbb C[\Delta]$ or
  $H^\prime$-module by the $S$-twisted action: if $F\in\K$
  \[
P(\Delta,\partial_\tau).F=S\left(P(\Delta,\partial_\tau)\right)F.
\]
It is no longer true, as on $\Czpol$, that $\Delta$ acts locally
nilpotent on $\Ksing$. We would like to write
  \begin{equation}
    \partial_\tau=-\log(T\inv)=-\sum_{k>0}(\overline \Delta)^k/k,\label{eq:partialtauexpansion}
      \end{equation}
but the infinite sum on the right does not make sense in
$\Ksing$. However, we can think of $\Ksing$ as consisting of
distributions on $\Czpol$, via the embedding $F\in
\K\mapsto\mathcal{D}_F$, where
\[
\mathcal{D}_F(G)=\Tr(FG),\quad G\in \Czpol.
\]
Now we have for any distribution on $\Czpol$ the expansion
\[
\mathcal{D}=\sum_{\ell\ge0}\mathcal{D}(\tau[\ell])\frac
1{\tau[\ell+1]},
\] 
and \eqref{eq:partialtauexpansion} makes sense in the dual of $\Czpol$
and similarly we can expand
\[
\Delta=1+\Delta=1+\frac1{1-\overline\Delta}
\] as a power series in $\overline \Delta$. This is sufficient to
check \eqref{eq:expansioneiDelta[n]} and the Lemma
\ref{lem:Extensiondistributions} applies to $\Ksing$. Now
$\K=\Czpol\oplus\Ksing$, so by combining with Example
\ref{example:CZpolExtension} we see that we can replace on $\K$ the
exponential operators for $H_T$ or $\hat H_T$ by the simpler ones for
$\mathbb C[\Delta]$.\qed
\end{example}

\begin{example}\label{example:LocFinExtension}
  We can generalize Example \ref{example:CZpolExtension} as follows. Let $M$
  be any $H_T$-module that is locally finite, i.e., $M$ decomposes into
  generalized eigenspaces for $T$ and $T\inv$:
  \[
  M=\bigoplus_{\lambda\in \mathbb C^\times} M_\lambda,
  \]
  where
\[    M_\lambda=\{ m\in M\mid (T-\lambda)^k m=0,\quad k\gg0\}.
\]
Then by a slight variation of the proof of Lemma
\ref{lem:Extensiondistributions} we see that the exponential operators
of $\mathbb C[\Delta]$ and $H_T$ or $\hat H_T$ on locally finite $M$ give the same
distributions (on $H_T$ or $\hat H_T$).\qed
\end{example}

\begin{example}
  Let $M=H_T B$ be the free $H_T$-module generated by a single element
  $B$. In this case we cannot have \eqref{eq:expansioneiDelta[n]}, in
  particular $\overline \Delta$ can not be expressed in terms of
  $\Delta$. Now let $N=\mathbb C[\Delta] B$, the free $\mathbb
  C[\Delta]$-module generated by $B$. Then $N$ has a completion $\bar
  N$ consisting of infinite sums $\sum_{k\ge0} c_k\Delta[k]B$. Then we
  can define in $\bar N$ an action of $H_T$ or of $\hat H_T$, such
  that the conditions of Lemma \ref{lem:Extensiondistributions} are
  satisfied, and we can extend the exponential distribtions from
  $\mathbb C[\Delta]$ to $H_T$ or $\hat H_T$.\qed
\end{example}

\section{Localization of (Twisted) Coproducts and Dirac Distributions}
\label{sec:LocCoproductsDeltaDistributions}

We defined the twisted coproduct of $\phi\in H^*$ as the
distribution on $H\otimes H$ given by
\begin{equation}
\mStwod(\phi)(h_1\otimes h_2)= \phi(h_1Sh_2).\label{eq:defmStwod}
\end{equation}
We are interested in localizations $\K$ of $H^*$, and we
would like to extend the definition of $\mStwod$ to $\K$. We can not
use \eqref{eq:defmStwod}, as $\phi\in K$ is in general a
\emph{singular} function on $H$, and not defined on all $h\in H$.

Now we expressed the twisted coproduct $\mStwod$ on $H^*$ in terms of
twisted exponentials $\mathcal L^{S}, \mathcal R^S$. These expressions
make sense for $F\in \K$, $\K$ being an $H$-module, but one finds that
the two calculations give in this case in general a \emph{different}
result. So define for $F\in \K$ the obstruction of extending $\mStwod$
to $F$ as the \emph{Dirac Distribution} associated to $F$ by
\[
\delta(F)=  \mathcal R^{S}(F)-\mathcal L^{S}(F).
\]
In case $F\in H^*$ we have $\delta(F)=0$, of course.

\begin{example}\label{ex:HDtwistedexponentialsdelta}
We continue the example \ref{ex:HDexponentials}. We consider the
localization $K=\mathbb C[[t]][t\inv]$. Then
\begin{align*}
  \delta(\frac 1t)
  &=\exp(-t_2\partial_1)\frac1{t_1}-S\exp(t_1\partial_2)\frac1{t_2}\\
  &=\sum_{k\ge0} \frac{t_2^k}{t_1^{k+1}}+\sum_{k\ge0}
  \frac{t_1^k}{t_2^{k+1}}=\sum_{n\in\mathbb Z}t_1^nt_2^{-n-1}.
\end{align*}
This is the usual formal delta distribution, usually denoted by
$\delta(t_1-t_2)$ or $\delta(t_1,t_2)$, see \cite{MR99f:17033}.
\end{example}

\begin{example}
  \label{ex:HT-twisted-exponentialsdelta}
  We continue with Example \ref{ex:HTexponentials} and consider the
  localization $\K=M\inv\Czpol$ of $\Czpol\subset\Cz$. By Example
  \ref{example:KsingExtension} the exponential operators of
  $\mathbb C[\Delta]$ and of $H^\prime=H_T,\hat H_T$ give the same
  distribution on $K$. So define
  \begin{equation}\label{eq:defrholambda}
\rho(\frac1\tau)=\mathcal{R}_{1\otimes
  S}^\Delta(\frac1\tau)=\sum_{\ell\ge0}\frac{\tau_2(\ell)}{\tau_1(\ell+1)},\quad
\lambda(\frac1\tau)=\mathcal{L}_{1\otimes
  S}^\Delta(\frac1\tau)=-\sum_{\ell\ge0}\frac{\tau_1(\ell)}{\tau_2(\ell+1)}.
\end{equation}
then we also have
\begin{align*}
\rho(\frac1\tau)&=\sum_{n\in\mathbb Z} S(\Delta[n])\frac1\tau\otimes \Delta^*[n]
     =\sum_{n,k\in\mathbb Z,k\ge0 } S(e_{n,k})\frac1\tau\otimes
     e_{n,k}^*,\\
\lambda(\frac1\tau)&=-\sum_{n\in\mathbb Z}  \Delta^*[n]\otimes S(\Delta[n])\frac1\tau
     =-\sum_{n,k\in\mathbb Z,k\ge0 } e_{n,k}^*\otimes  S(e_{n,k})\frac1\tau,
  \end{align*}
and so we get 3 forms of the Dirac distribution
$\delta(\frac1\tau)=\rho(\frac1\tau)-\lambda(\frac1\tau)$:
\begin{align*}
  \delta(\frac1\tau)&=\sum_{n\in\mathbb Z} \tau_1(n)\tau_2(-n-1)\\
                       &=\sum_{n\in\mathbb Z}
                       S(\Delta[n])\frac1\tau\otimes\Delta^*[n]
                       +\sum_{n\in\mathbb Z} \Delta^*[n]\otimes S(\Delta[n])\frac1\tau\\
                       &=\sum_{n,k\in\mathbb Z,k\ge0}
                       S(e_{n,k})\frac1\tau\otimes e_{n,k}^*
                       +\sum_{n,k\in\mathbb Z,k\ge0}
                       e_{n,k}^*\otimes S(e_{n,k})\frac1\tau.
\end{align*}
\qed
\end{example}

\section{Difference Operators}
\label{sec:DiffOps}
Let $H$ be a Hopf algebra and let $V$ be vector space. Denote by $VH$
the vector space of $V$-valued $H$-operators: elements of $VH$ are
finite sums of the form
\begin{equation}
  P=\sum_{i} p_i e_i,\quad p_i\in V,\label{eq:DifOpCoefV}
\end{equation}
using a basis $\{e_i\}$ for $H$. Then $VH$ is an $H$-module: if $h\in
H, P\in V$ then we have
\[
h.P=\sum h^\prime_V(P)h^{\prime\prime},
\]
where the coproduct of $h$ is given by $\sum h^\prime\otimes h^{\prime\prime}$.

In case $H=H_T$ we will call elements of $VH_T$ difference operators
with values in $V$. In case $V$ is an $H_T$-Leibniz algebra, see
section \ref{sec:HTLeibnizExpo}, $VH_T$ will be an algebra in the
obvious way, but we will not need this at this point. 

The following operations make sense for any space $VH$. First we have
the adjoint map of $VH$ given, for $P=\sum p_{i} e_i$, by
\begin{equation}
P^*=\sum S(e_i).p_i.\label{eq:defadjoinedoperator}
\end{equation}
In case $V$ happens to be an $H$-Leibniz algebra this is an anti-involution.
Secondly there is the antipodal operator
\begin{equation}
  \label{eq:antipodaloperator}
  P^S=\sum_{i} p_i S(e_i).
\end{equation}
Thirdly we will need later the map 
\begin{equation}
VH\to V, \quad P\mapsto P_V=\sum_{i}S(e_i)(p_i).\label{eq:defP_V}
\end{equation}

\section{Expansions of Distributions}
\label{sec:ExpDistributions}

Let $W$ be a vector space. We denote by $W_n$ the space
of all $W$-valued distributions on $K^{\otimes n}$. For simplicity consider
first $\mathcal{D}\in W_1$. It has a decomposition
\[
\mathcal{D}=\mathcal{D}_\Hol+ \mathcal{D}_\Sing,
\]
where
\[
\mathcal{D}_\Hol:\Ksing\to W,\quad \mathcal{D}_\Sing:\Czpol\to W,
\]
i.e., the holomorphic part of the distribution vanishes on $\Czpol$
and the singular part on $\Ksing$.

Now we have a basis $\{\tau[\ell]\}_{\ell\ge0}$ for $\Czpol$ and a
basis $\{S(e_{n,k})\frac1\tau\}$ for $\Ksing$. This implies that we
can represent the distribution $\mathcal{D}$ and its components by
kernels
$\mathcal{D}(\tau)=\mathcal{D}_\Hol(\tau)+\mathcal{D}_\Sing(\tau)$,
where
\begin{align}
  \label{Dholkernel}
  \mathcal{D}_\Hol(\tau) &=\sum \mathcal{D}\left(S(e_{n,k}\frac1\tau)\right)e_{n,k}^*,\\
  \mathcal{D}_\Sing(\tau)&=\sum \mathcal{D}(\tau[\ell])\frac1{\tau[\ell+1]}.\notag
\end{align}
Then we have for all $F\in \K$
\[
\mathcal{D}(F)=\Tr(\mathcal{D}(\tau)F).
\]
To make sense of this formula, let 
\begin{equation}
  \label{eq:defhatK}
  \hat \K=\hatCzpol\oplus \Ksing,\quad\hatCzpol=\oplus \mathbb C e_{n,k}^*\subset (\hat H_T)^*.
\end{equation}
Then the trace extends to $\hat \K$. Now $\mathcal{D}(\tau)F$ is an
infinite sum of terms in $W\otimes \hat \K$, of which at most a finite
number have a second component with nonvanishing trace. The trace then
also extends then to such infinite sums. One also sees that we can let
$F\in\hat \K$, and $\mathcal{D}(F)$ is well defined in that case.

Now note that in the expression \eqref{Dholkernel} for
$\mathcal{D}_\Hol$ the $\hat H_T$ form of the distribution
$\rho(\frac1\tau)$ appears, see Example
\ref{ex:HT-twisted-exponentialsdelta}. So we can write
$\mathcal{D}_\Hol$ in the form adapted to $\mathbb C[\Delta]$:
\[
\mathcal{D}_\Hol(\tau)=\sum_{\ell\ge0}\mathcal{D}(\frac1{\tau[\ell+1]})\tau[\ell],
\] 
where we think of $\mathcal{D}_\Hol$ as a distribution on $\mathbb
C[\Delta]\simeq S(\mathbb C[\Delta])\frac1\tau$ that can be extended
to a distribution on all of $\hat H_T\simeq \Ksing$, as discussed in
section \ref{sec:extExpDistr}. This means that we have a convenient
uniform expression for the kernel of any distribution in one variable:
\begin{equation}
\mathcal{D}(\tau)=\sum_{n\in\mathbb
  Z}\mathcal{D}(\tau(n))\tau(-n-1),\label{eq:expansiondistribution}
\end{equation}
using the notation \eqref{eq:deftau-n-1}.
If $F\in\K$ it defines \emph{a rational distribution} $\mathcal{D}_F$ with
$F$ as kernel:
\begin{equation}
\mathcal{D}_F(G)=\Tr(FG).\label{eq:DefRationalDistribution}
\end{equation}
This means that any $F\in\K$ (or rather its image in $K^*$) has a \emph{formal expansion} 
\begin{equation}
  F=\sum_{n=-N}^\infty F_n\tau(-n-1),\quad F_n=\Tr(F\tau(n)).\label{eq:formalexpansionK}
\end{equation}
For example we have the following infinite expansion of the rational
function $\frac1{\tau^2}\in\Ksing$:
\[
\frac1{\tau^2}=\sum_{k\ge0}\frac{k!}{\tau(k+2)},
\]
if we think of $\frac1{\tau^2}$ as a distribution (on $\Czpol$).

In the same way we have for $\mathcal{D}\in W_2$ a kernel
\begin{equation}
\mathcal{D}(\tau_1,\tau_2)=\sum_{m,n\in\mathbb
  Z}\mathcal{D}(\tau(m)\otimes\tau(n))\tau_1(-m-1)\tau_2(-n-1),\label{eq:ExpansionTwovariableDist}
\end{equation}
with for all $F,G\in \K$
\[
\mathcal{D}(F\otimes
G)=\Tr_{\tau_1}\Tr_{\tau_2}(\mathcal{D}(\tau_1,\tau_2)F(\tau_1)G(\tau_2)).
\]
We will frequently identify a distribution with its kernel. We will also write
\[
\rho(\frac1\tau)=\rho(\tau_1,\tau_2),\quad
\lambda(\frac1\tau)=\lambda(\tau_1,\tau_2),\quad
\delta(\frac1\tau)=\delta(\tau_1,\tau_2),
\]
and more generally we will use the notation $a(\tau_i,\tau_j)$ for a
distribution depending on variables $\tau_i(n),\tau_j(n)$, where
\[
\tau_i(n)=1\otimes\dots \otimes\tau(n)\otimes \dots 1,
\]
with $\tau(n)$ on the $i$-th position.

Similarly any distribution  $\mathcal{D}\in W_n$ has an
expansion in $\tau_i$, $i=1,\dots,n$. We will write for this reason
$W_n=W[[\tau_1^{\pm1},\tau_2^{\pm1},\dots,\tau_n^{\pm1}]]$, and we
will denote by $W[[\tau]]\subset W[[\tau^{\pm1}]]$ the subspace of
holomorphic distributions (i.e., distributions on $\Ksing$), and by
$W[[\tau\inv]]\tau\inv$ the subspace of singular distributions, those on
$\Czpol$. 

\section{Some Properties of Distributions}
\label{sec:Propdist}

Distributions on $\K^{\otimes n}$ can be multiplied by difference
operators, via the action of difference operators on $K$. For instance
if $P\in VH_T$ and $\mathcal{D}\in \mathbb
C[[\tau_1^{\pm1},\tau_2^{\pm1}]]$ then we define
\begin{align*}
  P_1\mathcal{D}(F\otimes G)&=\mathcal{D}(P^*F\otimes G),\\
  P_2\mathcal{D}(F\otimes G)&=\mathcal{D}(F\otimes P^* G),
\end{align*}
using the adjoint operation on $VH_T$, see
\eqref{eq:defadjoinedoperator}.  Here we extend the $\mathbb C$-valued
distribution $\mathcal{D}$ on $\K\otimes\K$ to one on
$V\otimes\K\otimes\K$ or $K\otimes V\otimes \K$ by linearity in the
obvious way, to obtain a $V$-valued distribution.

We can not multiply distributions except in special cases. 
Any distribution $\mathcal{D}\in W_n$ can be multiplied by a rational
distribution in $n$ ways. For example, for $n=2$, $F,L,M\in\K$ and
$\mathcal{D}_F$ the rational distribution
(\ref{eq:DefRationalDistribution}) we have
\[
(\mathcal{D}_{F})_1\mathcal{D}(L\otimes M)=\mathcal{D}(FL\otimes
M),\quad 
(\mathcal{D}_{F})_2\mathcal{D}(L\otimes M)=\mathcal{D}(L\otimes F
M).
\]
In general we would like to define the product of a distribution
$\mathcal{D}$ in one variable with a two variable distribution
$\mathcal{E}$ by
\begin{equation}
  \label{eq:hypoproductdistr}
  \mathcal{D}_1.\mathcal{E}(F\otimes
  G)=\Tr_{\tau_1}\Tr_{\tau_2}\left(
\mathcal{D}(\tau_1)\mathcal{E}(\tau_1,\tau_2)F(\tau_1)G(\tau_2)\right).
\end{equation}
However, for general $\mathcal{D}$ and $\mathcal{E}$ the products of
kernels doesn't make sense.  We have the special two variable
distributions $\rho(\frac1\tau)$, $ \lambda(\frac1\tau)$ and
$\delta(\frac1\tau)$. The claim now is that in case
$\mathcal{E}=\rho,\lambda$ or $\delta$ the definition
\eqref{eq:hypoproductdistr} \emph{does} make sense for any
$\mathcal{D}$. This can be easily checked.

To investigate the properties of such products we need the notion of
the trace of a distribution in several variable with respect to one of
the variables. For instance, let $\mathcal{D}(\tau_1,\tau_2)$ be a
$W$-valued distribution on $\K\otimes \K$. We can produce
distributions in a single variable from $\mathcal{D}$ by
\begin{align*}
\Tr_{\tau_1}(\mathcal{D}(\tau_1,\tau_2))(G)&=\mathcal{D}(\tau_1,\tau_2)(1_{\K}\otimes G), \\
\Tr_{\tau_2}(\mathcal{D}(\tau_1,\tau_2))(F)&=\mathcal{D}(\tau_1,\tau_2)(F\otimes
1_{\K}).
\end{align*}
For instance 
\[
\Tr_{\tau_1}(\delta(\tau_1,\tau_2))=\Tr_{\tau_2}(\delta(\tau_1,\tau_2)))=1,
\]
where $1=1_{\K}$ is the rational distribution on $\K$ with kernel $1$,
so
\[
\Tr_{\tau_1}\left(\delta(\tau_1,\tau_2)\right)(F)=\Tr(F)=1_{\K}(F).
\]
\begin{lem}\label{lem:TrDrhoRlambda}
  Let $\mathcal{D}(\tau)$ be a $W$-valued distribution on $\K$. Then
  \begin{align*}
    \Tr_{\tau_1}(\mathcal{D}(\tau_1)\rho(\tau_1,\tau_2))&=\mathcal{D}_\Hol(\tau_2),&\quad
        \Tr_{\tau_2}(\mathcal{D}(\tau_2)\rho(\tau_1,\tau_2))&=\mathcal{D}_\Sing(\tau_1),\\
    \Tr_{\tau_1}(\mathcal{D}(\tau_1)\lambda(\tau_1,\tau_2))&=-\mathcal{D}_\Sing(\tau_2),&\quad 
    \Tr_{\tau_1}(\mathcal{D}(\tau_1)\lambda(\tau_1,\tau_2))&=-\mathcal{D}_\Hol(\tau_1) 
  \end{align*}
so that
\[
    \Tr_{\tau_1}(\mathcal{D}(\tau_1)\delta(\tau_1,\tau_2))=\mathcal{D}(\tau_2),\quad
    \Tr_{\tau_2}(\mathcal{D}(\tau_2)\delta(\tau_1,\tau_2))=\mathcal{D}(\tau_1).\]
\end{lem}

\begin{proof}
  This follows from expanding $\mathcal{D}$ as in
  (\ref{eq:ExpansionTwovariableDist}) and using the explicit forms for
  $\rho,\lambda$, see section
  \ref{sec:LocCoproductsDeltaDistributions}, combined with
  orthogonality relations from section \ref{sec:orth-relat}.
\end{proof}

We use this to give explicit formulas for the product of an arbitrary
distribution $\mathcal{D}$ with $\rho,\lambda,\delta$.
\begin{lem}\label{lem:propertiesrholamdadelta}
  For all $F,G\in\K$ we have
  \begin{Myenumerate}
  \item \label{lem:item:rhoFHolGsing} $\rho(\frac1\tau)(F\otimes G)=\Tr(F_\Hol
    G_\Sing)$, $\lambda(\frac1\tau)(F\otimes G)=-\Tr(F_\Sing G_\Hol)$.
  \item \label{lem:item:deltaFG}$\delta(\frac1\tau)(F\otimes G)=\Tr(FG)$.
  \item \label{lem:item:Drholambdadelta}For all $\mathcal{D}\in W[[\tau^{\pm1}]]$ 
    \begin{align*}
      \mathcal{D}(\tau_1)\rho(\frac1\tau)(F\otimes G)&= \mathcal{D}(F
      G_\Sing),&\quad\mathcal{D}(\tau_2)\rho(\frac1\tau)(F\otimes G)&= \mathcal{D}(F_\Hol
      G)\\
\mathcal{D}(\tau_1)\lambda(\frac1\tau)(F\otimes G)&=- \mathcal{D}(F
      G_\Hol),&\quad
\mathcal{D}(\tau_2)\lambda(\frac1\tau)(F\otimes G)&=- \mathcal{D}(F_\Sing
      G)\\
\mathcal{D}(\tau_1)\delta(\frac1\tau)(F\otimes G)&= \mathcal{D}(F
      G),&\quad \mathcal{D}(\tau_2)\delta(\frac1\tau)&=\mathcal{D}(F
      G)
    \end{align*}
  \end{Myenumerate}
\end{lem}

\begin{proof}
  Part \eqref{lem:item:rhoFHolGsing} follows from the explicit forms
  for $\rho,\lambda$ in Example \ref{ex:HT-twisted-exponentialsdelta},
  combined with the expansions of $F$ and $G$ (see
  \eqref{eq:formalexpansionK})) and with the orthogonality relations
  of section \ref{sec:orth-relat}. The part \eqref{lem:item:deltaFG}
  follows from this. For 
  part \eqref{lem:item:Drholambdadelta} we calculate for instance
  \begin{align*}
    \mathcal{D}.\rho(\frac1\tau)(F\otimes G)&=\Tr_{\tau_1}\Tr_{\tau_2}
\left(\mathcal{D}(\tau_1)\rho(\tau_1,\tau_2)F(\tau_1)G(\tau_2)\right)\\
&=\Tr_{\tau_1}\left(\mathcal{D}(\tau_1)F(\tau_1)G_\Sing(\tau_1)\right)\\
&=\mathcal{D}(FG_\Sing),
  \end{align*}
since by Lemma \ref{lem:TrDrhoRlambda} we have
\[
\Tr_{\tau_2}\left(\rho(\tau_1,\tau_2)
  G(\tau_2)\right)=G_\Sing(\tau_1),
\]
identifying $G$ with the rational distribution of which it is the
kernel. The other identities follow similarly.
\end{proof}

We can not only multiply $\rho(\frac1\tau),\lambda(\frac1\tau)$ and
$\delta(\frac1\tau)$ by distributions, but also by distribution valued
difference operators. For instance, let $V=W[[\tau^{\pm1}]]$ and
consider distributions defined by, if $P=\sum P_n(\tau)\Delta[n]\in
VH_T$, 
\begin{equation}
  \label{eq:VHTdelta}
  \begin{split}
    P(\tau_1)\delta(\tau_1,\tau_2)(F\otimes G)&=\sum_{n\in\mathbb
      Z}P_n(F\Delta[n]G),\\ 
    P(\tau_2)\delta(\tau_1,\tau_2)(F\otimes G)&=\sum_{n\in\mathbb
      Z}P_n(\Delta[n]\langle F\rangle G).
  \end{split}
\end{equation}
\begin{lem}
  For $P(\tau)\in VH_T$, where $V=W[[\tau^{\pm1}]]$ we have
  \begin{equation}
    \label{eq:dualDifferencedelta}
    P(\tau_1)\delta(\tau_1,\tau_2)=P^*(\tau_2)\delta(\tau_1,\tau_2).
  \end{equation}
\end{lem}

\section{Distributions of the form $P\frac1\tau$}
\label{sec:DistributionsForm1tau}

We will need, if $P\in VH_T$, expressions of the form
\begin{equation}
  \label{eq:VHTdistributions}
  P^S\frac1\tau=\sum_nP_n S(\Delta[n])\frac1\tau.
\end{equation}
This is a $V$-valued distribution on $\K$ (in fact on
$\Czpol$). 
Since $V$ is an $H_T$-module, we can let the exponential operator
$\mathcal{R}$ of $V$ act on expressions
\eqref{eq:VHTdistributions}. Recall from \eqref{eq:oppositeactionKsing}
\begin{equation}
  \Delta^*[m].S(\Delta[n])\frac1\tau=0,\label{eq:oppsitesignstau}
\end{equation}
in case of opposite sign $m>0>n$ or $m<0<n$. In case of the same sign
($m,n\ge0$ or $m,n\le0$) we have
\begin{equation}
  \Delta^*[m].S(\Delta[n])\frac1\tau=\tau[m]S(\Delta[n])\frac1\tau.\label{eq:samesigntau}
  \end{equation}
We use this to note that for $n\ne0$ we have
\begin{equation}
  \label{eq:ExpPDist}
  \mathcal{R}_V P_n S(\Delta[n])\frac1\tau=
  \begin{cases}
    \sum\limits_{k\ge0}\Delta_V[k]P_n \tau[k]S(\Delta[n])\frac1\tau & n>0\\
    \sum\limits_{k\ge0}\overline{\Delta_V[k]}P_n \overline{\tau[k]}S(\Delta[n])\frac1\tau & n<0
  \end{cases}
\end{equation}
Recall, see section \ref{sec:ExpDistributions}, the notion of the singular part of a distribution
$\mathcal{D}$, We write $\Sing(\mathcal{D})=\mathcal{D}_\Sing$.
\begin{lem}\label{lem:PstarrhoSing}
  \begin{align}
(P^*)^S
\frac1\tau&=\Sing\left(\mathcal R_V(\tau)P\frac1\tau\right),\label{eq:P*S}\\    
(P^*)
\frac1\tau&=\Sing\left(\mathcal R^S_V(\tau)P^S\frac1\tau\right).    \label{eq:P*}
  \end{align}
\end{lem}
\begin{proof} 
  For $P=P_0\Delta[0]$ \eqref{eq:P*S} is clear. Let
  $P_+=P_n\Delta[n]$, for $n>0$. Then, using the coproduct
  \eqref{eq:coproddiffbasis}, the factorization
  (\ref{eq:factorizationtau[k]}) and \eqref{eq:ExpPDist} we find
  \begin{align*}
    (P^*)^S\frac1\tau 
    &=\left(S(\Delta[n])P_n\right)^S\frac1\tau=(-1)^n\sum_{s=0}^n
     \overline{\Delta_V}[s]P_n(T^{-s}\overline{\Delta[n-s]})^S\frac1\tau\\
      &=(-1)^n\sum_{s=0}^n
      \overline{\Delta_V}[s]P_nT^{s}S(\overline{\Delta[n-s]})\frac1\tau\\
      &=(-1)^n\sum_{s=0}^n
      \overline{\Delta_V}[s]P_nT^{s}\frac1{\overline{\tau(1+n-s)}}\\
      &=\sum_{s=0}^n
      \overline{\Delta_V}[s]P_n\overline{\tau[s]}\Delta[n]\frac1{\tau}\\
      &=\Sing(\mathcal{R}_VP_+\frac1\tau).
  \end{align*}
  In the same way one calculates \eqref{eq:P*S} for $P_-=P_n\Delta[n]$
  for $n<0$. The proof of \eqref{eq:P*} is similar.
\end{proof}

\begin{cor}\label{cor:TbTaSing}
  For all $a,b\in \mathbb Z$ we have
\[
T^{-b}\langle T^{-a}\langle\frac1{\tau_1}\rangle
\frac1{\tau_2}\rangle=\Sing_{\tau_2}\left(\mathcal
  R_{W}(\tau_{2})(T_1^{-a}\frac1{\tau_1}T^b)\frac1{\tau_2} \right),
\]
where $W$ is the $H_T$-module generated by $\frac1{\tau_1}$
\end{cor}

\begin{proof}
  This is \eqref{eq:P*S} in case $P=T^{-a}_1\langle \frac1{\tau}\rangle  T^b$.
\end{proof}

\section{Rationality of Distributions on $\K$}
\label{sec:RatDist}

Recall the notion of a rational distribution, see
(\ref{eq:DefRationalDistribution}). We say that a distribution
$\mathcal{D}$ on $\K$ has rational singularities in case the singular
part $\mathcal{D}_\Sing$ is rational, i.e., there is
$\mathcal{K}_{\mathcal{D}}^\Sing\in \Ksing$ such that for all $G\in\Czpol$
\[
\mathcal{D}_\Sing(G)=\Tr(\mathcal{K}_{\mathcal{D}}^\Sing G).
\]
\begin{lem}
  \label{lem:RationalSing}
A $W$-valued formal distribution $\mathcal{D}$ has rational singularities if and
only if $F.\mathcal{D}_\Sing=0$ for some $F=\prod_{(n,d)}(\tau-n)^d\in\Czpol$.
\end{lem}
\begin{proof}
  
If $\mathcal{D}$ has rational singularities with singular kernel
  $\mathcal{K}_{\mathcal{D}}^\Sing\in\Ksing$, then there is some
  $F=\prod_{(n,d)}(\tau-n)^d\in\Czpol$ such that
  $F.\mathcal{K}_{\mathcal{D}}^\Sing\in\Czpol$. Then 
\[
F.\mathcal{D}_\Sing(G)=\Tr(\mathcal{K}_{\mathcal{D}}^\Sing F G)=0,
\]
for all $G\in\Czpol$, since the argument of $\Tr$ is here nonsingular.

Conversely, let
\[\mathcal{D}_\Sing=\sum_{k\ge0} \frac{d_k}{\tau(k+1)},\quad d_k\in W,
\]
such that $F.\mathcal{D}_\Sing=0$ (as a distribution on $\Czpol$), for
\[
F=\prod_{(n,d)}(\tau-n)^d=\tau(N)+\sum_{n=1}^N f_k\tau(N-k).
\]
Now the action of $\Czpol$ on singular distributions is given by
\begin{equation}\label{eq:actionKT}
  \begin{split}
\tau(m+\ell).  \frac 1{\tau(\ell)}&=  T^{-\ell}\tau(m)=0 \mod \Czpol\\
\tau(m).  \frac 1{\tau(\ell+m)}&=  T^{-\ell}\frac
1{\tau(m)}=(1-\overline\Delta)^\ell \frac 1 {\tau(m)}.
  \end{split}
\end{equation} 
By using the last equation of
\eqref{eq:actionKT} we can write  $F.\mathcal{D}_\Sing=0$
as a collection of recursion equations for the coefficients of
$\mathcal{D}_\Sing$:
\[
d_{N+p}+ \sum_{k=1}^N\gamma_{p,k} d_{N+p-k}=0,\quad p\ge0,
\quad \gamma_{p,k}\in \mathbb C.
\]
This leaves $d_0,d_1,\dots,d_{N-1}$ arbitrary, and determines uniquely
$d_{p+N}$, for $p\ge0$. So the subspace of distributions that solve
$F.\mathcal{D}=0$ is $n$ dimensional. Now consider the rational
element
\begin{equation}
  \label{eq:ratfunctionf}
A=\sum_{(n,d)}\sum_{i=1}^d
\frac{a_{(n,i)}}{(\tau-n)^{i}},
\end{equation}
containing $N=\deg(F)$ arbitrary constants $a_{(n,i)}$. In the
summation $(n,d)$ runs over the same set of pairs that occur in
$F$. Then $F.A=0\mod \Czpol$. We can expand $A$ (or rather its image
in $\K^*$) as 
\[
A=\sum_{n\ge0}a_{(n)}\frac1{\tau(n1)},
\]
and the first $N$ coefficients are determined by a system of equations
of the form
\[
a_{(j)}=\sum_{(n,d)}\sum_{i=1}^d c_j^{(n,i)}a_{(n,i)},\quad
j=0,1,\dots, N-1,
\]
where $c_j^{(n,i)}\in\mathbb C$. One easily checks that the $N\times
N$ matrix $(c_j^{(n,i)})$ is 
invertible, so that any solution of the
equation $F.\mathcal{D}_\Sing=0$ determines uniquely a rational function of the form
\eqref{eq:ratfunctionf}.
\end{proof}

\section{Rationality of Multivariable Distributions}
\label{sec:RatMultiDistr}

We define
\begin{equation}
  \tau_{1\otimes S}[\ell]=\mStwod(\tau[\ell]).
\end{equation}

\begin{lem}\label{lem:coefficientDeltaexpansion} Let
  $P(\tau)=\sum_{n}P_n(\tau)\Delta[n]$ be a difference operator with
  values in the space of $W$-valued distributions on $\K$. For all $m\ge0$ we have
  \begin{equation}
    \label{eq:coeffDistr}
  \Tr_{\tau_1}\left(\tau_{1\otimes S}[m]P(\tau_2)\delta(\tau_1,\tau_2)\right)=P_m(\tau_2).  
  \end{equation}
\end{lem}

\begin{proof}
  Let $X$ be the left hand side in \eqref{eq:coeffDistr}. We have
  \begin{align*}
    X&=\sum_{\ell\ge0}\Tr\left( \tau_{1\otimes
        S}[m] P_\ell(\tau_2)S(\Delta_1[\ell])\delta(\tau_1,\tau_2)\right)\\
    &=\sum_{\ell\ge0}\Tr\left( \Delta_1[\ell]\langle\tau_{1\otimes
        S}[m]\rangle P_\ell(\tau_2)\delta(\tau_1,\tau_2)\right).
  \end{align*}
  Now by using Lemma \ref{lem:HpropertiesExponentials} we have
  \[
  \Delta_1[\ell]\tau_{1\otimes S}[m]=\mStwod(\Delta[\ell]\tau[m]),
  \]
  which is zero for $\ell>m$. In this case there is no contribution to
  $X$. So assume $m=\ell+s$, $s\ge0$. Then
  \[ 
  \Delta[\ell]\tau[\ell+s]=\gamma \tau[s],
  \] 
  for some non zero constant $\gamma$. Now, by the multiplicativity of
  exponentials, see section \ref{sec:MultTwistedExp}, we have
  \[
  \tau_{1\otimes
    S}[s]\delta(\tau_1,\tau_2)=\delta\left(\frac{\tau[s]}{\tau}\right),
  \]
  which is zero unless $s=0$.
\end{proof}

Recall the formal expansion of $f\in\K$, see
\eqref{eq:formalexpansionK}. Similarly
we will often consider distributions that have a formal expansion
\begin{equation}
  \label{eq:expansiondelta}
  a(\tau_1,\tau_2)=\sum_{j=0}^\infty a_j(\tau_2)\Delta_2[j]\delta(\tau_1,\tau_2).  
\end{equation}
The coefficients of such an expansion are traces by Lemma \ref{lem:coefficientDeltaexpansion}:
\begin{equation}
  \label{eq:coefastrace}
  a_j(\tau_2)=\Tr_{\tau_1}\left(a(\tau_1,\tau_2) \tau_{1\otimes S}(j)\right).
\end{equation}
Let $V^{\text{conv}}[[\tau_1^{\pm},\tau_2^{\pm}]]$ be the space of $V$-valued distributions
$a(\tau,\tau_2)$ such that
\[
\pi(a(\tau_1,\tau_2))=\sum_{j=0}^\infty a_j(\tau_2) \Delta_2(j)
\delta(\tau_1,\tau_2)\]
converges, where $a_j(\tau)$ is given by \eqref{eq:coefastrace}.

Recall from section \ref{sec:hatH_Tanddual} that the singular part
$\Ksing$ is generated by $\frac 1\tau$ over $\hat
H_T=H_T[\partial_\tau]$. If
\[
\kappa=(-\partial_\tau)^d S(\Delta[n])\frac 1\tau \in \Ksing,
\]
then, by Lemma \ref{lem:HpropertiesExponentials},
\begin{equation}
  \label{eq:deltakappa}
\delta(\kappa)=\Delta_2(n)(\partial_{\tau_2})^d\delta(\tau_1,\tau_2).
\end{equation}
Recall also the notion of an element $\sigma\in V[[\tau^{\pm1}]]$
having rational singularities, see section \ref{sec:ExpansionsK}.
Similarly we say that a distribution $a(\tau_1,\tau_2)\in
V^{\text{conv}}[[\tau_1^{\pm},\tau_2^{\pm}]]$ has \emph{rational
  singularities} (or is rational) if there is a $f\in\Czpol$ such that
\[
f_{1\otimes S}a(\tau_1,\tau_2)=0,
\]
where $f_{1\otimes S}=\md_{1\otimes S}(f)\in \Czpol\otimes\Czpol$.
\begin{lem}\label{lem:rationaldistrtwovariables}
  A distribution $a(\tau_1,\tau_2)\in V^{\mathrm{conv}}[[\tau_1^{\pm},\tau_2^{\pm}]]$ is rational if an only if 
it is a finite sum
\begin{equation}
  \label{eq:rationaldistribution}
a(\tau_1,\tau_2)=\sum_{n\in\mathbb
  Z,d\ge0}c_{n,d}(\tau_2)\Delta_2(n)(\partial_{\tau_2})^d\delta(\tau_1,\tau_2).
\end{equation}
\end{lem}
The proof of this Lemma is similar to that of Lemma \ref{lem:RationalSing}.

Note that the expansion \eqref{eq:expansiondelta} of a rational
$a(\tau_1,\tau_2)$ will in general be infinite.

\chapter{Singular Hamiltonian Structures}
\label{chap:SingHamStruc}

\section{Overview}
\label{sec:IntroKVPA}
Let $V$ be a commutative $H_T$-Leibniz algebra, see section
\ref{sec:HTLeibnizExpo}. We will show in section
\ref{sec:H_LeibnizAlgSpec} that the affine variety of $V$ can be
identified with the set $V(H^*)=\Hom_{H-\mathrm{alg}}(V,H^*)$ of
homomorphisms of $H$-Leibniz algebras from $V$ to $H^*$. We refer to
$V(H^*)$ as the set on regular points of V (or rather of the variety
asssociated to $V$). If we replace $H^*$ by some other $H$-Leibniz
algebra $L$ we obtain a set $V(L)=\Hom_{H-\mathrm{alg}}(V,L)$, which
we refer to as a set of singular points of $V$. The notion of
singularity is distinct from the usual concept of a singular point of
an algebraic variety and motivated by the example of $H=H_D$. In this
case, if $V$ is a free differential algebra with one generator, so
that $V=\mathbb C[v^{(i)}], Dv^{(i)}=v^{(i+1)}$, the set $V(H^*)$ is
isomorphic to $H^*=\mathbb C[[t]]$, the set of nonsingular functions
on the formal disk, whereas if we take $\K=\mathbb C[[t]][t\inv]$ the
set $V(K)$ is isomorphic to $K$, a set of singular functions on the
formal disk.

In section \ref{sec:Classicalfields} we will define classical
fields as generating series of functions on $V(L)$, for the case
$H=H_T$. The rest of this chapter develops then the theory of
(singular) Poisson brackets between classical fields.
 
\section{Spectrum and Singular Points}
\label{sec:H_LeibnizAlgSpec}
Let $V$ be a commutative $H$-Leibniz algebra. The $H$-action on $V$ induces an
$H$-action on $V^*$ and by duality we obtain a coaction
\[
m_V^*\colon  V\to (H\otimes V^*)^*, \quad v\mapsto \{ h\otimes v^*\mapsto v^*(hv)\}.
\]
We can expand $\md_V(v)$ in an infinite sum of elements of $V\otimes
H_T^*$. As we discussed in section \ref{sec:exp_op} we can
calculate $\md_V$ using the exponential operator $\mathcal R$
introduced there: introducing a basis $\{e_i\}$ for $H$ and functions
$e_i^*$ such that $e_i^*(e_j)=\delta_{ij}$ we have
\[
m^*_V(v)=\mathcal R(v)=\sum e_i(v)\otimes e_i^*.
\]
Then $\md_V$ is  a multiplicative map (see Lemma
\ref{lem:multiplicativityexps}) and satisfies $H$-covariance
(see Lemma \ref{lem:HpropertiesExponentials}). Since we
assume $V$ to be a commutative algebra we can define its spectrum
$\Specm(V)$: a point $\sigma\in\Specm(V)$ is a $\mathbb C$-algebra
homomorphism $\sigma\colon V\to \mathbb C$. Given such $\sigma$ we
construct, using the exponential map $\md_V$, a $H$-Leibniz algebra
homomorphism
\[
\hat \sigma\colon V\to H^*, \quad v\mapsto \sigma\md_V(v)=\sum
\sigma(e_i v)e^*_i .
\]
Indeed, this is an algebra homomorphism:
\begin{align*}
  \hat\sigma(v_1v_2)&= (\sigma\otimes 1)\md_V(v_1v_2)
  = (\sigma\otimes 1)\md_V(v_1)\md_V(v_2),\\
  &= \hat \sigma(v_1)\hat\sigma(v_2).
\end{align*}
Similarly, using part \eqref{Part:1} of Lemma
\ref{lem:HpropertiesExponentials}, we see that $\hat\sigma$ is an
$H$-module morphism:
\[
\hat\sigma(hv)=(\sigma\otimes 1)(1\otimes h_{H^*})\md_V(v)=(1\otimes
h_{H^*})(\sigma\otimes1) \md_V(v)=h_{H^*}\hat\sigma(v).
\]
Denote by $V(H^*)$ the set of all $H$-Leibniz algebra morphisms
$\tau\colon V\to H^*$. So we have constructed a map $\Specm(V)\to V(H^*)$,
$\sigma\mapsto \hat\sigma$. Conversely we get a map $V(H^*)\to \Specm(V)$
by composing $\tau\in V(H^*)$ with the canonical algebra homomorphism
\[
i^*\colon H^*\to \mathbb C,
\]
dual to the unit $i\colon \mathbb C \to H$. Then one easily checks
that these maps are inverses of each other:
\[
\hat{(i^*\circ \tau)}=\tau,\quad i^*\circ
\hat\sigma=\sigma.
\]
Hence as sets we find an isomorphism
\[
V(H^*)\simeq \Specm(V),
\]
and we refer to $V(H^*)$ as the set of nonsingular points of the
(variety associated to the) $H$-Leibniz algebra $V$.

\section{Classical Fields}
\label{sec:Classicalfields}

From now on we restrict ourselves for convenience to the case $H=H_T$
and consider the localization $\K=M\inv \Czpol$ described in
section \ref{sec:localizationHTcirc}.

Let $V$ be a commutative $H_T$-Leibniz algebra and let $V(\K)$ be the
set of $H_T$-Leibniz algebra homomorphisms $\hat\sigma\colon V\to \K$.
(By taking the non singular part of $\hat \sigma$ we get a map $\hat
\sigma_\Hol\colon V\to H_T^*$, i.e., $\hat \sigma_\Hol$ corresponds to
a point of $\Specm (V)$.)

For $f\in V$ we we define the \emph{classical field} $C(f)$ as the
distribution on $\K$ with values in the space of functions $V(\K)\to
\mathbb C$, defined tautologically by:
\[
C(f)(F)(\hat\sigma)=\Tr\left(\hat\sigma(f)F\right),\quad
F\in\K,\hat\sigma\in V(\K).
\]
For each $\hat \sigma$ we get a rational distribution $C(f)(\hat
\sigma)$ on $\K$ (with values in $\mathbb C$). Rational distributions
can be multiplied, see Section \ref{sec:Propdist}. This allows us to
define the (commutative, associative) product $C(f).C(g)$ of classical
fields associated to $f,g\in V$ by
\[
C(f).C(g)(\hat \sigma)=C(f)(\hat\sigma)C(g)(\hat\sigma).
\]
Then we have
\[
C(f).C(g)=C(fg).
\]
Also we have an $H_T$ action on distributions and this turns the space
$V(\tau)$ of classical fields (spanned by $C(f), f\in V$) into a
commutative $H_T$-Leibniz algebra.

\section{Expansions of Classical Fields}
\label{sec:ExpClassFields}
Any distribution has an expansion (\ref{eq:expansiondistribution}). In
the case of classical fields this becomes:
\[
C(f)=\sum_{n\in\mathbb Z} f_{(n)}\tau(-n-1),
\]
where $f_{(n)}=C(f)(\tau(n))$ is the function
\begin{equation}
f_{(n)}\colon V(\K)\to \mathbb C, \quad\hat\sigma\mapsto
\Tr(\hat\sigma(f)\tau(n)).\label{eq:Deff_(n)ClassField}
\end{equation}
We will write $f(\tau_i)$ for the expansion of $C(f)$ in the
variable $\tau_i$. 

Consider the holomorphic part
$C(f)_\Hol=\sum_{k\ge0}f_{(-k-1)}\tau(k)$ of the classical field
$C(f)$. We can identify the components $f_{(-k-1)}$ of $C(f)_\Hol$
with the elements $\Delta[k]f\in V$. Indeed, for any $\hat \sigma\in
V(K)$ we have
\begin{align*}
  f_{(-k-1)}(\hat\sigma)&=\Tr\left(\hat
    \sigma(f)\frac1{\tau(k+1)}\right)=\Tr\left(\hat \sigma_\Hol(f)S(\Delta[k]\frac1{\tau}\right)\\
  &=\Tr\left(\hat \sigma_\Hol(\Delta[k]f)\frac1{\tau}\right)\\
  &=\sigma(\Delta[k]f),
\end{align*}
if $\sigma\in\Specm(V)$ corresponds to $\hat\sigma_\Hol\in V(H_T^*)$.
This mean that we can identify the holomorphic part of the classical
field $C(f)$ with the distribution
$\mathcal{R}_V(f)=\sum_{k\ge0}\Delta[k]f\otimes \tau(k)$. In
particular the constant term of a classical fields $C(f)=f(\tau)$ can
be identified with the element $f\in V$:
\begin{equation}
  \label{eq:constanttermClassicalField}
f_{(-1)}=f_{\{\frac1\tau\}}=  \Tr(C(f)\frac1\tau)\simeq f.
\end{equation}

\section{Classical Fields and Affinization}
\label{sec:ClassFieldsAffin}

If $M$ is an $H_T$-module, its \emph{affinization} is the $H_T$-module
$LM=M\otimes \K$,where the tensor product is given the $H_T$-module
structure by using the coproduct: for $h\in H_T$, $m\in M$ and $k\in
\K$,
\[
h(m\otimes k)=\sum (h^{\prime}m)\otimes (h^{\prime\prime}k).
\]
Recall the counit $\epsilon\colon H_T\to \mathbb C$. This is an algebra
homomorphism, so that 
\[
\mathfrak m=\epsilon\inv(0)\subset H_T
\]
is an ideal (called the \emph{augmentation ideal} ); it is the ideal generated
by $\Delta$, or equivalently the ideal generated by
$\overline{\Delta}$. Define for any $H_T$-module $M$ the \emph{module
  of coinvariants} 
\begin{equation}
  \label{eq:DefMm}
M_{\mathfrak m}=M/\mathfrak m M.
\end{equation}
So $M_{\mathfrak m}$ is the quotient by total differences. Note that $M_{\mathfrak m}$ is a
trivial $H_T$-module. Therefore we will denote the canonical
projection by 
\begin{equation}
\Tr_M\colon M\to M_{\mathfrak m}.\label{eq:DefTraceM}
\end{equation}
Often we will suppress the subscript $M$ and write just $\Tr$.
Denote by $\mathcal LM$ the quotient $LM/\mathfrak{m}LM$. Elements of
$\mathcal{L}M$ will be written as
\[
f_{\{p\}}=\Tr(f\otimes p),\quad f\in M,\quad p\in\K.
\]
We have an action of $\hat H_T$ on $\mathcal{L}M$: if $h\in\hat H_T$
we put
\begin{equation}
h.f_{\{p\}}=f_{\{S(h)p\}}.\label{eq:defHactionLM}
\end{equation}
In this section we consider $M=V$, where $V$ is a commutative
$H_T$-Leibniz algebra. Consider $f\otimes p\in LV$. It defines a
function on $V(K)$ via
\[
f\otimes p [\hat \sigma]=\Tr_K\left(\hat\sigma(f)p(\tau)\right).
\]
Suppose that $f\otimes p\in\mathfrak{m}LV$, so that
\[
f\otimes p=\Delta\sum f_i\otimes k_i=\sum \Delta^\prime
f_i\otimes \Delta^{\prime\prime}k_i.
\]
Then, for any $\hat\sigma\in V(\K)$,
\begin{align*}
  f\otimes p[\hat \sigma]&=\Tr(\hat\sigma(f)p(\tau)=
  \Tr\left(\sum \sigma(\Delta^\prime
    f_i)\Delta^{\prime\prime}k_i\right)\\
  &=\Tr\left(\Delta(\sum\sigma(f_i)k_i)\right)=0.
\end{align*}
Hence we see that in fact any representative in $LV$ of the element
$f_{\{p\}}\in\mathcal{L}V$ will define the same function on $V(K)$.
Also we note that $f_{\{\tau(n)\}}$ defines the same function on
$V(K)$ as the coefficient $f_{(n)}$ of the classical field $C(f)$, see
(\ref{eq:Deff_(n)ClassField}). More generally $f_{\{p\}}$ defines the
same function as the value $C(f)(p)$ of the classical field, where we
think of the classical filed as a distribution on $K$ with values in
functions on $V(K)$.  We will from now on identify the coefficients of
$C(f)$ with elements of $\mathcal{L}(V)$ so that classical fields will
be generating sequences for elements of $\mathcal LV$.
\section{Multivariable Classical Fields}
\label{sec:MultClassFields}

A classical field $C(f)$ is a distribution on $K$, so we can multiply
it by the Dirac distribution, in two ways, see Section
\ref{sec:Propdist}. More generally, we can consider classical field
valued difference operators $P=\sum_{n\in\mathbb
  Z}P_n(\tau)\Delta[n]$, where $P_n(\tau)\in V(\tau)$, with $V(\tau)$
the space of classical fields. Then we get products
$P_1\delta(\tau_1,\tau_2)$, $P_2\delta(\tau_1,\tau_2)$ and we have
\begin{equation}
\label{eq:defadjoint}
P_1\delta(\tau_1,\tau_2)=P_2^*\delta(\tau_1,\tau_2),
\end{equation}
using the involution (\ref{eq:defadjoinedoperator}). 
Let $V(\tau_1,\tau_2)$ be the $V(\tau)H_T\otimes V(\tau)H_T$-module generated by
$\delta(\tau_1,\tau_2)$ (cf., section
\ref{sec:RatDist}). We think of elements of
$V(\tau_1,\tau_2)$ as classical fields depending on two variables.
Similarly we define the module of three variable classical fields
$V(\tau_1,\tau_2,\tau_3)$ as the $VH_T^{\otimes 3}$ module generated
by
\begin{align*}
  \delta(\tau_1,\tau_2,\tau_3)&=\delta(\tau_1,\tau_2)\delta(\tau_2,\tau_3)\\
                              &=\delta(\tau_1,\tau_3)\delta(\tau_2,\tau_3).
\end{align*}
We have then relations between the various $VH_T$ actions similar to
\eqref{eq:defadjoint}, for instance
\begin{align*}
P_1\delta(\tau_1,\tau_2,\tau_3)&=(P^*)_2\delta(\tau_1,\tau_2,\tau_3)\\
                               &=(P^*)_3\delta(\tau_1,\tau_2,\tau_3).  
\end{align*}

\section{Derivations and Poisson Structures}
\label{sec:DerPois}
Let $A$ be a commutative associative unital algebra, $M$ an
$A$-module. Denote by $\Der(A,M)$ the space of all derivations from
$A$ to $M$; this are linear maps $X\colon A\to M$ such that $X(ab)=aX(b)+bX(a)$.

There is an $A$-module $\Omega^1(A)$ together with a universal derivation
$d\colon A\to \Omega^1(A)$, such that for any $X\in\Der(A,M)$ there
exists a unique $A$-module morphism $\gamma_X\colon \Omega^1(A)\to M$ making
the following diagram commute:
\begin{equation}
  \label{eq:Universalitydif}
\begin{diagram}
   &        & A               &               &              \\
   &\ldTo<X &                 & \rdTo>d       &              \\
M  &        & \lTo_{\gamma_X} &               &\Omega^1(A)
\end{diagram}
\end{equation}
The \emph{module of differentials} $\Omega^1(A)$ is generated by
elements $da$, $a\in A$, subject to the relations $d(ab)=adb+bda$, for
$a,b\in A$. There is a pairing 
\[
\langle\,,\rangle\colon   \Omega^1(A)\otimes \Der(A,M)\to M,
\]
given by $\omega\otimes X\mapsto \langle
\omega,X\rangle=\gamma_X(\omega)$. In particular $\langle da,
X\rangle=X(a)$.

The space $\Der(A,M)$ is an $A$-module via 
\[
a.X(b)=a\left(X(b)\right),\quad a,b\in A.
\]
Then consider derivations $\Gamma\colon A\to \Der(A,M)$. Define the bracket
associated to $\Gamma$ by
\[
\{a,b\}=\Gamma(a)(b),\quad a,b\in A.
\]
Then we call $\Gamma$ \emph{Hamiltionian} if
\begin{itemize}
\item $\Gamma$ is skew-symmetric: $\{a,b\}=-\{b,a\}$,
\item $\Gamma$ is a Lie map: $\Gamma(\{a,b\})=[\Gamma(a),\Gamma(b)]$.
\end{itemize}

In this case we call the bracket of $\Gamma$ the Poisson bracket; $(A,\{\})$ is
then a Lie algebra. The bracket factors through the differentials:
\[
\{a,b\}=\langle \gamma_\Gamma(da),db\rangle.
\]
In case $A$ is a polynomial algebra we get simple explicit formulae
for the Poisson bracket in terms of generators: let $A=\mathbb
C[v_\alpha]_{i\in I}$. Then
\begin{align*}
  \{a,b\}&=\sum_{i,j\in I}\frac{\partial a}{\partial
    v_i}\frac{\partial b}{\partial v_j} \{v_i,v_j\},
\end{align*}
so that the Poisson bracket is uniquely determined by the Poisson
bracket between the generators.

\section{Evolutionary Derivations and Singular Poisson Brackets}
\label{sec:EvDerFieldPoissonBracket}

Let $V$ be an $H_T$-Leibniz algebra. The module of differentials
$\Omega^1(V)$ for any algebra $V$ a $V$-module, but it is in this case
in fact a $VH_T$-module, via
\[
h(fdg)=\sum(h^\prime f)d(h^{\prime\prime}g), \quad h\in H_T, f,g\in V.
\]
We call, for $M$ an $VH_T$-module, a derivation $X\colon V\to M$
\emph{evolutionary} in case
\[
X(h f)=hX(f), \quad h\in H_T,f\in{V}.
\]
In this case the $V$-module morphism $\gamma_X\colon \Omega^1(V)\to M$ of
\eqref{eq:Universalitydif} is in
fact a morphism of $VH_T$-modules. 

Now we would like to generalize the notion of a Hamiltonian map to the
case of evolutionary derivations. An obstacle is that the space
$\Derev(V,V)$ of evolutionary derivations (from $V$ to itself) is not
(in an obvious way) a module over $V$ or $VH_T$: if $X$ is
evolutionary then
\[
f.X(hg)=f.(hX(g))\ne h f.X(g),
\]
as $h\in H,f\in V$ do not commute in $VH_T$, in general.
In case $V$ is free as $H_T$-module there is a $V$-module structure on
$\Der(V,V)$, using a basis for $V$, see section
\eqref{eq:paringDerevForms} below.

Therefore we consider instead the space $\mathcal
D_2=\Derev(V,V(\tau_1,\tau_2))$ of evolutionary derivations from $V$
to the module $V(\tau_1,\tau_2)$ of two variable classical
fields introduced in section \ref{sec:MultClassFields}.  So an
element $X\in \mathcal{D}_2$ acts on $f\in V$ by
\[
X(f)=\sum_{n\in\mathbb Z} X(f)_n(\tau_2)
\Delta[n]_2\delta(\tau_1,\tau_2),\quad X(f)_n\in V,
\]
and satisfies
\begin{equation}
\begin{split}
  X(fg)&=f(\tau_2)X(g)+g(\tau_2)X(f), \quad f,g\in V,\\
  X(hf)&=h_2X(f)     ,\quad h\in H_T.
\end{split}\label{eq:DefEvolder}
\end{equation}
In other words, we use the $1\otimes VH_T$-module structure on the
target to define the words evolutionary and derivation.

Now note that we can define a $VH_T$ module structure on
$\mathcal{D}_2$ by using the other $VH_T$ module structure on the
target: if $X$ is an evolutionary derivation then define for $P\in VH_T$
\[
P.X(f)=P_1(X(f)),\quad  f\in V.
\]
The two actions of $VH_T$ on $V(\tau_1,\tau_2)$ commute:
\[
P_1P_2 X(f)=P_2P_1X(f),
\]
from which we see that $P.X$ is again evolutionary:
\[
P.X(h f)=P_1h_2X(F)=h_2 P.X(f).
\]
Any $X\in\mathcal{D}_2$ factors through $\Omega^1(V)$:
\begin{equation}
\begin{diagram}
     V           &   \rTo^d         & \Omega^1(V)    \\
  \dTo<X         & \ldTo_{\gamma_X} &                \\
V(\tau_1,\tau_2) &                  &                  
\end{diagram}
\end{equation}
where $\gamma_T\colon \Omega^1(V)\to V(\tau_1,\tau_2)$ is a morphism of
$VH_T$-modules (for the $1\otimes VH_T$ structure on the target).

Instead of considering derivations with target $V(\tau_1,\tau_2)$ we
could use distributions in any pair of variables as target. In this
case we will write
\[
X(f)_{ij}=\sum_n X(f)_n(\tau_j) \Delta_j(n)\delta(\tau_i,\tau_j).
\]
We can extend any $X\in\Derev(V,V(\tau_i,\tau_k))$ to a map
\[
X_{ijk}\colon  V(\tau_j,\tau_k)\to V(\tau_i,\tau_j,\tau_k),
\]
by
\begin{multline*}
  \sum_{n\in\mathbb Z}p_n(\tau_k)\Delta_k(n)\delta(\tau_j,\tau_k)\mapsto
  \sum_{n\in\mathbb Z}X(p_n)_{ik}\Delta_k(n)\delta(\tau_j,\tau_k)=\\
  \sum_{n\in\mathbb Z}X(p_n)_{m}(\tau_k)\Delta_k(m)\delta(\tau_i,\tau_k)\Delta_k(n)\delta(\tau_j,\tau_k).
\end{multline*}
Then $X_{ijk}$ is an evolutionary derivation, for the
$(VH_T)_k$-module structures on $V(\tau_j,\tau_k)$ and $ V(\tau_i,\tau_j,\tau_k)$.

Now we will define the evolutionary analog of the Hamiltonian maps introduced in
section  \ref{sec:DerPois}. Consider a map
\[
\Gamma\colon V\to \Derev(V,V(\tau_1,\tau_2)),
\]
and introduce the \emph{Poisson bracket} of $f,g\in V$ by
\begin{equation}
\{f(\tau_1),g(\tau_2)\}=\Gamma(f)(g)\in
V(\tau_1,\tau_2).\label{eq:PoissonBracketfg}
\end{equation}
\begin{defn}\label{def:SingHamMap}
We call $\Gamma$  a \emph{singular Hamiltonian map} in case
\begin{Myenumerate}
\item \label{item:7}$\Gamma$ is an evolutionary derivation,
\item \label{item:8}$\Gamma$ is
  skew-symmetric: $\{f(\tau_1),g(\tau_2)\}=-\{g(\tau_2),f(\tau_1)\}$,
\item \label{item:9}$\Gamma$ is a Lie map:
  \begin{multline*}
    \{f(\tau_1),\{g(\tau_2),k(\tau_3)\}\}=\{\{f(\tau_1,g(\tau_2)\},k(\tau_3)\}\\
    + \{g(\tau_2),\{f(\tau_1),k(\tau_3)\}\}.
  \end{multline*}
\end{Myenumerate}
\end{defn}
Here we extend in \eqref{item:9} the derivation $\{f(\tau_1),
\cdot\}=\Gamma(f)$ from a map on $V$ to a map
$\Gamma(f)_{123}\colon V(\tau_2,\tau_3)\to V(\tau_1,\tau_2,\tau_3)$, as
discussed above, and similar for the other terms.

We call, in case $\Gamma$ is a singular Hamiltonian map for a
$H_T$-Leibniz algebra $V$, the pair $(V,\Gamma)$ an $H_T$-Poisson algebra.

The singular Poisson bracket $\{f(\tau_1),g(\tau_2)\}$ is an
evolutionary derivation in both $f$ and $g$ and so factors through
$\Omega^1(V)^{\otimes 2}$. We can write
\[
\{f(\tau_1),g(\tau_2)\}=\langle \gamma_\Gamma(df),dg\rangle,
\]
where $\langle\,,\,\rangle\colon \mathcal{D}_2\otimes \Omega^1(V)\to
V(\tau_1,\tau_2)$ is the canonical pairing, and $\gamma_\Gamma$ is the unique $VH_T$-module
morphism $\Omega^1(V)\to \mathcal{D}_2$ making the diagram
\begin{equation}  
\begin{diagram}
                           &              & V                    &               &              \\
                           &\ldTo<\Gamma  &                      & \rdTo>d       &              \\
\mathcal D_2 &              & \lTo_{\gamma_\Gamma}& &\Omega^1(V)\\
\end{diagram}
\end{equation}
 commute.

\section{Projection from $V(\tau_1,\tau_2)$ to $V$}
\label{sec:projection}

Let   $F(\tau_1,\tau_2)\in
V(\tau_1,\tau_2)$ be a two variable classical field, see section
\ref{sec:MultClassFields}. It can be 
can be uniquely expanded as
\[
F(\tau_1,\tau_2)=\sum_{k\ge0}
F_k(\tau_2)\Delta_2[k]\delta(\tau_1,\tau_2),
\]
see (\ref{eq:expansiondelta}). Define a map
\[
\Pi\colon  V(\tau_1,\tau_2)\to V, \quad F(\tau_1,\tau_2)\mapsto (F_0)_{(-1)}
\]
so that 
\[
\Pi(F(\tau_1,\tau_2))=\Tr_{\tau_1,\tau_2}(F(\tau_1,\tau_2)\frac1{\tau_2}),
\]
see Lemma \ref{lem:coefficientDeltaexpansion} and Section \ref{sec:ExpClassFields}.
\begin{lem} \label{LemmaPiVH_T}
  Let $P, Q\in VH_T$, so that
\[
P=\sum_{n\in\mathbb Z}p_n \Delta[n], \quad Q=\sum_{n\in\mathbb Z}q_n
\Delta[n].
\]
Then
\begin{Myenumerate}
\item $\Pi(P_2 Q_2\delta(\tau_1,\tau_2))=P(q_0)$.\label{item:3}
\item $\Pi(P_1 Q_2\delta(\tau_1,\tau_2))=Q(P_V)$.\label{item:4} 
\end{Myenumerate}
\end{lem}

Here we write for instance $P_i$ for the classical field valued
difference operator $P_i=\sum_{n\in\mathbb Z}p_n(\tau_i) \Delta_i[n]$,
for $i=1,2$.

In particular, it follows from part  \eqref{item:4} that
\begin{equation}
  \label{eq:PiDelta_1vanishes}
  \Pi(\Delta_1F(\tau_1,\tau_2))=0, \quad \forall F(\tau_1,\tau_2)\in V(\tau_1,\tau_2).
\end{equation}
\begin{proof}
  We use two facts about $H_T$. First of all the product of elements
  of the difference operator basis has an expansion
  \begin{equation}
    \label{eq:productdifop}
    \Delta[n]\Delta[m]=\sum_{k\in\mathbb Z} a_{nm}^k \Delta[k], \quad
    a_{nm}^k\in\mathbb Z,
  \end{equation}
  with the coefficient $ a_{nm}^0$ of $\Delta[0]=1$ zero, unless
  $n=m=0$. Secondly the coproduct of $\Delta[n]$ has the form
  \begin{equation}
    \label{eq:coproductform}
    \pi(\Delta[m])=\Delta[m]\otimes 1 +\sum_i h_i^\prime\otimes h_i^{\prime\prime},
  \end{equation}
  with $h_i^{\prime\prime}\in\mathfrak m$.
  It suffices to check part \eqref{item:3} in case $P=p\Delta[m]$,
  $Q=q\Delta[n]$, and in this case we have
  \begin{align*}
    \Pi&\left(P_2Q_2\delta(\tau_1,\tau_2)\right)= \Pi\left(
      p(\tau_2)\Delta_2(n)[q(\tau_2)\Delta_2(m)\delta(\tau_1,\tau_2)]\right)\\
    &=\Pi\left(p(\tau_2)\left[\Delta_2(n)(q(\tau_2))\Delta_2(m)+ \sum_i
    h_{i,2}^\prime(q(\tau_2))h_{i,2}^{\prime\prime}\right]\delta(\tau_1,\tau_2)\right)\\
    &=\delta_{m0}p\Delta[n](q)\\
    &=P(Q_0).
  \end{align*}
  Then part \eqref{item:4} follows from \eqref{item:3}: we have 
  \[
  \Pi\left(P_1Q_2\delta(\tau_1,\tau_2)\right)=\Pi(Q_2
  (P^S)_2\delta(\tau_1,\tau_2))=Q(P^S_0),
  \]
  but 
  \[
  P^S_0=\sum S\Delta[n](p_n)=P_V.
  \]
\end{proof}
\begin{lem}
Let $X$ be an evolutionary derivation $X\colon V\to V(\tau_1,\tau_2)$.
Then the composition
\[
X_0=\Pi\circ X\colon  V\to V
\]
is an evolutionary derivation of $V$.  
\end{lem}
\begin{proof}
We have for $f,g\in V$ by definition of evolutionary derivation
with values in $V(\tau_1\tau_2)$, see (\ref{eq:DefEvolder}),
\[
  X_0(fg)=\Pi\bigl( f(\tau_2)X(g)+ g(\tau_2)X(f)\bigr)=f  X_0(g)+g  X_0(f),
\]
using Part \ref{item:3} of Lemma \ref{LemmaPiVH_T} with $P=f$, $Q=1$. In case $h\in{
  H_T}$ similarly with $P=h$ and $Q=1$
\[
  X_0(h.f)=\Pi( h_2X(f)=h  X_0(f).
\]
\end{proof}
So we obtain a map
\begin{equation}
\Pi^{\text{ev}}\colon  \Derev\left(V,V(\tau_1,\tau_2)\right)\to \Derev(V,V),\quad
X\mapsto X_0.\label{eq:projEvDer}
\end{equation}
Now let $\Gamma\colon  V\to \Derev(V, V(\tau_1,\tau_2)$ be a singular
Hamiltonian structure, see section
\ref{sec:EvDerFieldPoissonBracket}. Recall that if $P\in VH_T$ then 
$\Gamma(P.f)(g)= P_1 \Gamma(f)(g)$. We define a map by composition
with $\Pi^{\text{ev}}$:
\[
\Gamma_0=\Pi^{\text{ev}}\circ \Gamma\colon V\to \Derev(V,V), 
\quad \Gamma_0(f)=\left(\Gamma(f)\right)_0.
\]
Recall that $\Gamma$ factors through the differentials: if $f\in V$ then
\[
\Gamma(f)=\gamma_\Gamma(df),
\]
where $\gamma_\Gamma\colon \Omega^1(V)\to \Derev(V,V(\tau_1,\tau_2)$ is a
morphism of $VH_T$-modules. In particular, if $f\in\mathfrak m V$
(i.e., $f$ is a total difference, $f=\Delta g$), then
\[
\Gamma_0(f)=\Pi\circ \gamma_\Gamma(df)=\Pi\circ \gamma_\Gamma(\Delta
dg)=\Pi\left(\Delta_1\gamma_\Gamma(dg)\right)=0,
\]
by \eqref{eq:PiDelta_1vanishes}. This means that $\Gamma_0$ factors
through, $\Omega^1(V)_{\mathfrak m}$, the module of coinvariants.
Then we obtain the following commutative diagram
\begin{equation}
\begin{diagram}
     V           &          &  \rTo^d  &                       & \Omega^1(V)    \\
                 &          &          &\ldTo(4,4)_{\Pi\circ\gamma_X} &         \\
  \dTo^{\Gamma_0}&          &          &                       &\dTo_{\Tr_{\Omega^1(V)}}\\
                 &          &          &                       &\\
\Derev(V,V)      &          &   \lTo_B &                       &\Omega^1(V)_{\mathfrak m}
\end{diagram}
\end{equation}
We write $\delta f=\Tr_{\Omega^1(V)}(df)$ for the map
$V\to\Omega^1(V)_{\mathfrak m}$ (variational differential). For $V$ a
free $H_T$-Leibniz algebra $B\colon \Omega^1(V)_{\mathfrak m}\to
\Derev(V,V)$ is a matrix of difference operators, i.e., an matrix with
coefficients in $VH_T$, as we will explain in the next section.

\section{Free Difference Algebras}

Let $V$ be the free difference algebra generated by $v_\alpha,
\alpha\in\mathcal A=\{1, 2, \dots, k\}$. So $V$ is the polynomial
algebra 
\[
V=\mathbb C[v_\alpha^{(n)};\alpha\in\mathcal A,n\in\mathbb Z],
\] 
where $v_\alpha^{(0)}=v_\alpha$ and the action of $H_T$ is determined
by 
\[
\Delta[n]v_\alpha=v_\alpha^{(n)},
\]
together with linearity and the Leibniz rule for $H_T$. In this case
$\Derev(V,V)$ and $\Omega^1(V)_{\mathfrak m}$ are free $V$-modules of
rank $k=\lvert \mathcal A\rvert$, which simplifies the discription of
singular Hamiltonian structures on $V$.

Any $X\in \Derev(V,V)$ is uniquely determined by its values
$X_\alpha=X(v_\alpha)$ on the generators of $V$. Introduce
notation
\[
X_\alpha\frac{\delta{\ }}{\delta v_\alpha}=\sum_{n\in\mathbb Z}
\Delta[n](X_\alpha)\frac {\partial}{\partial
  v_\alpha^{(n)}}\in\Derev(V,V),
\]
so that
\[
X=\sum_{\alpha\in\mathcal A}X_\alpha \frac{\delta}{\delta v_\alpha}.
\]
Define an action of $V$ on evolutionary derivations by 
\[
f.X=\sum_{\alpha\in\mathcal A} (fX_\alpha)\frac{\delta}{\delta
  v_\alpha}, \quad f\in V.
\]
Then $\Derev(V,V)$ is a free $V$-module:
\[
\Derev(V,V)=\bigoplus_{\alpha\in\mathcal A} V \frac{\delta}{\delta
  v_\alpha},
\]
generated by the (evolutionary) derivations $\frac{\delta}{\delta v_\alpha}$.

The module of differentials of $V$ is a free $V$-module:
\[
\Omega^1(V)=\bigoplus_{\alpha\in\mathcal A, n\in\mathbb Z} V
dv_\alpha^{(n)}.
\]
Now the $H_T$-module structure on $\Omega^1(V)$ is such that 
\[
dv_\alpha^{(n)}=\Delta[n]dv_\alpha,
\]
so that $\Omega^1(V)$ is also a free $VH_T$-module:
\[
\Omega^1(V)= \bigoplus_{\alpha\in\mathcal A} VH_T dv_\alpha.
\]
In particular for $f\in V$ we can express the differential in terms of
difference operators:
\[
df= \sum_{\alpha\in\mathcal A}d_\alpha f dv_\alpha,
\]
where the $d_\alpha f$ are difference  operators 
\[
d_\alpha f=\sum_{n\in\mathbb Z}\frac{\partial f}{\partial
  v_\alpha^{(n)}}\Delta[n]\in V H_T.
\]
Inside $\Omega^1(V)$ we have a $V$-submodule 
\[
\Omega^1(V)_0=\bigoplus_{\alpha\in\mathcal A} V dv_\alpha.
\] 
We have a pairing 
\begin{equation}
\langle \,, \rangle\colon  \Derev (V,V)\otimes \Omega^1(V)_0\to V, \quad
X\otimes \omega\mapsto \sum X_\alpha\omega_\alpha,
\label{eq:paringDerevForms}
\end{equation}
when $X=\sum X_\alpha\frac{\delta}{\delta v_\alpha}$, $\omega=\sum
\omega_\alpha dv_\alpha$.
It is known (see \cite{MR86m:58070}, chapter 2) that
\begin{equation}
  \label{eq:intersectMOmega_0}
  \Omega^1(V)_0\cap \mathfrak m\Omega^1(V)=\{0\}.
\end{equation}
Recall the module of coinvariants $\Omega^1(V)_{\mathfrak
  m}=\Omega^1(V)/\mathfrak m \Omega^1(V)$. We have for all $P\in VH_T$
\[
Pdv_\alpha\equiv P_V dv_\alpha \mod \mathfrak m \Omega^1(V),
\]
(where $P_V$ is defined in \eqref{eq:defP_V}).  This implies, combined
with \eqref{eq:intersectMOmega_0}, that
\[
\Omega^1(V)=\mathfrak m \Omega^1(v)\bigoplus \Omega^1(V)_0.
\]
So we have a splitting of the exact sequence
\[
\begin{diagram}
  0&\rTo &\mathfrak m\Omega^1(V)&\rTo& \Omega^1(v)&\pile{\lTo^s\\
    \rTo_{\mathsf p}}&\Omega^1(V)_{\mathfrak m}&\rTo & 0,
\end{diagram}
\]
where the splitting map $s\colon  \Omega^1(v)_{\mathfrak m} \to \Omega^1(V)$
is 
\[
[\sum
P_\alpha dv_\alpha]\mapsto \sum (P_\alpha)_V dv_\alpha.
\]
In particular we obtain a variational differential
$\delta=s\circ\mathsf p\circ d\colon  V\to \Omega^1(V)_0$ given by
\[
\delta f=\sum_\alpha \frac{\delta f}{\delta v_\alpha} dv_\alpha, 
\]
where the variational derivative is
\begin{equation}
\frac{\delta f}{\delta v_\alpha}= (d_\alpha f)_V=\sum_n S\Delta[n]\left( \frac{\partial
  f}{\partial v_\alpha^{(n)}}\right).\label{eq:defvarderivative}
\end{equation}
The variational derivative $\delta f$ is the unique element of
$\Omega^1(V)_0$ such that for all $X\in \Derev$
\[
X(f)\equiv \langle X, \delta f\rangle \mod \mathfrak m V,
\]
using the pairing \eqref{eq:paringDerevForms}. Indeed, if $\delta=\sum
\delta_\alpha dv_\alpha\in \Omega^1(V)_0$ also 
satisfies this then we have for all $X\in \Derev$ $\langle X,\delta f-\delta\rangle\equiv 0
\mod \mathfrak m V$, or for each $\alpha$
\[
X_\alpha(\frac{\delta f}{\delta v_\alpha}-\delta_\alpha)\in \mathfrak m V,
\quad \forall X_\alpha  \in  V,
\]
but it is known that this forces $\frac{\delta f}{\delta
  v_\alpha}=\delta_\alpha$, see \cite{MR86m:58070}, Lemma 17.

The module of coinvariants is a free $V$-module:
\[
\Omega^1(V)_{\mathfrak m}=\bigoplus_\alpha V dv_\alpha.
\]
We can therefore represent an element $\omega=\sum_\alpha
\omega_\alpha dv_\alpha$ of $\Omega^1(V)_{\mathfrak m}$ as a column
vector
$\omega=
\begin{pmatrix}
  \omega_1\\ \vdots \\ \omega_k
\end{pmatrix}
$ and similarly if $X=\sum_\alpha X_\alpha \frac {\delta}{\delta
  v_\alpha}\in \Derev(V,V)$ then we write this as a column vector
$X=
\begin{pmatrix}
  X_1\\ \vdots \\ X_k
\end{pmatrix}
$.
Now suppose we are given a singular Hamiltonian 
structure$\Gamma\colon 
V\to \Derev(V,V(\tau_1,\tau_2)$ on $V$. We
have seen (see section \ref{sec:projection}) that this gives rise 
to a map $B\colon \Omega^1(V)_{\mathfrak m}\to \Derev(V,V)$. In our case $B$ is a map
between free rank $k$ $V$-modules, and we will now see that $B$ is given
by an $k\times k$ matrix of difference operators. 

Introduce difference operators $B_{\beta\alpha}\in VH_T$ by
\[
\Gamma(v_\alpha)(v_\beta)=B_{\beta\alpha,2}\delta(\tau_1,\tau_2)\in
V(\tau_1,\tau_2).
\]
Let $X_f=\Pi\circ \Gamma(f)$ be the evolutionary derivation on $V$
associated to $f$. The we have $X_f=\sum_\alpha
X_{f,\alpha}\frac{\delta}{\delta v_\alpha}$ and
\begin{align*}
  X_{f,\alpha}&=X_f(v_\alpha)=\Pi\left(\Gamma(f)(v_\alpha)\right)\\
              &=\Pi\left( (df_\beta)_1
                \Gamma(v_\beta)(v_\alpha)\right)\\
              &=\Pi\left( (df_\beta)_1
                B_{\alpha\beta,2}\delta(\tau_1,\tau_2)\right)\\
              &=\sum_\beta B_{\alpha\beta}\left(\frac{\delta f}{\delta
                v_\beta}\right),
\end{align*}
by Lemma \ref{LemmaPiVH_T} and \eqref{eq:defvarderivative}. So
\[
X_f
\begin{pmatrix}
  v_1\\ \vdots\\ v_k
\end{pmatrix}= \mathcal B 
\begin{pmatrix}
  \frac{\delta f}{\delta v_1}\\\vdots\\  \frac{\delta f}{\delta v_k}
\end{pmatrix}, \quad \mathcal B_{\alpha\beta}=B_{\alpha\beta},
\]
so that $\mathcal B$ is the matrix of the map $B$.

\chapter[Hamiltonian Structure for Infinite Toda][Infinite Toda]{Hamiltonian Structure for Infinite Toda}
\label{chap:HamstrucInfToda}

\section{The Finite Toda Lattice}
\label{sec:finTodaLattice}

Recall, see \cite{MR971987}, that the finite non periodic Toda lattice
with $N$ particles is a Hamiltonian system on a phase space with
global coordinates $p_i, q_i$ for $i=1,\dots, N$ and Hamiltonian
$H=\frac12\sum_{i=1}^Np_i^2 + \sum_{k=2}^N \exp(q_{k-1}-q_k)$. So the
equations of motion are
\begin{align*}
  \dot p_i&=\exp(q_{i-1}-q_i)-\exp(q_{i}-q_{i+1}),\\
  \dot q_i&=p_n.
\end{align*}
Introduce new variables $  B_i=p_i$, $ i=1,\dots,N$ and $C_j=\exp(q_{j-1}-q_j)$, $j=2,\dots,N$.
Then we obtain the following system in the new variables:
\begin{equation}
\begin{split}
  \dot B_i&=C_i-C_{i+1}       \qquad\quad  i=1,\dots N, C_{N+1}=0\\
  \dot C_j&=C_j(B_{j-1}-B_{j}) \quad  j=2,\dots N.
\end{split}\label{eq:BCequation}
\end{equation}
Introduce a Lax matrix
\begin{equation}
L=
\begin{pmatrix}
  B_1 & 1     & 0      &0             &0         &0\\
  C_2 & B_2   &  1     &\cdots        &\cdots    &0\\
   0  & C_3   &  B_3   &\cdots        &\cdots    &0\\
\vdots&\vdots & \vdots &\ddots   &\vdots    &\vdots\\
\vdots&\vdots &\vdots  &C_{N-1}  & B_{N-1}  &1  \\
0     &0      &  0     &0        & C_N      &B_N  
\end{pmatrix}.\label{eq:LaxMatrix}
\end{equation}
Then the system (\ref{eq:BCequation}) is equivalent to the Lax
equation
\begin{equation}
\dot L=[L_-, L],\label{eq:LaxEquation}
\end{equation}
where $L_-$ is the strictly lower triangular part of $L$. We also have
for the Hamiltonian
\begin{equation}
H=\frac12 \Tr(L^2)=\frac12\sum_{i=1}^NB_i^2+\sum_{j=2}^N C_j.\label{eq:HamiltonianTrace}
\end{equation}
Now the system (\ref{eq:BCequation}) make still sense if we let $i,j$
run over all integers, and we can introduce a doubly infinite Lax
matrix such that this infinite system is stil given by
(\ref{eq:LaxEquation}). Now, however, the meaning of
(\ref{eq:HamiltonianTrace}) becomes less clear (but see
\cite{MR1304086}). Kuperschmidt, for instance, advocates in
\cite{MR86m:58070} an other way of dealing with infinite systems,
where one replaces the infinite Lax matrix by a difference operator.
We follow in this chapter Kuperschmidt's approach, expressed in the
formalism of singular Hamiltonian maps, and also indicating the
connection with the formalism of $r$-matrices.

\section{Lax Operator and $r$-Matrix}
\label{sec:LaxOp}

Consider the difference operator (``Lax Operator'')
\[
L=A T + B + CT\inv\in VH_T,
\]
where $V$ is the $H_T$-Leibniz algebra freely generated by
$A,B,C$. So $V$ is the polynomial algebra
\[
V=\mathbb C[A^{(n)},B^{(n)},C^{(n)}],\quad Y^{(n)}=\Delta[n]Y,
\quad Y=A,B,C.
\]
We will also write the Lax operator as $L=\sum_\alpha v_\alpha
T^\alpha$, so that $v_1=A,v_0=B,v_{-1}=C$.

Any evolutionary derivation $X$ of $V$ defines a time-flow on the Lax
operator by 
\[
\dot L=X_L=X(A)T+X(B)+X(C)T\inv. 
\]
Conversely we will use $L$ to define singular Hamiltonian structures
on $V$.

We have a map, if $P=\sum_\alpha p_\alpha T^\alpha$,
\[
\Sp \colon VH_T\to V, \quad P\to p_{0}.
\]
For all $P,Q\in VH_T$ we have
\begin{equation}
  \label{eq:traceproperty}
  \Sp([P,Q])\equiv 0 \mod \mathfrak{m}V.
\end{equation}
Indeed, it suffices to check this for $P=pT^\alpha,Q=qT^{-\alpha}$,
for $\alpha\in \mathbb Z$, $p,q\in V$. But in that case
\[
[pT^\alpha,qT^{-\alpha}]=(T^\alpha-1)(qT^{-\alpha}\langle p\rangle),
\]
and $T^\alpha-1\in\mathfrak m$.

Recall the pairing \eqref{eq:paringDerevForms}. In our situation we
have encoded the evolutionary derivation $X$ in the difference
operator $X_L$. Similarly we can encode an element $\omega=\sum
\omega_\alpha dv_\alpha\in \Omega^1(V)_0$ in a difference operator
 $\omega_L=
\sum T^{-\alpha}\omega_\alpha$. Then we have
\[
\langle X,\omega\rangle=\Sp(X_L\omega_L).
\]
For $f\in V$ we get a difference operator $\delta_Lf=\sum
T^{-\alpha}\frac{\delta f}{\delta v_\alpha}$ such that for
all $X\in\Derev$
\begin{equation}
X(f)\equiv \Sp(X_L \delta_Lf) \mod \mathfrak m
V.\label{eq:defdelta_Lf}
\end{equation}
Note that the difference operator $\delta_L f$ is \emph{not} uniquely
determined by this condition: we are free to add to $\delta_L f$ any
difference operator of the form $\delta_{\ge 2}+\delta_{\le -2}$,
where $\delta_{\ge 2}=\sum_{\alpha\ge 2}T^\alpha \delta_\alpha$ and
$\delta_{\le -2}= \sum_{\alpha\le -2}T^\alpha \delta_\alpha$, since
$\delta_{\ge 2}$ and $\delta_{\le -2}$ are orthogonal to $X_L$ with
respect to $\Sp$. 

We have projections
\[
\Pi_{\le k}, \Pi_{\ge k}, \Pi_{> k}, \Pi_{< k}, \Pi_k\colon  VH_T\to VH_T,
\]
where, if $P=\sum_{\alpha\in\mathbb Z}p_{\alpha}T^\alpha\in VH_T$, for instance
\[
\Pi_{\le k}(P)=\sum_{s\le k}p_{\alpha}T^\alpha, \quad \Pi_k P=p_{k},
\]
and similar for the other projections.

We use $\Sp$ also to define the dual of a map $\sigma\colon VH_T\to VH_T$ by
\[
\Sp(\sigma(P)Q)=\Sp(P\sigma^*(Q)).
\]
Then, for instance, we have
\[
\Pi_{\ge 0}^*=\Pi_{\le 0},\quad \Pi_{<0}^*=\Pi_{>0}.
\]
Define the $r$-matrix $\rho\colon VH_T\to VH_T$ for $VH_T$ by 
$P\mapsto (\Pi_{\ge 0}-\Pi_{<0})(P)$.

Then $\rho$ is not skew-adjoint, but we have $\rho=S+A$, where
$A=\Pi_{>0}-\Pi_{<0}$ is skew and $S=\Pi_0$ is symmetric, so that
$\rho^*=S-A$. 

Besides the Lax operator $L$ (with coefficients in $V$) we also have
the Lax operator $L(\tau_i)$ with coefficients classical fields
depending on $\tau_i$:
\[
L(\tau_i)=A(\tau_i)T_i+B(\tau_i)+C(\tau_i)T_i\inv.
\]
Define distributions associated to the $r$-matrix $\rho$ and its
adjoint $\rho^*$:
\begin{align*}
  r(\tau_1,\tau_2)&= \sum_{\alpha,m\in\mathbb Z}
  \tau_1(-m-1)T_1^\alpha \rho\left(T_2^{-\alpha}\right)\tau_2(m)
=\rho \left(T_2^{-\alpha}\right)\delta(\tau_1,\tau_2)T_1^\alpha\\
  r^*(\tau_1,\tau_2)&= \sum_{\alpha,m\in\mathbb Z} \rho\left(T_1^{-\alpha}\right)\tau_1(m)
   \tau_2(-m-1)T_2^\alpha=\rho(T_1^{-\alpha}\delta(\tau_1,\tau_2)T_2^\alpha.
\end{align*}
We define commutators of the (field-valued) Lax operator with these
distributions:
\[
[L_2, r]= [L(\tau_2),  \rho\left(T_2^{-\alpha}\right)\delta(\tau_1,\tau_2)]T_1^\alpha,
\]
and similarly
\[
[L_1, r^*]=\sum_{\alpha,m\in\mathbb Z}[L(\tau_1), \rho\left(T_1^{-\alpha}\right)\delta(\tau_1,\tau_2)]
T_2^\alpha.
\]
If $h\in H_T$, $x\in V$, we write $h\langle x\rangle$ for the action
of $h$ on $x$, and $h x$ for the multiplication in $VH_T$. For
instance, $T^k x= T^k\langle x\rangle T^k$.

\begin{lem}\label{lem:lax-operator-r-comm}
  \begin{align*}
[L_2, r]-[L_1,r^*]&=C(\tau_2)(T_2\inv-1)\langle\delta(\tau_1,\tau_2)\rangle
T_2\inv+\\&\quad+(T_2-1)\langle C(\tau_2)
\delta(\tau_1\tau_2)\rangle T_1\inv.
  \end{align*}
\end{lem}
\begin{proof}
  First of all
  \begin{align*}
[L(\tau_i),T_i^k \delta(\tau_1,\tau_2)]&=
\left (A(\tau_i)T^{k+1}\langle \delta(\tau_1,\tau_2)\rangle-
T_i^k\langle A(\tau_i)\delta(\tau_1,\tau_2)\rangle \right ) T_i^{k+1}\\
&+\left (B(\tau_i)T_i^{k}\langle \delta(\tau_1,\tau_2)\rangle-T_i^k
\langle B(\tau_i)\delta(\tau_1,\tau_2)\rangle \right ) T_i^{k}\\
&+\left (C(\tau_i)T_i^{k-1}\langle \delta(\tau_1,\tau_2)\rangle-T_i^k
\langle C(\tau_i)\delta(\tau_1,\tau_2)\rangle \right ) T_i^{k-1}.
  \end{align*}
  We will write all expressions in the form
  $\sum_{\alpha,\beta}\ell_{\alpha,\beta}T_1^\alpha T_2^\beta$, i.e.,
  with the powers of $T$ to the right.  Consider the coefficient
  of $T_1^{-k}T_2^{k+1}$ for $k>0$ in $[L_2, r]-[L_1,r^*]$.  It is
\begin{align*}
  A(\tau_2)T_2^{k+1}\langle \delta(\tau_1,\tau_2)\rangle
  -T_2^k\langle  A(\tau_2)\delta(\tau_1,\tau_2)\rangle + &\\
  A(\tau_1)T_1^{-k}\langle\delta(\tau_1,\tau_2)\rangle-
  T_1^{-k-1}\langle \delta(\tau_1,\tau_2)A(\tau_1)\rangle&=0
\end{align*}
by the property \eqref{eq:defadjoint} of the action of $VH_T\otimes
VH_T$ on two variable classical fields.  In the same manner one checks
that the coefficients of $T_1^{-k}T_2^{k},T_1^{k}T_2^{-k}$, for $k\ge
0$ and of $T_1^{-k}T_2^{k-1}, T_1^{k}T_2^{1-k}, T_1^{k}T_2^{-k-1}$ for
$k>0$ all vanish. There remain only 2 nonzero coefficients. For
$T_1^0T_2^1$ the coefficient is
\[
C(\tau_2)T_2\inv \delta(\tau_1,\tau_2)-C(\tau_2)\delta(\tau_1,\tau_2)
=C(\tau_2)(T_2\inv-1) \delta(\tau_1,\delta_2).
\]
Similarly the coefficient of $T_1\inv T_2^0$ is found to be
\begin{align*}
  \tau(-m-1)[ C(\tau_1)T_1\inv \langle \tau_1(m)\rangle&-C(\tau_1)\tau_1(m)]\\
  &= C(\tau_1)T_1\inv\delta(\tau_1,\delta_2)-C(\tau_1)\delta(\tau_1,\tau_2)\\
  &=(T_1\inv-1) C(\tau_1)\delta(\tau_1,\delta_2)\\
  &= C(\tau_2)(T_2-1)\delta(\tau_1,\delta_2).
\end{align*}
\end{proof}
Given a singular Hamiltonian structure on $V$ we define the Poisson
bracket of the field Lax operator $L(\tau)=\sum_\alpha
v_\alpha(\tau)T^\alpha$ by
\[
\{L(\tau_1)\poissontensor L(\tau_2)\}=\sum_{\alpha,\beta}\{v_\alpha(\tau_1),v_\beta(\tau_2)\}
T_1^\alpha T_2^\beta.
\]
Conversely we define on $V$ a singular Hamiltonian structure by
putting
\begin{equation}
  \label{eq:poissonLaxtensor}
  \{L(\tau_1)\poissontensor L(\tau_2)\}=[L_2, r]-[L_1,r^*].
\end{equation}
By Lemma \ref{lem:lax-operator-r-comm} this means that we have
specified the Poisson brackets of the generators of $V$ to be
\begin{align}\label{eq:fundamentalpoissonbracketsToda}
  \{B(\tau_1),C(\tau_2)\}&=C(\tau_2)(T_2\inv-1)\delta(\tau_1,\tau_2),\\
  \{C(\tau_1),B(\tau_2)\}&=(1-T_2)C(\tau_2)\delta(\tau_1,\tau_2),
\end{align}
and all other Poisson brackets are zero, in particular
$\{A(\tau_1),F(\tau_2)\}=0$, for all $F\in V$.
Then one checks that this can be extended consistently to a Poisson
brackets on all of $V$, essentially since
\begin{multline*}
  \{B(\tau_1),\{B(\tau_2),C(\tau_3)\}\}
  +\{B(\tau_2),\{C(\tau_3),B(\tau_1)\}\} +\\
  \{C(\tau_3),\{B(\tau_1),C(\tau_2)\}\} =0.
\end{multline*}
The projection \eqref{eq:projEvDer} of $\Gamma(f)\in
\Derev(V,V(\tau_1,\tau_2)$ gives an evolutionary derivation
$X(f)\in\Derev(V,V)$ which reads in column form
\[
X(f)
\begin{pmatrix}
  A\\B\\C
\end{pmatrix}=
\begin{pmatrix}
  0&0&0\\0&0&(1-T)C\\0&C(T\inv-1)&0
\end{pmatrix}
\begin{pmatrix}
  \frac{\delta f}{\delta A}\\  \frac{\delta f}{\delta B}\\  \frac{\delta f}{\delta C}
\end{pmatrix}
\]
In terms of the Lax operator we have
\[
X_L(f)=0T+ (1-T)\langle C\frac{\delta f}{\delta c}\rangle T^0 + C(T-1)\langle\frac{\delta
  f}{\delta b}\rangle T\inv.
\]
We can express this in terms of the difference operator $\delta_Lf$
and the classical $r$-matrix.

\begin{lem}
  \[
X_L(f)=[L,\rho(\delta_Lf)]+\rho^*[L,\delta_Lf].
\]
\end{lem}
\begin{proof}
  Writing $L=\sum v_\alpha T^\alpha$ we have
  \begin{equation}\label{eq:GammafLaxoperator}
    \begin{split}
      \Gamma(f)\left(L(\tau_2)\right)&=\{f(\tau_1),v_\alpha(\tau_2)\}
      T_2^\alpha\\
      &=(d_\beta f)_1 \{v_\beta(\tau_1),v_\alpha(\tau_2)\}T_2^\alpha.
  \end{split}
\end{equation}
Now $\{v_\beta(\tau_1),v_\alpha(\tau_2)\}T_2^\alpha$ is the
$T_1^\beta$-component of $\{L(\tau_1)\poissontensor L(\tau_2)\}$, so is given by
 \eqref{eq:poissonLaxtensor}. Therefore we can write
 \eqref{eq:GammafLaxoperator}
 as a sum of two terms
(involving $r$ and $r^*$ respectively). The first term is (suppressing
the summation over the repeated indices here and below)
\begin{equation}
  \label{firsttermcom}
  \begin{split}
    &(d_\beta f)_1[L(\tau_2),\rho(T_2^{-\beta})\delta(\tau_1,\tau_2)]\\
    &=(d_\beta f)_1\Bigl(v_\alpha(\tau_2) T^{\alpha-\beta}_2
      \langle \delta(\tau_1,\tau_2)\rangle-T_2^{-\beta}\langle
      v_\alpha(\tau_2)\delta(\tau_1,\tau_2)\rangle\Bigr) 
      T_2^{\alpha}\rho(T_2^{-\beta}).
\end{split}
\end{equation}
Now recall that the projection $\Pi$ from section
\ref{sec:projection} satisfies for all difference operators $P,Q$ the
relation $\Pi\left(P_1Q_2\delta(\tau_1,\tau_2)\right)=Q(P_V)$ and that
we have $(d_\beta f)_V=\frac{\delta f}{\delta v_\beta}$. Therefore, by
applying $\Pi$ to the coefficients of \eqref{firsttermcom} we get as
first contribution to $X_f(L)$ (writing $T$ for $T_2$)
\begin{align*}
  \bigl(v_\alpha T^{\alpha-\beta}\langle \frac{\delta f}{\delta
      v_\beta}\rangle- T^{-\beta}\langle v_\alpha \frac{\delta
      f}{\delta v_\beta}\rangle \bigr ) T^\alpha\rho(T^{-\beta})
&=[L, \rho(T^{-\beta})\frac{\delta f}{\delta v_\beta}]\\
&=[L, \rho(\delta_L f)].
\end{align*}
Next consider
\begin{align*}
  [r^*,L(\tau_1)]
  &=[\rho(T_1^{\beta-\alpha})\delta(\tau_1,\tau_2),
      v_\alpha T_1^\alpha]T_2^{\alpha-\beta}\\
  &=\Bigl (\rho(T_1^{\beta-\alpha}
  )\langle v_a(\tau_1)\delta(\tau_1,\tau_2)\rangle -\\
&\qquad\qquad - v_a(\tau_1)T_1^\alpha\rho(T_1^{\beta-\alpha})\langle
  \delta(\tau_1,\tau_2)\rangle\Bigr) T_1^\beta T_2^{\alpha-\beta}.
\end{align*}
Taking the coefficient of $T_1^\beta$ of the last expression we see
that the second term of \eqref{eq:GammafLaxoperator} is
\begin{equation}
(d_\beta f)_1\left(\rho(T_1^{\beta-\alpha})\langle
  v_\alpha(\tau_1)\delta(\tau_1,\tau_2)\rangle
  -v_\alpha(\tau_1)T_1^\alpha\rho(T_1^{\beta-\alpha}) \langle
  \delta(\tau_1,\tau_2)\rangle \right)
T_2^{\alpha-\beta}.\label{eq:secondtermGammaf}
\end{equation}
Now the $r$-matrix satisfies 
\[
\rho(T_1^k)\delta(\tau_1,\tau_2)=\rho^*(T_2^{-k})\delta(\tau_1,\tau_2),
\]
so that we can rewrite \eqref{eq:secondtermGammaf} as
\begin{multline*}
  (d_\beta f)_1\Bigl(v_\alpha(\tau_2)\rho^*(T_2^{\alpha-\beta})\langle
    \delta(\tau_1,\tau_2)\rangle
    -\\T_2^{-\alpha}\rho^*(T_2^{\alpha-\beta}) \langle
    v_\alpha(\tau_2)\delta(\tau_1,\tau_2)\rangle \Bigr)
  T_2^{\alpha-\beta}.
\end{multline*}
Applying $\Pi$ to the coefficients we get as second contribution to
$X_f(L)$ (again writing $T$ for $T_2$)
\begin{align*}
( v_\alpha T^{\alpha-\beta}&
\langle \frac{\delta f}{\delta v_\beta}\rangle-
T^{-\beta}
\langle v_\alpha \frac{\delta f}{\delta v_\beta}\rangle
)\rho^*(T^{\alpha-\beta})\\
&=\rho^*[v_\alpha T^\alpha, T^{-\beta}\frac{\delta f}{\delta
  v_\beta}]=\rho^*[L,\delta_L f],
\end{align*}
so that indeed 
\[
X_f(L)=[L,\rho(\delta_Lf)]+\rho^*[L, \delta_Lf].
\]

\end{proof}
Recall that in \eqref{eq:defdelta_Lf} the difference operator
$\delta_Lf$ was determined up to addition of terms $\delta_{\ge2}$ and
$\delta_{\le-2}$. We have the same freedom in the vector field
$X_f(L)$: we have
\begin{align*}
  [L,\rho(\delta_{\ge2})]+\rho^*[L, \delta_{\ge2}]&=0\\
  [L,\rho(\delta_{\le-2})]+\rho^*[L, \delta_{\le-2}]&=0,
\end{align*}
so that $X_f(L)$ is independent of addition of these terms.

Next consider the Hamiltonians 
\[
H_n=\Sp(L^n/2n)\in V.
\]
We have for any $X\in \Derev(V,V)$
\[
X(H_n)\equiv\frac12\Sp(X(L)L^{n-1})\mod \mathfrak m V,
\]
so that we can represent the differential $\delta_L H_n$ of $H_n$ by the difference
operator $L^{n-1}/2$. Therefore the Hamiltonian vector field of $H_n$
is given by the Lax-equation
\begin{align*}
  X_{H_n}(L)&=[L, \rho(L^{n-1}/2)] + \rho^*[L, L^{n-1}/2]\\
            &=[L, \rho(L^{n-1}/2)] =[L/2, L_{0,+}^{n-1}-L_{-}^{n-1}] \\
            &=[L, L_{0,+}^{n-1})].
\end{align*}
In particular for $n=2$ we get the following equations:
\[
\dot A= 0,\quad\dot B = CT\inv A- AT C,\quad \dot C = C(T\inv -1) B.
\]
Putting $A=0$ we obtain an infinite analog of the finite Toda equation
(\ref{eq:BCequation}):
\begin{equation}
  \label{eq:InfTodaExplicit}
  \dot B=(1-T)C,\quad \dot C=C(T\inv-1)B.
\end{equation}
We want now clarify the motivation for the terminology $\Tr$ for the
maps $\K\to \mathbb C$, (\ref{eq:deftraceK}) and $LV\to \mathcal{L}V$,
see section \ref{sec:ClassFieldsAffin}. In the finite Toda lattice the
variables $q_i, p_i$ or $B_i, C_i$ give \emph{local} functions on the
phase space (attached to a site $i$ of the finite lattice, roughly
speaking), and the Hamiltonians, global functions on phase space, are
obtained as traces over the powers of the Lax matrix
(\ref{eq:LaxMatrix}). Similarly, for the infinite Toda lattice the
elements of $V$, difference polynomials in $B,C$, define
\emph{densities} of functions on the phase space $V(\K)$, and to
obtain honest functions one has to take the trace.

\chapter{$H_T$-Conformal Algebras}
\label{chap:H_TConformalAlg}

\section{Introduction}
\label{sec:Intro_Conformal}
 
In chapter \ref{chap:SingHamStruc} we introduced the notion of a
singular Hamiltonian structure on an $H_T$-Leibniz algebra $V$; this
was given by an evolutionary derivation $\Gamma\colon V\to \Derev(V,
V(\tau_1,\tau_2))$ satisfying the conditions of Definition
\ref{def:SingHamMap}. By projecting using the map $\Pi\colon
V(\tau_1,\tau_2)\to V$, see section \ref{sec:projection}, we obtain
for all $f\in V$ an evolutionary derivation 
$\Pi\circ\Gamma(f)=X_f\colon V\to V$.

In this chapter we will see that $X_f$ is just one of an
infinite collection of  products labeled by $h^*\in H_T^*$:
\[
\{h^*\}\colon  V\otimes V\to V,\quad f\otimes g\mapsto f_{\{h^*\}}g.
\]
In this notation $X_f=f_{\{1\}}$, where $1=1_{H_T^*}$, the identity of
the algebra $H_T^*=\Cz$, see section \ref{sec:dualseq}. The collection
of all \emph{conformal products} $f_{\{h^*\}}g$ for all
$h^*\in H^*_T$ will give $V$ the structure of what we will call an
$H_T$-Conformal Algebra.  These are analogs of the conformal algebras
from the usual theory of vertex algebras, see \cite{MR99f:17033}, cf.,
also \cite{math.QA/0209310} where they are called vertex Lie algebras.

We will define for any $H_T$-conformal algebra $C$ a \emph{singular
  vertex operator} $\Ysing(f,\tau)$ for all $f\in C$, see Section
\ref{sec:singvertop} below.  This is a generating series for the
conformal products $f_{\{h^*\}}$ and we will formulate the axioms of an
$H_T$-conformal algebra both in terms of the components $f_{\{h^*\}}$ and
also in terms of the singular vertex operator.

The singular vertex operators will be generalized in the next Chapter
\ref{chap:VertexAlg} to the vertex operators of an $H_T$-vertex
algebra.

\section[$H_T$-conformal and $H_T$-vertex Poisson algebras]{$H_T$-Conformal Algebras and $H_T$-Vertex Poisson
  Algebras.}\label{sec:H_TConfAlg}

If $\Gamma$ is a singular Hamiltonian map for an $H_T$-Leibniz
algebra $V$ we obtain for $f,g\in V$ a distribution $\Gamma(f)(g)\in
V(\tau_1,\tau_2)$ that we can uniquely write as
\begin{equation}
  \label{eq:DefGammainD}
  \Gamma(f)(g)=\mathsf{D}(f\otimes g)(\tau_2)\delta(\tau_1,\tau_2),
\end{equation}
where $\mathsf{D}(f\otimes g)$ is a difference operator in $ VH_T$. 

The properties of the Poisson brackets
$\{f(\tau_1),g(\tau_2)\}=\Gamma(f)g$ of a singular Hamiltonian
structure can be expressed in terms of the difference operators
$\mathsf{D}(f\otimes g)$. The fact that $\Gamma$ is evolutionary
implies that
\begin{equation}
  \label{eq:evolDfotimesg}
  \mathsf{D}(h.f\otimes g)=\mathsf{D}(f\otimes g)S(h), \quad h\in H_T,
\end{equation}
and the skew symmetry translates to
\begin{equation}
  \label{eq:skewsymdifference}
  \mathsf{D}(f\otimes g)=-\mathsf{D}(g\otimes f)^*,
\end{equation}
using the adjoint in the algebra $VH_T$, see
\eqref{eq:defadjoinedoperator}. The Lie property of $\Gamma$ gives a
more complicated condition on the $\mathsf{D}(f\otimes g)$ that we
will not need explicitly in the sequel.

To get explicit formulas it is convenient to choose a basis for
$H_T$. We can either use the translation
operator basis $\{T^n\}$ of $H_T$, with dual basis
$\{\delta_n\}_{n\in\mathbb Z}$ or the difference operator basis
$\{\Delta[n]\}$, with dual basis $\{\Delta^*[n]\}_{n\in\mathbb Z}$.
We  expand
\begin{equation}
\mathsf{D}(f\otimes g)=\sum_{n\in\mathbb Z} f_{n}g
T^n=\sum_{n\in\mathbb Z} f_{\langle n\rangle}g \Delta[n],
\label{eq:expansionsDfotimesg}
\end{equation}
with $f_ng, f_{\langle n\rangle}g\in V$. Here we use abbreviations:
\begin{equation}
f_n=f_{\{\delta_n\}}, \quad f_{\langle
  n\rangle}=f_{\{\Delta^*[n]\}},\label{eq:defconformalcomponentsTDelta}
\end{equation}
as the coefficients of the expansion \eqref{eq:expansionsDfotimesg}
should be labeled by elements of the dual of $H_T$.

Since $\Gamma(f)(-)$ is a derivation, we obtain in
this way derivations 
\[
f_n, f_{\langle n\rangle}\colon V\to V,
\]
for all $n\in \mathbb Z$. We think of the operations $n, \langle n\rangle\colon
V\otimes V\to V$ as giving (non-commu\-tative, non-associative)
multiplications on $V$. The properties of the multiplications are
simplest in terms of the translation operator basis, since the
multiplication and comultiplication of $H_T$ is simplest in that
basis.

\begin{lem} \label{lem:propertiesoperationSingHamStruct}
  For a singular Hamiltonian structure $\Gamma$ on a commutative
  $H_T$-Leibniz algebra $V$ the operations $n\colon V\otimes V\to V$
  satisfy
  \begin{Myenumerate}
  \item (Finiteness) For all $f,g\in V$ $f_{n}g=0$  
for all but a finite number of $n\in \mathbb Z$. \label{item:finiteness}
\item ($H_T$-covariance) For all $f\in V$, $n,p\in \mathbb Z$ we have
  \label{item:HTcovariance}
$(T^pf)_{n}=f_{n+p}$.
  \item (Skew-symmetry) For all $f,g\in V$, $n\in \mathbb Z$ we have \label{item:skewsymmetry}
$f_{n}g=-T^n\langle g_{-n}f\rangle$.

\item (Commutator) we have\label{item:commutator}
$[f_{m},g_{n}]=(f_{m-n}g)_{n}$ for all $f,g\in V$, $m, n\in \mathbb Z$.
  \end{Myenumerate}
\end{lem}

\begin{proof}
  The finiteness property of $f_{n}g$'s is clear since $\mathsf{D}(f\otimes
  g)\in VH_T$ 
  and any element of $VH_T$ has a finite expansion in
  $T^n$.

Taking $h=T^p$ in \eqref{eq:evolDfotimesg} we calculate
\[
\mathsf{D}(T^pf\otimes g)= \sum_{n\in\mathbb Z}(T^pf)_{n}g
T^n=\sum_{n\in\mathbb Z}f_{n}g T^{n-p}=\sum_{n\in\mathbb
  Z}f_{n+p}g T^n,
\]
so that
\[
(T^p f)_{n}g=f_{n+p}g.
\]
For skew-symmetry we use \eqref{eq:skewsymdifference}. Now
\[
  \mathsf{D}(g\otimes f)^*=\sum_{n\in\mathbb
    Z}T^{-n}g_{n}f=\sum_{n\in\mathbb Z}T^{-n}\langle g_{n}f\rangle
  T^{-n}=\sum_{n\in\mathbb Z}T^{n}\langle g_{-n}f\rangle T^n,
\]
so that
\[
f_{n}g=-T^n\langle g_{-n}f\rangle.
\]
For the commutator formula we calculate some Poisson brackets. For
$f,g,k\in V$ we have 
\begin{align*}
  \{f(\tau_1),\{g(\tau_2),k(\tau_3)\}\}&=\sum_{n,p\in\mathbb Z}
  f_{n}(g_{p}k)(\tau_3)T_3^n\delta(\tau_1,\tau_3)T_3^p\delta(\tau_2,\tau_3)\\
  &=\sum_{n,p\in\mathbb Z}
  f_{n}(g_{p}k)(\tau_3)T^{-n}_1T_2^{-p}\delta(\tau_1,\tau_3)\delta(\tau_2,\tau_3).\\
  \{g(\tau_2),\{f(\tau_1),k(\tau_3)\}\}&=\sum_{n,p\in\mathbb Z}
  g_{p}(f_{n}k)(\tau_3)T^{-n}_1T_2^{-p}\delta(\tau_1,\tau_3)\delta(\tau_2,\tau_3).\\
  \{\{f(\tau_1),g(\tau_2)\},k(\tau_3)\}
  &=\sum_{q,r\in\mathbb Z}(f_{q}g)_{r}T_1^{-q}T_3^r\delta(\tau_1,\tau_2)\delta(\tau_2,\tau_3)\\
  &=\sum_{q,r\in\mathbb Z}(f_{q}g)_{r}T_1^{-q-r}T_2^{-r}\delta(\tau_1,\tau_2)\delta(\tau_2,\tau_3).
\intertext{Making the substitution $q+r=n$, $r=p$ we obtain }
  &=\sum_{n,p\in\mathbb Z}(f_{n-p}g)_{r}T_1^{-n}T_2^{-p}\delta(\tau_1,\tau_2)\delta(\tau_2,\tau_3),
\end{align*}
so that by the Lie property of $\Gamma$ we find
\[
[f_{n},g_{p}]=(f_{n-p}g)_{p}.
\]
\end{proof}

We now observe that the operations $f_{n}$ (which are derivations of
the commutative algebra structure of $V$) with the properties
\eqref{item:finiteness}-\eqref{item:commutator} of Lemma
\ref{lem:propertiesoperationSingHamStruct} make sense also if $V$ is
just an $H_T$-module, not an $H_T$-Leibniz algebra. So we define:
\begin{defn}\label{def:HTconformal}
An \emph{$H_T$-Conformal algebra} is an $H_T$-module $C$ with
operations
\[
n\colon C\otimes C\to C,\quad n\in \mathbb Z,
\]
satisfying
\begin{Myenumerate}
  \item (Finiteness) For all $f,g\in C$ $f_{n}g=0$ for all but a
    finite number of $n\in \mathbb Z$.
\item ($H_T$-covariance) For all $f\in C$, $n,p\in \mathbb Z$ we have $(T^pf)_{n}=f_{n+p}$.
  \item (Skew-symmetry) For all $f,g\in C$, $n\in \mathbb Z$ we have $f_{n}g=-T^n\langle g_{-n}f\rangle$.
\item (Commutator) We have $[f_{m},g_{n}]=(f_{m-n}g)_{n}$ for all $f,g\in C$, $m, n\in \mathbb Z$.
  \end{Myenumerate}
\end{defn}

As a first result we have the Leibniz rule for conformal products.
\begin{lem}\label{lem:conformalLeibniz}
  Let $C$ be an $H_T$-conformal algebra, and $f, g\in C$. If $h\in
  H_T$ has coproduct $\sum h^\prime\otimes h^{\prime\prime}$ then for
  all $n\in\mathbb Z$ we have
  \begin{equation}
h(f_ng)=\sum (h^\prime
f)_n(h^{\prime\prime}g).\label{eq:conformalLeibniz}
\end{equation}
\end{lem}
\begin{proof}
  By $H_T$-covariance we have $f_ng=(T^nf)_0g$ so it suffices to check
  \eqref{eq:conformalLeibniz} in case $n=0$. Also, it suffices to
  check for $h=T^p$, $p\in\mathbb Z$, in which case the coproduct is
  just $T^p\otimes T^p$. Then
  \begin{align*}
    T^p(f_0g)&=T^p\left((T^pf)_{-p}g\right)=-(g_p(T^pf))\\
    &=-(T^pg)_0(T^pf)=(T^pf)_0(T^pg),
  \end{align*}
using skew-symmetry and $H_T$-covariance repeatedly.
\end{proof}
\begin{example}\label{ex:CToda}
  Let $\CToda$ be the free $H_T$-module generated by two elements $B$
  and $C$. We define products $x_{n}y$ for $x,y\in \CToda$,
  $n\in\mathbb Z$ by putting on generators
  \begin{equation}
    \label{eq:bnbcncproducts}
    B_{n}B=C_{n}C=0, \quad n\in\mathbb Z,
  \end{equation}
and 
\begin{align}
  \label{bncproducts}
  B_{0}C&=-C,&\quad B_{-1}C&=C\\
  C_{0}C&=C,&\quad C_{1}B&=-TC,\label{eq:cnbproducts}
\end{align}
and all other products $B_{n}C$ and $C_{n}B$ are zero. Since the
fundamental Poisson brackets \eqref{eq:fundamentalpoissonbracketsToda}
define a singular Hamiltonian structure on $V=\mathbb
C[B^{(n)},C^{(n)}]$ the products \eqref{eq:bnbcncproducts},
\eqref{bncproducts} and \eqref{eq:cnbproducts} define an
$H_T$-conformal structure on $\CToda$, called the \emph{first Toda
  conformal structure} on $\CToda$, cf.,  \cite{MR86m:58070}. We plan
to discuss higher Toda conformal structures in a future paper.

Note that $\CToda$ is a finite rank $H_T$-conformal subalgebra of $V$
and in a sense $V$ is the symmetric $H_T$-conformal algebra of
$\CToda$. Also note that $\CToda$ is a non trivial extension: if we
let $\mathrm{CC}$, $\mathrm{CB}$ be the rank 1 free Abelian
$H_T$-conformal algebras generated $C$, resp $B$ then we have an exact
sequence
\[
0\to \mathrm {CC}\to \CToda\to \mathrm{CB}\to 0.
\]\qed
\end{example}

\begin{example}\label{ex:AffineHTconformalalg}
  Let $\mathfrak g$ be a Lie algebra, and let $\Cg$ be the free
  $H_T$-module generated by $\mathfrak g$. Then define, for
  $x,y\in\mathfrak g$ and $p,q,n\in\mathbb Z$ products
\[
(T^px)_n(T^qy)= T^q[x,y]\delta_{n+p-q,0}.
\]
Then one checks that this gives $\Cg$ the structure of $H_T$-conformal
algebra, the commutator axiom follows from the Jacobi identity for
$\mathfrak g$. We call $\Cg$ the affine $H_T$-conformal algebra of the
Lie algebra $\mathfrak g$.\qed
\end{example}

The notion of a singular Hamiltonian structure on an $H_T$-Leibniz
algebra gives rise to a special class of $H_T$-conformal algebras.

\begin{defn}\label{defn:VertexPoisson}
  An \emph{$H_T$-Vertex Poisson Algebra} is a (commutative) $H_T$-Leibniz algebra
  $V$ with a compatible $H_T$-conformal algebra structure: the
  operations $f_{n}$ for $f\in V$, $n\in\mathbb Z$ are derivations
  of the algebra structure of $V$.
\end{defn}

So a singular Hamiltonian map $\Gamma\colon V\to
\Derev(V,V(\tau_1,\tau_2))$ gives an $H_T$-Leibniz algebra $V$ the
structure of $H_T$-vertex Poisson algebra. Conversely, let $V$ be an
$H_T$-vertex Poisson algebra. Then we can define an evolutionary
derivation $\Gamma\colon V\to \Derev(V,V(\tau_1,\tau_2))$ by
\[
\Gamma(f)(g)=\sum_{n\in\mathbb
  Z}f_{n}g(\tau_2)T_2^{n}\delta(\tau_1,\tau_2).
\]
One then checks that $\Gamma$ is in fact a singular Hamiltonian map.

\section{The Lie Algebra of a Conformal Algebra}
\label{sec:LieAlgConfAlg}

To obtain from a singular Hamiltonian structure $\Gamma$ on $V$ a
conformal algebra structure we used the expansion of the difference
operators $\mathsf{D}(f\otimes g)$ in the basis $\{T^n\}$ of $H_T$. 
The operations $\langle n\rangle\colon V\otimes V\to V$ related to the difference
basis, see \eqref{eq:expansionsDfotimesg}, satisfy properties similar to
those of ${n}$, but they are more complicated. 
We will not need explicitly the operations $\langle n\rangle$ in the sequel,
except for the zeroeth one, which is given by
\begin{equation}
f_{\langle 0\rangle}g=\sum_{n\in\mathbb Z}f_{n}g.\label{eq:DefZeroComponent}
\end{equation}
In the notation of \eqref{eq:defconformalcomponentsTDelta} we have
$f_{\langle0\rangle}=f_{\{1_{H^*}\}}$.

Recall the augmentation ideal $\mathfrak m\subset H_T$ and the
submodule $\mathfrak m M$ of total differences of an $H_T$-module $M$, see
Section \ref{sec:ClassFieldsAffin},
\begin{lem}
  In any $H_T$-conformal algebra $V$ the operation $\langle 0\rangle$ given by
  \eqref{eq:DefZeroComponent} satisfies
  \begin{Myenumerate}
  \item 
$(Tf)_{\langle 0\rangle}=f_{\langle 0\rangle}$.
\item $f_{\langle 0\rangle}(Tg)=T\langle{f_{\langle 0\rangle}g}\rangle{}$.
\item $f_{\langle 0\rangle}g-g_{\langle 0\rangle}f\in{\mathfrak m}C$.
\item $[f_{\langle 0\rangle},g_{\langle 0\rangle}]=(f_{\langle 0\rangle}g)_{\langle 0\rangle}$.
  \end{Myenumerate}
\end{lem}

\begin{proof}
  We have by $H_T$ covariance
\[
(Tf)_{\langle 0\rangle}g=\sum_{n\in\mathbb Z}(Tf)_{n}g=\sum_{n\in\mathbb
  Z}(f)_{n+1}g=f_{\langle 0\rangle}g.
\]
By skew-symmetry and $H_T$-covariance we get
\begin{align*}
  f_{\langle 0\rangle}(Tg)&=\sum_{n\in\mathbb Z}f_{n}(Tg) 
  =-\sum_{n\in\mathbb Z}T^n\langle (Tg)_{-n}f\rangle=-\sum_{n\in\mathbb Z}TT^{n-1}\langle g_{1-n}f\rangle
  \\&=T\sum_{n\in\mathbb Z}\langle f_{n-1}g\rangle=T\langle f_{\langle 0\rangle}g\rangle.
\end{align*}
Using the commutator we find
\[
[f_{\langle 0\rangle},g_{\langle 0\rangle}]q=\sum_{m,n\in \mathbb
  Z}[f_{m},g_{n}]=\sum_{m,n\in \mathbb Z}(f_{m-n}g)_{n}=\sum_{p,q\in
  \mathbb Z}(f_{p}g)_{q} =(f_{\langle 0\rangle}g)_{\langle 0\rangle}.
\]
\end{proof}

Recall the trace map $\Tr:M\to M/\mathfrak m M$, see section
\ref{sec:ClassFieldsAffin}, and write $\tilde f$ for the image
$\Tr(f)$ in $C_{\mathfrak m}$ of $f\in C$, see \eqref{eq:DefTraceM}.
Then the lemma shows that in case $C$ is an $H_T$-conformal algebra
the bilinear map
\begin{equation}
  [\,,]\colon C_{\mathfrak m}\otimes C_{\mathfrak m}\to C_{\mathfrak m}, \quad
  [\tilde f, \tilde g]= \widetilde{f_{\langle
  0\rangle}g}\label{eq:DefComConfLiealg}
\end{equation}
is well defined, and defines a Lie algebra structure on $C_{\mathfrak  m}$.

\begin{example}
  In case of the Toda conformal algebra $\CToda$ of Example
  \ref{ex:CToda} the module of coinvariants is two dimensional:
\[
\CToda_{\mathfrak{m}}=\CToda/\mathfrak{m}\CToda=\mathbb C \tilde B\oplus
\mathbb C \tilde C,
\]
with \emph{vanishing} Lie bracket: $[\tilde B,\tilde C]=0$.\qed
  \end{example}

  \begin{example}
    For the affine $H_T$-conformal algebra $\Cg$ of a Lie algebra
    $\mathfrak g$, see \ref{ex:AffineHTconformalalg}, the module of
    coinvariants $\Cg_{\mathfrak m}\simeq \mathfrak g$ and the Lie
    bracket \eqref{eq:DefComConfLiealg} is just the standard bracket
    of $\mathfrak g$.\qed
  \end{example}
\section{Singular Vertex operators}
\label{sec:singvertop}
In this section we reformulate the axioms of an $H_T$-conformal
algebra $C$, see Definition \ref{def:HTconformal}, in terms of the
generating series of the conformal products of $f\in C$.  Define the
\emph{singular vertex operator} of $f$ by
\[
\Ysing(f,\tau)=\sum_{n\in \mathbb
  Z}f_{\langle n\rangle}S(\Delta[n])\frac{1}\tau=\sum_{n\in \mathbb Z}f_{n}S(T^n) \frac{1}\tau.
\]
\begin{lem}\label{lem:PropertiesSingVertex}
  Let $C$ be an $H_T$-conformal algebra. The singular vertex operator
  for $f\in C$ satisfies
  \begin{Myenumerate}
  \item (Finiteness) $\Ysing(f,\tau)g$ belongs to $CH_T\frac1\tau$ for
    all $g\in C$.
\item ($H_T$-covariance) For all $h\in H_T$ we have $\Ysing(hf,\tau)=h_K\Ysing(f,\tau)$.
\item (Skew-symmetry) For all $g\in C$
\[
\Ysing(f,\tau)g=-\Sing_\tau[\mathcal R_C(\tau)\YsingS(g,\tau)f].
\]
\item (Commutator) For all $g\in C$
\[
[\Ysing(f,\tau_{1}),\Ysing(g,\tau_{2})]
=\Sing_{\tau_{2}}\left(\Ysing(\mathcal{R}^S_W(\tau_{2}).\Ysing(f,\tau_{1})g,\tau_{2})\right),
\]
where $W=(H_T)\frac1{\tau_{1}}$.
  \end{Myenumerate}
\end{lem}

\begin{proof}
Finiteness is clear. For $H_T$-covariance note that
\[
\Ysing(f,\tau)g=\mathsf{D}(f\otimes g)^S\frac1\tau.
\]
So by \eqref{eq:evolDfotimesg} we have
\[
\Ysing(hf,\tau)g=\mathsf{D}(hf\otimes
g)^S\frac1\tau.=\mathsf{D}(f\otimes g)^Sh\frac1\tau=h_{\K}\Ysing(f,\tau)g.
\]
For skew-symmetry we use \eqref{eq:skewsymdifference} to calculate:
\begin{align*}
  \Ysing(f,\tau)g&=\mathsf{D}(f\otimes
  g)^S\frac1\tau=-(\mathsf{D}(g\otimes f)^*)^S\frac 1\tau\\
  &=-\Sing(\mathcal{R}_V(\tau)\mathsf{D}(g\otimes f)\frac1\tau)\\
  &=-\Sing(\mathcal{R}_V(\tau)\YsingS(g,\tau)f),
\end{align*}
using Lemma \ref{lem:PstarrhoSing}, where
\[
\YsingS(g,\tau)f=\mathsf{D}(g\otimes f)\frac1\tau.
\]
Finally
\begin{align*}
  [\Ysing(f,\tau_{1}),\Ysing(g,\tau_{2})]
  &=\sum_{m,n}[f_{m},g_{n}]T^{-m}\langle\frac1{\tau_{1}} \rangle T^{-n}\frac1{\tau_{2}} \\
  &=\sum_{n,m}(f_{m-n}g)_{n}T^{-m}\langle\frac1{\tau_{1}}\rangle T^{-n}\frac1{\tau_{2}} \\
  &=\sum_{a,b}(f_{a}g)_{b}T^{-b}( T^{-a}\langle\frac1{\tau_{1}}\rangle\frac1{\tau_{2}})\\
  &=\sum_{a,b}(f_{a}g)_{b}P_{a,b}^*\frac1{\tau_{2}}, 
\intertext{where $P_{a,b}=T^{-a}\langle \frac1{\tau_1}\rangle T^b$.
    Then by Corollary \ref{cor:TbTaSing} we have}
  &=\Sing_{\tau_2}\left(\mathcal{R}^S_W (\tau_2)(f_ag)_b P_{a,b}^S \frac 1{\tau_2}\right)\\
  &=\Sing_{\tau_2}\left(\mathcal{R}^S_W (\tau_2)(f_ag)_b
    (T_1^{-a}\frac1{\tau_1})T_2^{-b}
    \frac 1{\tau_2}\right)\\
  &=\Sing_{\tau_2} (\Ysing\left(\mathcal{R}^S_W
    (\tau_2)\Ysing(f,\tau_1)g,\tau_2\right)).
\end{align*}
\end{proof}

Of course, we can reverse our procedure, and \emph{define} an
$H_T$-conformal algebra as an $H_T$-module $C$ together with for all
$f\in C$ a singular vertex operator $\Ysing(f,\tau)$, satisfying the
properties of Lemma \ref{lem:PropertiesSingVertex}. Expanding the
singular vertex operator we find that the components $f_{n}$ satisfy
the properties of Definition \ref{def:HTconformal}.

\section{Holomorphic $H_T$-vertex algebras}
\label{sec:holVertex}
Let $V$ be a commutative $H_T$-Leibniz algebra. For $f\in V$ define
the \emph{holomorphic vertex operator} by
\begin{equation}
\Yhol(f)=\mathcal R(f)=\sum_{n\in\mathbb Z} \Delta[n](f)\otimes \Delta^*[n],
\label{eq:defYholHTLeibniz}
\end{equation}
where $\mathcal R$ is the untwisted exponential operator for $H_T$
introduced in section \ref{sec:exp_op}, see
\eqref{eq:HTexponentials}, so that $\Yhol$ is just a suggestive
notation for the coaction map $\md_V: V\to V\hat \otimes H_T^*$
introduced in section \ref{sec:H_LeibnizAlgSpec}. In particular
$\Yhol$ is multiplicative, see section \ref{sec:MultTwistedExp}. We define 
\[
\YholS(f)=\sum \Delta[n](f)\otimes
S(\Delta^*[n])=\mathcal{R}^{S}(f_1).
\]
By Lemma \ref{lem:InverseExp} we have
\[
\YholS(\Yhol(f))=\Yhol(\YholS(f))=f\otimes1_{H_T^*}.
\]
It then follows that we have \emph{Skew-symmetry}:
\[
\Yhol(f)g=\mathcal R\left(\YholS(g)f\right).
\]
The $H_T$-covariance of exponential operators, see section
\ref{sec:HinvarianceTwistedExp}, implies that
\[
\Yhol(hf)=h_{H_T^*}(\Yhol(f)),\quad \YholS(hf)=S(h)_{H_T^*}(\YholS(f)).
\]

\begin{defn}\label{def:holovertexalgebra}
  A holomorphic $H_T$-vertex algebra is an $H_T$-module, together with
  a distinguished element $1_V\in V$ and a map
\[
\Yhol\colon  V\otimes V\to V\hat\otimes H_T^*, \quad f\otimes g\mapsto \Yhol(f)g,
\]
such that
\begin{Myenumerate}
\item (Vacuum) For all $f,g\in V$ 
\[\Yhol(1_V)g=g\otimes1_{H_T^*}, \quad \Yhol(f)1_V|_0=f.
\]
\item ($H_T$-covariance)  We have $\Yhol(hf)=h_{H_T^*}\Yhol(f)$ for all $f\in V$, $h\in H_T$.
\item (Skew-symmetry) We have $\Yhol(f)g=\mathcal
  R\left(\YholS(g)f\right)$ for all $f,g\in V$.
\item (Commutator) We have $[\Yhol(f),\Yhol(g)]=0$ for all $f,g\in V$.
\end{Myenumerate}
\end{defn}
It is clear that any commutative $H_T$-Leibniz algebra $V$ is a
holomorphic $H_T$-vertex algebra, via the choice
\eqref{eq:defYholHTLeibniz} of holomorphic vertex operator.
Conversely, if $V$ is an holomorphic $H_T$-vertex algebra, then we
define a product on $V$ by evaluating at zero:
\[
f.g=\Yhol(f)g|_0.
\]
Then one checks that this gives $V$ the structure of a commutative
$H_T$-Leibniz algebra. 

\section{Extension of $H_T$-conformal Structure 
to the Affinization}
\label{sec:extAff}
Recall the notion of the affinization $LM=M\otimes K$ of an
$H_T$-module $M$ and the canonical trace map $\Tr\colon LM\to
\mathcal{L}M=LM/\mathfrak{m}LM$, see Section
\ref{sec:ClassFieldsAffin}. 
In this section we take $M=C$, where $C$ is an $H_T$-conformal
algebra, and show that the affinization $LC$ inherits an
$H_T$-conformal structure from the $H_T$-conformal structure of $C$
and the holomorphic vertex algebra structure of $\K$. In particular
this gives $\mathcal{L}C$ a Lie algebra structure, by Section
\ref{sec:LieAlgConfAlg}. In case $C=V$ is an $H_T$-vertex Poisson
algebra, see Section \ref{sec:H_TConfAlg}, the Lie bracket on
$\mathcal LV$ is the Poisson bracket of the components of the
classical fields of Chapter \ref{chap:SingHamStruc}, see Example \ref{ex:VertexPoissonLiealg}.

For $f\in C$ we have the singular vertex operator 
\[
\Ysing(f,\tau)g=\sum_{n\in\mathbb Z}f_{\langle n\rangle}gS(\Delta[n])\frac1\tau
\]
For $p\in K$ we have the holomorphic vertex operator $\Yhol(p)$, see
\eqref{eq:defYholHTLeibniz}, and by combining $\Ysing$ on $C$ with
$\Yhol$ on $K$ we define a singular vertex operator on $LC$:
\[
\YsingLC(f\otimes
p,\tau):=\Sing\left(\Ysing(f,\tau)\otimesd\Yhol(p,\tau)\right).
\]
Here $\otimesd$ is a combination of tensor product and multiplication:
we have explicitly
\begin{equation}
  \label{eq:otimesdot}
\Ysing(f,\tau)\otimesd\Yhol(p,\tau)=
\sum_{m,n\in\mathbb Z}(f_{[m]}\otimes \Delta[n]\langle
p\rangle)\Delta^*[n]S(\Delta[m])\frac1\tau).  
\end{equation}
Now by \eqref{eq:oppsitesignstau} we have a vanishing contribution to
\eqref{eq:otimesdot}, and hence to $\YsingLC$, from the terms in which
$m,n$ have opposite sign. In case $0\le m<n$
\[
\Delta^*[n]S(\Delta[m])\frac1\tau,\quad \overline{\Delta^*[n]}S(\overline{\Delta[m]})\frac1\tau
\]
is nonsingular, so these terms will also not contribute to
$\YsingLC$. Finally, by \eqref{eq:samesigntau} and the factorization
\eqref{eq:factortau} we have, if $m=n+s$, $n,s\ge0$
\begin{align*}
  \Delta^*[n]S(\Delta[n+s])\frac1\tau&=\tau[n]\frac1{\tau[n+s+1]}=T^{-n}S(\Delta[s])\frac1\tau\\
  \Delta^*[-n]S(\Delta[-n-s])\frac1\tau&=\overline{\tau[n]}\frac1{\overline{\tau[n+s+1]}}=T^{n}S(\Delta[-s])
  \frac1\tau
\end{align*}
Hence
\begin{multline}\label{eq:ExplicitVertexLC}
    \YsingLC(f\otimes
    p,\tau)=\sum_{m>0}\sum_{\ell=0}^m(f_{\langle m\rangle}\otimes\Delta[\ell]p)\otimes
    T^{-\ell}S(\Delta[m-\ell])\frac1\tau+\\
    \sum_{m>0}\sum_{\ell=0}^m(f_{\langle-m\rangle}\otimes\overline{\Delta[\ell]}p)\otimes
    T^{\ell}S(\overline{\Delta[m-\ell]})\frac1\tau.
\end{multline}
\begin{lem}
  $\YsingLC$ gives $LC$ the structure of $H_T$-conformal algebra.
\end{lem}
\begin{proof}
We need to check the properties of Lemma
\ref{lem:PropertiesSingVertex}. Finiteness is manifest from
\eqref{eq:ExplicitVertexLC}, as in any $H_T$-conformal algebra
$f_{\langle m\rangle}g=0$ for all but a finite number of $m\in\mathbb Z$.

For $h\in H_T$ we have, if $\pi(h)=\sum h^\prime\otimes h^{\prime\prime}$ is the
coproduct, using the $H_T$-covariance of $\Ysing$ and $\Yhol$,
\begin{align*}
\YsingLC(h.(f\otimes p),\tau)
&=\sum\Sing_{\tau}\left( \Ysing(h^\prime f,\tau)\otimesd
  \Yhol(h^{\prime\prime} p,\tau)\right) \\
&=\sum\Sing_{\tau}\left(h^\prime_K \Ysing(f,\tau)\otimesd
h^{\prime\prime}_K \Yhol( p, \tau)\right)\\
&=\Sing_{\tau}\left( h_K\Bigl(\Ysing(f,\tau)\otimesd \Yhol( p,\tau)\Bigr)\right)\\
&=h_{\K}\Sing_{\tau}\left( \Ysing(f,\tau)\otimesd \Yhol( p,\tau)\right),\\
&=h_{\K}\YsingLC(f\otimes p,\tau),
\end{align*}
since for any distribution we have
\[
\Sing(h\mathcal{D})=h\Sing(\mathcal{D}).
\]
This proves $H_T$-covariance for $\YsingLC$.

For skew-symmetry we need the following simple fact about the
projection $\Sing$, cf., \cite{math.QA/0209310}: 
let $X,Y$ be vector spaces and let  $A\in X\otimes
\K$ and $B\in Y\otimes H_T$, then
\begin{equation}
\Sing( A\otimesd B)=\Sing \left(\Sing (A)\otimesd B\right).
\label{eq:SingSingidentity}
\end{equation}
Then we have
\begin{align*}
  \YsingLC(f&\otimes p,\tau)(g\otimes q)
= \Sing\left(\Ysing(f,\tau)g\otimesd
  \Yhol(p,\tau)q\right),
\intertext{and by skewsymmetry in the conformal algebra $C$ and the
  holomorphic vertex algebra $\K$:}
&= \Sing\left(\Sing(\mathcal R_C(\tau)\YsingS(g,\tau)f)\otimesd \mathcal
  R_K(\tau)\Yhol(q,\tau)p\right),\\
&= \Sing\left(\mathcal R_C(\tau)\YsingS(g,\tau)f \otimesd
  \mathcal R_K\YholS(q,\tau_{12})p\right),\quad\text{by \eqref{eq:SingSingidentity}}\\
\intertext{and by multiplicativity of exponentials, Section \ref{sec:MultTwistedExp},:}
&= \Sing\left(\mathcal R_{LC}(\tau)(\YsingS(g,\tau) f \otimesd
  \YholS(q,\tau)p)\right),\\
&= \Sing\left(\mathcal R_{LC}(\tau)\Sing\left (\YsingS(g,\tau)f
    \otimesd \YholS(q,\tau)p\right)\right),\quad\text{again by \eqref{eq:SingSingidentity}}\\
&= \Sing\left(\mathcal R_{LC}(\Ysing^{LC,S}(g\otimes q,\tau)f\otimes p)\right),
\end{align*}
proving skew-symmetry for $\Ysing^{LC}$.

To check commutator for $LC$ we first note that it follows from Lemma
\ref{lem:productExpandcoproduct} that 
\[
\mathcal{R}_{\K}(\tau_1)\langle
p\rangle\mathcal{R}_{\K}(\tau_2)\langle q\rangle=
\mathcal{R}_{\K}(\tau_2)\left(\mathcal{R}^S_{\K_{\tau_1}}(\tau_2)\mathcal{R}_{\K}(\tau_1)\langle
  p\rangle q\right).
\]
Here $p,q\in\K$ are rational functions of $\tau$, say. Then
$\mathcal{R}_{\K}(\tau_1)p$ is rational in $\tau$ and holomorphic in
$\tau_1$, and $\mathcal{R}^S_{\K_{\tau_1}}$ is the exponential
operator acting on rational functions of $\tau_1$. We will use similar
notation below without further elaboration.

Then we calculate
\begin{align*}
  [&\YsingLC(f\otimes p,\tau_1),\YsingLC(g\otimes
  q,\tau_2)]\\&=\Sing_{\tau_1,\tau_2}\left(
    [\YsingC(f,\tau_1),\YsingC(g,\tau_2)]\otimes\mathcal{R}_{\K}(\tau_1)\langle
    p\rangle\mathcal{R}_{\K}(\tau_2)\langle q\rangle\right)\\
  &=\Sing_{\tau_1,\tau_2}\left(
    \YsingC(\mathcal{R}^S_{\K_{\tau_1}}(\tau_2)\YsingC(f,\tau_1)g,\tau_2)\otimes\mathcal{R}_{\K}(\tau_1)\langle
    p\rangle\mathcal{R}_{\K}(\tau_2)\langle q\rangle\right)\\
  &=\Sing_{\tau_1,\tau_2}\left(
    \YsingC(\mathcal{R}^S_{\K_{\tau_1}}(\tau_2)\YsingC(f,\tau_1)g,\tau_2)\otimes
    \mathcal{R}_{\K}(\tau_2)\left(\mathcal{R}^S_{\K_{\tau_1}}(\tau_2)\mathcal{R}_{\K}(\tau_1)\langle
      p\rangle q\right)\right)\\
  &=\Sing_{\tau_2}\left(
    \YsingLC(\mathcal{R}^S_{\K_{\tau_1}}(\tau_2)\YsingLC(f\otimes
  p,\tau_1)g\otimes q,\tau_2)\right).
\end{align*}
\end{proof}

The Lie bracket of $\mathcal{L}C$ is induced by the $\langle 0\rangle$
bracket on $LC$. This is the coefficient of  $\sum_{n\in\mathbb
  Z}T^{n}\frac1\tau$ in \eqref{eq:ExplicitVertexLC}, hence
\[
(f\otimes p)_{\langle 0\rangle}g\otimes q=\sum_{n\in\mathbb
  Z}f_{\langle n\rangle}g
 {\Delta[n]\langle p\rangle}q.
\]
Write $f_{\{p\}}=\Tr(f\otimes p)\in\mathcal{L}C$. The Lie bracket on
$\mathcal LC$ is then given by
\begin{equation}
  \label{eq:LieLC}
  [ f_{\{ p\}}, g_{\{q\}}]=\sum_{n\in\mathbb Z}f_{\langle
  n\rangle}g_{\{ \Delta[n]\langle  p\rangle q\}}.
\end{equation}
Define the generating series of elements of $\mathcal{L}C$
(''currents'') for $f\in C$ by
\begin{equation}
f(\tau)=\Tr_{\tau_0}\left(f\otimes
  \delta(\tau_0,\tau)\right)=\sum_{n\in\mathbb
  Z}f_{\{\tau(n)\}}\tau(-n-1).
\label{eq:DefCurrentsLC}
\end{equation}
The currents are $\hat H_T$-covariant: for all $h\in \hat H_T$
\[
h_{\mathcal{L}C}.f(\tau)=\Tr_{\tau_0}\left(f\otimes
  S(h)_0\delta(\tau_0,\tau)\right)=\Tr_{\tau_0}\left(f \otimes
  h_1\delta(\tau_0,\tau)\right)=h_{\K}f(\tau).
\]
The current commutator is given by
\begin{equation}
\begin{split}
  [f(\tau_1),g(\tau_2)]&=[\Tr_{\tau_0}\left(f\otimes
    \delta(\tau_0,\tau_1)\right), \Tr_{\tau_0}\left(g\otimes
    \delta(\tau_0,\tau_2)\right)]\\
  &=\sum_{n\in\mathbb Z}\Tr_{\tau_0}\left(f_{\langle n\rangle}g\otimes
    \Delta[n]_0\delta(\tau_0,\tau_1)\delta(\tau_0,\tau_2) \right)\\
&=\sum_{n\in\mathbb Z}f_{\langle n\rangle}g(\tau_2)\Delta[n]_2\delta(\tau_1,\tau_2).
\end{split}\label{eq:conformalcurrentcommutator}
\end{equation}

\begin{example}
  \label{ex:affineLiealgebra}
Let $C\mathfrak g$ be the $H_T$-conformal algebra of the Lie algebra
$\mathfrak g$, see Example \ref{ex:AffineHTconformalalg}. Then let
$\mathcal{L}\mathfrak g=\mathcal{L}C\mathfrak g$ be the corresponding
Lie algebra. If $\{X^\alpha\}$ is a basis for $\mathfrak g$, then
$\mathcal{L}\mathfrak g$ is spanned by elements
\[
X^\alpha_{\{p\}}=\Tr(X^\alpha\otimes p),\quad p\in K,
\]
If $X, Y\in\mathfrak g$ the commutator in $\mathcal{L}\mathfrak g$ is 
\[
[X_{\{p\}},Y_{\{q\}}]=[X,Y]_{\{pq\}}.
\]
The corresponding currents have commutator
\[
[X(\tau_1),Y(\tau_2)] = [X,Y](\tau_2)\delta(\tau_1,\tau_2).
\]
\qed
\end{example}
\begin{example}\label{ex:TodaLiealg}
  Let $\CToda$ be the Toda conformal algebra, see Example
  \ref{ex:CToda}. Let $\LToda=\mathcal{L}\CToda$. Then elements of
  $\LToda$ are linear combinations of expressions of the form
  $X_{\{p\}}$, where $X=B,C$ and $p\in\K$, with commutator
  \begin{equation}
[B_{\{p\}},C_{\{q\}}]=-C_{\{\overline{\Delta}\langle p\rangle q\}},
    \label{eq:BCcommutator}
      \end{equation}
      and current commutator
\[
[B(\tau_1),C(\tau_2)]=-C(\tau_2)\overline{\Delta}_2\delta(\tau_1,\tau_2).
\]
Of course this looks like \eqref{eq:fundamentalpoissonbracketsToda}
and in the next example we explain the connection.\qed
\end{example}
\begin{example}
  \label{ex:VertexPoissonLiealg}
  Let $V$ be an $H_T$-vertex Poisson algebra. This is in particular an
  $H_T$-conformal algebra, and so we obtain a Lie algebra
  $\mathcal{L}V$. Elements of $\mathcal{L}V$ are of the form
  $f_{\{p\}}$, for $f\in V$ and $p\in \K$. The commutator of two such
  elements is given by
\[
[f_{\{p\}},g_{\{q\}}]=\sum_{n\in\mathbb Z} f_{\langle
  n\rangle}g_{\{\Delta[n]\langle p\rangle q}.
\]
The corresponding current is given by
\eqref{eq:conformalcurrentcommutator}. 

Recall form Section \ref{sec:ClassFieldsAffin} that the coefficients
of the classical field $C(f)$, for $f\in V$, can be identified with
elements of $\mathcal{L}V$:
\[
C(f)(p)\simeq  f_{\{p\}}=\Tr(f\otimes p),\quad p\in \K.
\]
Under this identification the Poisson bracket
(\ref{eq:PoissonBracketfg}) is nothing but the commutator
(\ref{eq:conformalcurrentcommutator}) in the conformal algebra $V$.

\end{example}

 \chapter{$H_T$-Vertex Algebras.}
\label{chap:VertexAlg}

\section{Introduction}
\label{sec:HTvertexIntro}

Recall that in chapter \ref{chap:SingHamStruc} we defined the
classical field $C(f)$, for $f\in V$, where $V$ is a commutative
$H_T$-Leibniz algebra. This is a function $V(L)\to L$, for some other
$H_T$-Leibniz algebra $L$, where
$V(L)=\Hom_{H_T\mathrm{-alg}}(V,L)$. Recall the localization $\K$ of
section \ref{sec:localizationHTcirc}.  In case $L=\K$, which we will
assume from now on, the classical field has an expansion
)\begin{equation}
C(f)=f(\tau)=\sum_{n\in\mathbb Z} f_{(n)}\tau(-n-1).\label{eq:expansionKcircclassfield}
\end{equation}
The coefficients $f_{(n)}$ are $\mathbb C$-valued functions (on
$V(\K)$). In particular these coefficients generate a commutative
$\mathbb C$-algebra. We can also think of a classical field as a
distribution with values in functions on $V(\K)$: for each $F\in \K$
we get the function $f_{\{F\}}=C(f)(F)$ given by
$f_{\{F\}}(\sigma)=\Tr(\sigma(f)F)$. In particular
$f_{(n)}=f_{\{\tau(n)\}}$, $n\in\mathbb Z$.

In this chapter we will quantize this situation, and define fields
(which will be called vertex operators) that have expansions like
\eqref{eq:expansionKcircclassfield}, but now the coefficients will be
linear maps on a vector space $V$, and they will generate a
non-commutative algebra.

This leads to the definition of an $H_T$-vertex algebra, which will be
a quantization of the singular Hamiltonian structures of chapter
\ref{chap:SingHamStruc}, which are the $H_T$-vertex Poisson algebras of chapter
\ref{chap:H_TConformalAlg} and a generalization of the $H_T$-conformal
algebras of chapter \ref{chap:H_TConformalAlg}.

\section{Fields}
\label{sec:fields}

Let $W$ be a vector space. Recall the vector space $W[[\tau^{\pm1}]]$
of $W$-valued distributions $\mathcal{D}$ on $\K$ in the variable
$\tau$ with decomposition
$\mathcal{D}=\mathcal{D}_\Hol+\mathcal{D}_\Sing$ in holomorphic and
singular part, see section \ref{sec:ExpDistributions}.

Let $V$ be a vector space. We
will call $f(\tau)\in\End(V)[[\tau^{\pm1}]]$ a \emph{field
  on }$V$ if the action of $f_{\mathrm {Sing}}(\tau)$ on $V$ is
rational: for all $v\in V$ the $V$-valued distribution obtained by
applying $f_\Sing$ to $v$ is has rational kernel: we have
\[
f_{\Sing}(\tau)v(G)=\Tr(F_v G), \quad G\in K
\]
for some $F_v\in V\otimes \Ksing$. So we can expand in a finite sum:
\begin{equation}
f_\Sing(\tau)v=\sum_{n,k}v_{n,k}S(e_{n,k})\frac1\tau,\quad
v_{n,k}\in V.
\label{eq:rationalexpansionftauv}
\end{equation}
Note that the rationality condition does \emph{not} imply that for all
$v\in V$ we have in the expansion \eqref{eq:expansiondistribution} of
the $V$-valued distribution $f_\Sing(\tau)v$ a finite sum for the
singular part (as is the case, mutatis mutandis, for fields in the
usual vertex algebras based on $H_D$).

In case $V$ is an $\hat H_T$-module and $f(\tau)$ is a field on $V$
the $\End(V)$-valued distributions $\ad^V_h f(\tau)$ and $h.f(\tau)$
are in fact again fields on $V$, for all $h\in \hat H_T$.

We expand the field $f(\tau)$ as in
\eqref{eq:expansiondistribution}:
\begin{equation}
f(\tau)=\sum_{n\in\mathbb Z} f_{(-n-1)}\tau(n),
\quad f_{(-n-1)}\in\End(V).
 \label{eq:ExpansionFieldtau}
\end{equation}
Consider $f_\Hol(\tau)$. This is an $\End(V)$-valued distribution on
$\Ksing$. By \eqref{eq:isomorphalpha} we can consider it a
distribution on $\hat H_T$. So for each $h\in\hat H_T$ we get a linear
map
\begin{equation}
f_{[h]}:V\to V, \quad f_{[h]}=\langle h,f_{\Hol}(\tau)\rangle=h_{\K}
f_{\Hol}(\tau)|_{\tau=0}.\label{eq:DefhsubV}
\end{equation}
We can also interpret $f_\Hol(\tau)$ as a distribution 
on $\mathbb
C[\Delta]$ or $H_T$ that extends to $\hat H_T$, so we have expansions
\begin{equation}
f_{\Hol}(\tau)=\sum_{k\ge0}f_{[\Delta[k]]}\tau[k]=\sum_{n\in\mathbb
  Z}f_{[\Delta[n]]}\Delta^*[n]=\sum_{n,\ell}f_{[e_{n,\ell}]}e_{n,\ell}^*,
\label{eq:expansionholovertexoperatorY(f,tau)}
\end{equation}
see section \ref{sec:extExpDistr}.

If $\mathcal{D}$ is a holomorphic field, i.e., a distribution that
vanishes on $\Czpol$, then we define 
\begin{equation}
  \label{eq:deftau=0}
\mathcal{D}|_{\tau=0}=\mathcal{D}(\frac1\tau),
  \end{equation}
  i.e., we define the \emph{constant term } of the holomorphic
  $\mathcal{D}$ as the value of $\mathcal{D}$ at $\frac1\tau$.

Since the field $f(\tau)$ is an $\End(V)$-valued distribution on $\K$ we get for
each $F\in\K$ a linear map
\begin{equation}
  f_{\{F\}}:V\to V,\quad f_{\{F\}}=\Tr_\tau\Bigl(f(\tau)F(\tau)\Bigr).\label{eq:defFcomponent}
  \end{equation}
So we have three different notations for the holomorphic coefficients
of $f(\tau)$, i.e., the coefficients of $f_\Hol(\tau)$: 
\[
f_{(-k-1)}=f_{\{\frac1{\tau(n+1)}\}}=f_{[\Delta[k]]}=\Tr\left(
  f(\tau)\frac1{\tau(k+1)}\right).
\]
Similarly for the singular components we have two notations:
\[
f_{(k)}=f_{\{\tau(k)\}}=\Tr\left(f(\tau)\tau(k)\right).
\]
In fact we can let $F\in\hat\K$, see \eqref{eq:defhatK}, for instance 
\[
f_\Sing(\tau)=\sum_{n,k}f_{\{e_{n,k}^*\}}S(e_{n,k}\frac1\tau)),
\]
so that in \eqref{eq:rationalexpansionftauv} we have $v_{n,k}=f_{\{e_{n,k}^*\}}v$.

If $f(\tau)$, $g(\tau)$ are fields on $V$, they are in particular
$\End(V)$-valued distributions, and we can calculate their commutator
distribution: this is the distribution on $\K\otimes \K$ with value at $F\otimes G$
\[
[f(\tau_1),g(\tau_2)](F\otimes G)=\left(f(F)g(G)-g(G)f(F)\right).
\]
If we have an expansion (infinite, in general) as in \eqref{eq:expansiondelta},
\[
[f(\tau_1),g(\tau_2)]=\sum_{\ell\ge0}c_{\ell}(\tau_2)\Delta[\ell]_2\delta(\tau_1,\tau_2),
\]
then the coefficients $c_\ell(\tau)$ are also fields, by \eqref{eq:coefastrace}.

We say that fields are \emph{mutually rational} if the 
distribution $[f(\tau_{1}),g(\tau_{2})]$ has rational singularities,
i.e., is killed by multiplication by some $\mStwod(F)$, for $F\in\Czpol$.
By Lemma \ref{lem:rationaldistrtwovariables} this implies that we have
a finite sum
\begin{equation}
  [f(\tau_{1}),g(\tau_{2})]=\hat {\mathsf{D}}(f\otimes g)(\tau_2)\delta(\tau_1,\tau_1)=\sum
c_{n,k}(\tau_{2})(e_{n,k})_2\delta(\tau_{1},\tau_{2}),\label{eq:expansionrationalcommutator}
\end{equation}
where $\hat{\mathsf{D}}(f\otimes g)$ is a $W$-valued differential
difference operator (i.e., an element of $ W\hat H_T$) where $W=
\End(V)[[\tau^{\pm1}]]$ and the right hand side sum is finite with the
$c_{n,k}(\tau)$ $\End(V)$-valued distributions (in fact fields) on $K$.

\section{Normal Ordered Products and Dong's Lemma}
\label{sec:dong}
In this section we collect some properties of fields on a vector
space $V$ that we will
later need to investigate $H_T$-vertex algebras (see Definition
\ref{Def:HTvertexalgebra} below) and construct examples.

\begin{lem}\label{lem:productoffields}
  Let $f(\tau), g(\tau)$ be fields on $V$. Then also
  \begin{equation}
f_{\Hol}(\tau)g_{\Hol}(\tau),f_{\Sing}(\tau)g_{\Sing}(\tau),
f_{\Hol}(\tau)g_{\Sing}(\tau)\label{eq:fieldproducts}
\end{equation}
are fields on $V$.
\end{lem}

\begin{proof}
  By the product formula \eqref{eq:prodtauellm} we see that
  $f_{\Hol}(\tau)g_{\Hol}(\tau)$ is a well defined element of
  $\End(V)[[\tau^{\pm1}]]$. The singular part of this distribution
  is zero so this product is a field for trivial reasons.
  
  Now for any $v\in V$ the $V$-valued distribution $g_\Sing(\tau)v$ is
  rational, and we can multiply it with any $\End(V)$-valued
  distribution, see section \ref{sec:ExpDistributions}. So both
  $f_\Hol(\tau)g_\Sing(\tau)v$ and $f_\Sing(\tau) g_\Sing(\tau)v$ are
  well defined $V$-valued distributions on $K$. Their singular part is
  clearly rational.
\end{proof}

\begin{remark}
  In general the product $f_{\Sing}g_{\Hol}$ is not well defined.
\end{remark}

We use these results on products of holomorphic and singular parts to
define normal ordered products. 

\begin{defn}\label{def:normalorderedprod}
Let $f(\tau)$ and $g(\tau)$ be fields on $V\!$, and $F\in \K$.
The \emph{$F$-th normal ordered product} of $f(\tau)$ and $g(\tau)$ is
\begin{multline}
  \label{eq:Fnormalorderedproduct}
  f(\tau_2)_{\{F\}}g(\tau_2)=\Tr_{\tau_1}\left(
    f(\tau_1)g(\tau_2)\mathcal{R}^S(\tau_2)
    (F(\tau_1))-\right.\\\left. g(\tau_2)f(\tau_1)\mathcal{L}^S(\tau_1) (F(\tau_2))\right).
\end{multline}
\end{defn}

We define in particular the \emph{normal ordered product}, using Lemma
\ref{lem:TrDrhoRlambda} and the definitions of $\rho(\frac1\tau)$ and
$\lambda(\frac1\tau)$, \eqref{eq:defrholambda}, as
\begin{equation}
  \begin{split}
: f(\tau_2)g(\tau_2):&=f(\tau_2)_{\{\frac1\tau\}}g(\tau_2)=\Tr_{\tau_1}\left(
  f(\tau_1)g(\tau_2)\rho_{12}-g(\tau_2)f(\tau_1)\lambda_{12}\right)\\
&=f_{\Hol}(\tau_2)g(\tau_2)+g(\tau_2)f_{\Sing}(\tau_2).
\end{split}\label{eq:normalorderedproduct}
\end{equation}
More generally, for $h\in H_T$ define
\begin{equation}
  \begin{split}
    f(\tau_2)_{[h]}g(\tau_2)&=f(\tau_2)_{\{S(h)\frac1\tau\}}g(\tau_2)\\
    &=\Tr_{\tau_1}\left(
      f(\tau_1)g(\tau_2)h_2\rho_{12}-g(\tau_2)f(\tau_1)h_2\lambda_{12}\right)\\
    &=h_K\langle f_{\Hol}(\tau_2)\rangle g(\tau_2)+
    g(\tau_2)h_k\langle f_{\Sing}(\tau_2)\rangle.
\end{split}\label{eq:hnormalorderedproduct}
\end{equation}
So
\begin{equation}
f(\tau_2)_{[h]}g(\tau_2)=:h_K\langle f(\tau_2)\rangle
g(\tau_2):.\label{eq:hnormalproducthK}
\end{equation}
We also refer to $f(\tau_2)_{[h]}g(\tau_2)$ as the \emph{$h$th normal ordered
product} of $f$ and $g$. 

In case $F\in\Czpol$ the $F$-th normal ordered product is a component
of the commutator. Indeed for $F\in\Czpol$ we have
\[
\mathcal{R}^S(\tau_1)(F(\tau_2))=\mathcal{L}^S(\tau_2)(F(\tau_1))=\mStwod(F),
\]
so that in this case
\[
  f(\tau_2)_{\{F\}}g(\tau_2)=\Tr_{\tau_1}\left([f(\tau_1),g(\tau_2)]\mStwod(F)\right).
\]
We have in particular, see Lemma \ref{lem:coefficientDeltaexpansion},
\begin{equation}
[f(\tau_1),g(\tau_2)]=\sum_{\ell\ge0}f(\tau_2)_{\{\tau[\ell]\}}g(\tau_2)\Delta_2[\ell]\delta(\tau_1,\tau_2).
\label{eq:infiniteexpansioncommutator}
\end{equation}
To get the components of the finite expansion
\eqref{eq:expansionrationalcommutator} in case $f(\tau_1)$ and
$g(\tau_2)$ are mutually rational we need to take $F$ in $\hatCzpol$,
see \eqref{eq:defhatK}: we have
\begin{equation}
c_{n,k}(\tau)=f(\tau)_{\{e_{n,k}^*\}}g(\tau).\label{eq:finitecomponentrationalcommutator}
\end{equation}
By Lemma \ref{lem:productoffields} the $h$th normal ordered product of
fields is again a field, for all $h\in \hat H_T$. For $F\in\hatCzpol$
the $F$-th normal ordered products are components of commutators of
fields so also fields.
  
\begin{lem}\label{lem:LeibnizFnormalproduct}
  Let $V$ be an $\hat H_T$-module and $f(\tau)$, $g(\tau)$ fields on
  $V$, $F\in\K$ and $h\in\hat H_T$ with coproduct $\sum
  h^\prime\otimes h^{\prime\prime}$. Then
  \begin{enumerate}
  \item $\ad^V_h \left(
      f(\tau)_{\{F\}}g(\tau)\right)=\sum \left(\ad^V_{h^\prime}
      f(\tau)\right)_{\{F\}}
\left(\ad^V_{h^{\prime\prime}}g(\tau)\right)$.
\item $h_{\K}\left( f(\tau)_{\{F\}}g(\tau)\right)=\sum
  \left(h^\prime_{\K}f(\tau)\right)_{\{F\}}\left(h^{\prime\prime}_{\K}g(\tau)\right)$.
  \end{enumerate}
\end{lem}
\begin{lem}(Dong's Lemma)
  \label{lem:Dongslemma}
  Let $a(\tau),b(\tau),c(\tau)$ be fields on $V$ that are mutually
  rational. Then also the $F$th normal ordered product
  $a(\tau)_{\{F\}}b(\tau)$ and $c(\tau)$ are mutually rational, for
  all $F\in\hat K$.
\end{lem}

We first need a simple Lemma on polynomials.
\begin{lem}
  \label{lem:polynomiallemma}
  Let $F\in\mathbb C[u]$, $H\in\mathbb C[v]$. Put $x=u+v$, $y=u-v$.
  Then there are non zero $q,r\in\mathbb C[u,v]$ such that
\[
K=F(u)q(u,v)+H(v)r(u,v)\in \mathbb C[x],
\]
i.e., such that $K$ is independent of $y$.
\end{lem}
\begin{proof}
  Write for $A\in\mathbb C[u,v]=\mathbb C[x,y]$ 
\[
A=\sum_{i\ge0} A_i(x)y^i,\quad A_i\in \mathbb C[x],
\]
so that
  \begin{equation}
    \label{eq:FqHr}
F(u)q(u,v)+H(v)r(u,v)=\sum_{i,j\ge0}(F_i(x)q_j(x)+H_i(x)r_j(x))y^{i+j}.
  \end{equation}
  Let $I\subset \mathbb C[x]$ be the ideal generated by
  $F_0(x),H_0(x)$. The $y^0$-term of \eqref{eq:FqHr} reads
\[
F_0(x)q_0(x)+H_0(x)r_0(x)\in I.
\]
We will assume that also $q_0,r_0$ belong to $I$, but for the rest
they are arbitrary. Then assume we have
found $q_j(x),r_j(x)\in I$ such that the coefficient of $v^j$ in
\eqref{eq:FqHr} vanishes, for $j=1,2, \dots,s-1$. Then the $v^s$ term
reads
\begin{equation}
  \label{eq:v^sterm}
  F_0q_s+H_0r_s +\underset{i+j=s}{\sum_{i=1}^s\sum_{j=0}^{s-1}}F_iq_j+H_ir_j.
\end{equation}
By assumption the double sum belongs to $I$, so that we can find $q_s,
r_s$ making \eqref{eq:v^sterm} vanish. By multiplying the $q_j, r_j$,
$j<s$, by an element of $I$, if necessary, we can assume $q_s, r_s\in
I$, too, and the lemma follows by induction.
\end{proof}

\begin{proof}[Proof (of Dong's Lemma, cf., \cite{MR99f:17033}, Lemma 3.2 )]
Recall the following notation: if $f\in\Czpol$ we write
$f_{12}=\mStwod(f)\in\mathbb C[\tau_1,\tau_2]$. Now by multiplicativity
of $\mStwod$ and the fact that $\mStwod(\tau)=\tau_{12}=\tau_1-\tau_2$
we have $f_{12}=f(\tau_{12})\in\mathbb C[\tau_{12}]$. In the same way we
define for $h,g\in\Czpol$ the elements $h_{13}\in\mathbb
C[\tau_{13}]$, $g_{23}\in\mathbb C[\tau_{23}]$, with $\tau_{ij}=\tau_i-\tau_j$.

By definition of rationality there are $f,g,h\in\Czpol$ such that
  \begin{equation}
    f_{12}[a(\tau_1),b(\tau_2)]=0,\quad    g_{23}[b(\tau_2),c(\tau_3)]=0,\quad
    h_{13}[a(\tau_1),c(\tau_3)]=0.
\label{eq:defrational}
  \end{equation}
We need to find $k\in\Czpol$ such that
\begin{equation}
  \label{eq:normalorderedcommutator}
  k_{23}[a(\tau_2)_{\{F\}}b(\tau_2),c(\tau_3)]=0.
\end{equation}
It suffices to show that
\begin{equation}
  \label{eq:kA=kB}
  k_{23}A=k_{23}B,
\end{equation}
where
\begin{align*}
  A&=
  \left(a(\tau_1)b(\tau_2)\mathcal{R}^S(\tau_2)F(\tau_1)-b(\tau_2)a(\tau_1)
    \mathcal{L}^S(\tau_1)F(\tau_2)\right)c(\tau_3)\\
  B&= c(\tau_3)\left(a(\tau_1)b(\tau_2)\mathcal{R}^S
    (\tau_2)F(\tau_1)-b(\tau_2)a(\tau_1)\mathcal{L}^S(\tau_1)F(\tau_2)\right).
\end{align*}
Indeed, taking $\Tr_{\tau_1}$ of \eqref{eq:kA=kB} gives
\eqref{eq:normalorderedcommutator}.

Now note that there is a $G\in\Czpol$ such that $FG\in\hatCzpol$ and
in that case we have
\[
\mathcal{R}^S(\tau_2)(F(\tau_1)G(\tau_1))=\mathcal{L}^S(\tau_1)(F(\tau_2)G(\tau_2))
                                         =\mStwod(FG)=(FG)_{12}.
\]
This implies that
\begin{equation}
  \label{eq:ftauAzero}
  \begin{split}
    (fG)_{12}A&=f_{12}[a(\tau_1),b(\tau_2)]c(\tau_3)(FG)_{12}=0,\\
    (fG)_{12}B&=c(\tau_3)f_{12}[a(\tau_1),b(\tau_2)](FG)_{12}=0,
\end{split}
\end{equation}
Also note that
\begin{equation}
  \label{eq:hgA}
  h_{13}g_{23}A=  h_{13}g_{23}B.
\end{equation}
Next we use Lemma \ref{lem:polynomiallemma}, with $u=-\tau_{12}=\tau_2-\tau_1$ and
$v=\tau_{13}=\tau_1-\tau_3$, $F(u)=(fG)_{12}$ and $H(v)=h_{13}$. 
Then there are $q$ and
$r$ such that 
\[
K(\tau_{23})=(fG)_{12} q+h_{13}r\in \mathbb C[\tau_{23}].
\]
and \eqref{eq:kA=kB} follows from \eqref{eq:ftauAzero} and \eqref{eq:hgA}
if we take $k=Kg$, so that $k_{23}=K_{23}g_{23}$.
\end{proof}

\section{State-Field Correspondence and Vacuum Axioms}
\label{sec:StateFieldVacuumAxioms}
Let $V$ be a vector space. A \emph{State-Field correspondence} for $V$
is a map that associates to all $f\in V$ a field $Y(f,\tau)$ on $V$,
called the vertex operator of $f$. We also write $f(\tau)$ for
$Y(f,\tau)$. 

Now let $V$ contain a distinguished vector $1_V$, called the
\emph{vacuum} of $V$. We say that a state-field correspondence
\emph{satisfies the Vacuum Axioms} (with respect to $1_V$) in case
\begin{equation}
  \label{eq:Vertexopvacuum}
  Y(1_V,\tau)=1_{\End(V)},\quad  Y(f,\tau)1_V=f_{\Hol}(\tau)1_V,
\end{equation}
with the constant term (see \eqref{eq:deftau=0})
\begin{equation}
  \label{eq:constantterm}
  f_{\Hol}(\tau)1_V|_{\tau=0}=f.
\end{equation}
So we have $f_{(-1)}1_V=f$, for all $f\in V$.

We fix from now on a state-field correspondence satisfying the vacuum
axioms for some $1_V\in V$. Let $f\in V$ and define for all $h\in \hat
H_T$ an operator $h_V:V \to V$ by
\begin{equation}
h_V f=f_{[h]}1_V,\quad f\in V.\label{eq:defhsubV}
\end{equation}
At this point we don't know that, given a state-field correspondence,
$h\mapsto h_V$ gives an $\hat H_T$-module structure to $V$. This is an
extra assumption about the state-field correspondence.

\section{Definition of $H_T$-Vertex Algebras}
\label{sec:DefH_TVertex}

We will need the notion of an antipodal vertex operator $Y^S(f,\tau)$
associated to $f\in V$. More generally if $\mathcal{D}$ is a
distribution on $\K$ we put 
\begin{equation}
\mathcal{D}^S(F)=\mathcal{D}(S(F)),\label{eq:antipodalD}
\end{equation}
using the antipodal homomorphism on $\K$ given by $S(\tau)=-\tau$.

\begin{defn}\label{Def:HTvertexalgebra}
  An $H_T$-vertex algebra is a vector space $V$ with a vacuum vector
  $1_V\in V$ and a state-field correspondence $f\mapsto
  f(\tau)=Y(f,\tau)$, satisfying the vacuum axioms
  \eqref{eq:Vertexopvacuum}, \eqref{eq:constantterm} and such that
  furthermore
\begin{itemize}

\item ($H_T$-covariance) For all $f\in V$ and $h\in H_T$
  \begin{equation}
h_{\K} Y(f, \tau)= Y(h_V f,\tau).\label{eq:HTcovariance}
\end{equation}
\item (Skew-Symmetry)
  \begin{equation}
Y(f,\tau)g=\mathcal{R}_V(\tau)Y^S(g,\tau)f.\label{eq:Skewsymmetry}
\end{equation}
\item (Mutual Rationality) The vertex operators $Y(f,\tau_1)$ and
  $Y(g,\tau_2)$ are for all $f,g\in V$ mutually rational.
\end{itemize}
\end{defn}

Note that from \eqref{eq:HTcovariance} it follows that the map
$h\in H_T\mapsto h_V\in\End(V)$ gives $V$ the structure of $H_T$-module.

The simplest examples of $H_T$-vertex algebras are of course the
holomorphic vertex algebras introduced in section
\ref{sec:holVertex}, corresponding to commutative $H_T$-Leibniz
algebras, with as vertex operator
$Y(f,\tau)=\Yhol(f,\tau)=\mathcal{R}_V(\tau_1)f$.

\begin{remark}
  We define similarly the notion of $\hat H_T$-vertex algebra, by
  imposing instead of \eqref{eq:HTcovariance}  a $\hat H_T$-covariance axiom. It seems that the
  non holomorphic $H_T$-vertex algebras are all in fact $\hat
  H_T$-vertex algebras.
\end{remark}

\section{First Properties of $H_T$-vertex Algebras}
\label{sec:firstprop}
From now on $V$ is an $H_T$-vertex algebra, with vacuum $1_V$.

\begin{lem}\label{lem:firstprop}
  \begin{enumerate}
  \item For all $h\in H_T$ \label{item:VacuumHaction}
    \[
    h_V1_V=\epsilon(h)1_V.
    \]
\item For all $f\in V$\label{item:VertexonVacuumExponentialOperator}
  \[f_{\Hol}(\tau)1_V=\mathcal{R}_V(\tau)f.
\] 
  \end{enumerate}
\end{lem}
\begin{proof}
  From the vacuum axioms and the definition
  \eqref{eq:expansionholovertexoperatorY(f,tau)} it follows that
  \[
  Y(1_V,\tau)1_V=\sum (1_{V})_{[\Delta[n]]}1_V\Delta^*[n]=1_{V}.
  \]
  Hence by definition \eqref{eq:defhsubV} we find
\[
\Delta[n]_V 1_V=\left \lbrace 
  \begin{aligned}
    1_V &\quad \text{if $n= 0$}\\
    0   &\quad \text{if $n\ne 0$}
  \end{aligned}\right \rbrace=\epsilon(\Delta[n])1_V,
\]
and part \ref{item:VacuumHaction}. follows from the linearity of  the
counit $\epsilon$ and linearity of distributions. 

Similarly, from \eqref{eq:expansionholovertexoperatorY(f,tau)} and (\ref{eq:defhsubV}) it
follows that 
\[
f_{\Hol}(\tau)=\sum_{n\in\mathbb
  Z}\Delta[n]_V f\Delta^*[n]=\mathcal{R}_V(\tau)f,
\]
by definition of the exponential operator, see
 \eqref{eq:HTexponentials}.
\end{proof}

\begin{lem} (Leibniz Rule)\label{lem:Leibniz} For all $f, g\in V$ and $h\in H_T$
\[
h_V (Y(f,\tau)g)=\sum Y(h^\prime_V f,\tau)h^{\prime\prime}_Vg.
\]  
\end{lem}
\begin{proof}
  First some general remarks. The $H_T$-covariance of the antipodal
  vertex operators is 
  \begin{equation}
Y^S(hf,\tau)=S(h)_{\K} Y^S(f,\tau).\label{eq:antipodalHTcovariance}
\end{equation}
For the exponential operator $\mathcal{R}_V(\tau)$ (acting on a space
of $V$-valued distributions, say), we have the following identities:
  \begin{align}
    h_{\K}\mathcal{R}_V(\tau)&=\mathcal{R}_V(\tau)\sum
    h^\prime_Vh^{\prime\prime}_K,\label{eq:hKR} \\
    h_{V}\mathcal{R}_V(\tau)&=\mathcal{R}_V(\tau)h_V,\label{hVR}
  \end{align}
the first by the Leibniz rule in $K$ and Lemma
\ref{lem:HpropertiesExponentials}, and the second by the commutativity
of $H_T$.

Then, by the $H_T$-covariance of $Y(f,\tau)$ and skew symmetry we have
  \begin{align*}
\sum Y(h^\prime_V f,\tau)h^{\prime\prime}_Vg
     &= \sum h^\prime_{\K}Y(f,\tau)h_V^{\prime\prime }g&&\\
     &= \sum h^\prime_{\K}\mathcal{R}_V(\tau)Y^S(h_V^{\prime\prime}g,\tau)f&&\\
     &= \sum h^\prime_{\K}\mathcal{R}_V(\tau)S(h^{\prime\prime})_{\K}Y^S(g,\tau)f&
     &\text{by \eqref{eq:antipodalHTcovariance}}\\
     &= \mathcal{R}_V(\tau)\sum h^\prime_{V}h^{\prime\prime}_{\K}
     S(h_V^{\prime\prime\prime})_{\K}Y^S(g,\tau)f,&&\text{by \eqref{eq:hKR}}\\ 
     &=  \mathcal{R}_V(\tau)h_VY^S(g,\tau)f,&&\\
     &= h_V Y(f,\tau)g,&&
  \end{align*}
by elementary identities in the Hopf algebra $H_T$ and \eqref{hVR}.
\end{proof}

As $Y(f,\tau)$ is an $\End(V)$-valued distribution on $\K$, we get for
all $F\in\K$ (in fact $\in\hat \K$) a linear map
\[
f_{\{F\}}:V\to V,\quad f_{\{F\}}=\Tr(Y(f,\tau)F(\tau)).
\]
We call $f_{\{F\}}g$ the $F$-product of $f,g$. Then Lemma
\ref{lem:Leibniz} implies the Leibniz rule for $F$-products:
\[
h_V(f_{\{F\}}g)=\sum (h_V^\prime)_{\{F\}}(h^{\prime\prime}g).
\]
Another reformulation of the Leibniz rule is the following $\ad$-covariance
of vertex operators.
\begin{lem}\label{lem:adcovariancevertxop}
For all $f\in V$, $h\in H_T$
\[
h_{\K}Y(f,\tau)=\ad^V_h(Y(f,\tau)).
\]
\end{lem}
\begin{proof}
  For any $g\in V$ we have by $H_T$-covariance and elementary identities
  \begin{align*}
    h_KY(f,\tau)g&=Y(h_Vf,\tau)g &&\\
                 &=\sum
                 Y(h^\prime_Vf,\tau)h^{\prime\prime}_VS(h^{\prime\prime\prime})_Vg
                 &&\\
                 &=\sum
                 h^\prime_V(Y(f,\tau)S(h^{\prime\prime})_Vg)&&\text{by
                 Lemma \ref{lem:Leibniz},}\\
                 &=\ad^V_h(Y(f,\tau))g.
  \end{align*}
\end{proof}

\begin{lem}\label{lem:holovertexexponential}
  For all $f\in V$, $h\in H_T$
\[
h_{\K} f_{\Hol}(\tau)=\mathcal{R}_V(\tau)\circ f_{[h]}\circ
\mathcal{R}^{S}(\tau).
\]
\end{lem}

\begin{proof}
  By taking $\tau=0$ in the holomorphic part of Lemma
  {}~\ref{lem:adcovariancevertxop} we find
\[
f_{[h]}=\ad^V_h(f_{[1]}).
\]
Hence, by $H_T$-covariance and the definitions 
\begin{align*}
  h_{\K}f_{\Hol}(\tau)&=\ad^V_h (f_{\Hol}(\tau)) &&\\
                      &=\ad^V_h \left(\sum f_{[\Delta[n]]}\Delta^*[n] \right)&&\\
                      &=\sum \ad^V_{\Delta[n]}  f_{[h]}\Delta^*[n]
                      &&\text{by commutativity of $H_T$}\\
                      &=\ad^V_{\mathcal{R}(\tau)}(f_{[h]}) &&\\
                      &=\mathcal{R}_V(\tau)\circ
                      f_{[h]}\circ\mathcal{R}_V^{S}(\tau),&&
\end{align*}
by Lemma {}~\ref{lem:Adjointactionexponentials}.
\end{proof}

\begin{lem}
  For all $f\in V$
\[
\mathcal{R}^{\K}(\tau_1)Y(f,\tau_2)=\mathcal{R}_V(\tau_1)\circ
Y(f,\tau_2)\circ \mathcal{R}_V^{S}(\tau_1).
\]
\end{lem}

\begin{proof}
  This follows from Lemma \ref{lem:adcovariancevertxop} and Lemma
  \ref{lem:Adjointactionexponentials}.
\end{proof}

\begin{lem}\label{lem:Covariancecomponents}
  For all $f\in V$, $h\in H_T$ and $F\in\hat K$ we have
\[
(h_Vf)_{\{F\}}=f_{\{S(h)F\}}=\ad_h^V f_{\{F\}}.
\]
\end{lem}
\begin{proof}
By (\ref{eq:defFcomponent}) and $H_T$-covariance
(\ref{eq:HTcovariance}) we have
\begin{align*}
  Y(h_Vf,\tau)(F)&=(h_Vf)_{\{F\}}=h_K.Y(f,\tau)(F)\\
  &=Y(f,\tau)(S(h)F)
\intertext{by the definition of the $H_T$-action
    on distributions, and so}
(h_Vf)_{\{F\}}  &=f_{\{S(h)F\}},
\end{align*}
proving the first equality. For the second use Lemma
\ref{lem:adcovariancevertxop}:
\[
(h_Vf)_{\{F\}}=Y(h_Vf,\tau)(F)=h_K.Y(f,\tau)(F)=\ad_h^V Y(f,\tau)(F)=\ad_h^V f_{\{F\}}.
\]
\end{proof}
\begin{lem}\label{lem:tracefmStwod}
  Let $F\in\hatCzpol$ and let $f(\tau)=Y(f,\tau)$ be a vertex operator. Then
\[
\Tr_{\tau_1}\left(
  f(\tau_1)\mStwod(F)\right)=\mathcal{R}_V(\tau_2)\circ
f_{\{F\}}\circ\mathcal{R}_V^S(\tau_2)=
(\mathcal{R}_V(\tau_2)f)_{\{F\}}.
\]
\end{lem}
\begin{proof}
  We have $\mStwod(F)=\mathcal{R}^S_{\K}(\tau_2)(F(\tau_1))$. So
  \begin{align*}
    \Tr_{\tau_1}\left( f(\tau_1)\mStwod(F)\right)&=\Tr_{\tau_1}\left(
      f(\tau_1)\mathcal{R}^S_{\K}(\tau_2)F(\tau_1)\right)&&\\
    &=f_{\{\mathcal{R}^S_{\K}(\tau_2)F\}}&&\text{by \eqref{eq:defFcomponent}}\\
    &=\ad_{\mathcal{R}_{V}(\tau_2)}f_{\{F\}}&&\text{by Lemma
      \ref{lem:Covariancecomponents}}\\
    &=\mathcal{R}_V(\tau_2)f_{\{F\}}\mathcal{R}_V^S(\tau_2)&&\text{by
      Lemma \ref{lem:Adjointactionexponentials}}\\
    &=(\mathcal{R}_V(\tau_2)f)_{\{F\}}&&\text{again by
      Lemma \ref{lem:Covariancecomponents}}.
  \end{align*}
\end{proof}
\begin{lem}\label{lem:FnormalproductvacuumExponential}
  For all $f,g\in V$, $F\in \hat H_T$,
\[
f(\tau)_{F}g(\tau)1_V=\mathcal{R}_V(\tau)f_{F}g.
\]
\end{lem}
\begin{proof}
  We distinguish two cases: $F\in\hatCzpol$ and $F\in\Ksing$. For
  $F\in\Ksing$ we write $F=S(h)\frac1\tau$, so that
\[
f(\tau)_{\{F\}}g(\tau)=f(\tau)_{[h]}g(\tau)=:h_{\K}\langle
f(\tau)\rangle g(\tau)\rangle:,
\]
by \eqref{eq:hnormalproducthK}. Using the vacuum axioms and
\eqref{eq:normalorderedproduct} we find
\begin{align*}
  f(\tau)_{\{F\}}g(\tau)1_V&=h_{\K}\langle f_\Hol(\tau)\rangle
  g_\Hol(\tau)1_V&&\\
  &=\mathcal{R}_V(\tau)\circ f_{[h]}\circ \mathcal{R}_v^S(\tau)\circ
  \mathcal{R}_V(\tau)g,&&\text{by Lemma
    \ref{lem:holovertexexponential}}\\
  &=\mathcal{R}_V(\tau) f_{[h]}g&&\text{by section
    \ref{sec:InverseExp}}\\
  &=\mathcal{R}_V(\tau)f_{\{F\}}g&&\text{for $F=S(h)\frac1\tau$}.
\end{align*}
For $F\in\hatCzpol$ we write $F_{12}=\mStwod(F)$ and we have
\begin{align*}
  f(\tau)_{\{F\}}g(\tau)1_V&=\Tr_{\tau_1}\left([f(\tau_1),g(\tau_2)]F_{12}\right) 1_V&&\\
  &=[f_{\{F_{12}\}},g(\tau_2)]1_V&&\text{by Lemma
    \ref{lem:tracefmStwod}}\\
  &=f_{\{F_{12}\}}g(\tau_2)1_V&&\text{as
    $f_{\{G\}}1_V=0$, $G\in\hatCzpol$}\\
  &=\mathcal{R}_V(\tau_2)f_{\{F\}}\mathcal{R}_V^S(\tau_2)g_\Hol(\tau_2)1_V
  &&  \text{by vacuum axioms}\\
  &=\mathcal{R}_V(\tau_2)f_{\{F\}}g&&\text{by Lemma \ref{lem:InverseExp}}.
\end{align*}
\end{proof}
\begin{cor}\label{cor:Fnormalproductvacuum}
  For all $f,g\in V$, $F\in\hat \K$
\[
f(\tau)_{\{F\}}g(\tau)1_V|_{\tau=0}=f_{\{F\}}g.
\]
\end{cor}
\begin{proof}
  This follows from Lemma \ref{lem:FnormalproductvacuumExponential}
  and \eqref{eq:deftau=0}.
\end{proof}
\section{Uniqueness and Normal Ordered Products }
\label{sec:UniqNormOrdProd}

\begin{thm}
  (\textbf{Goddard's Uniqueness}) Let $G(\tau)$ be a field on an
  $H_T$-vertex algebra that is mutually rational to all vertex
  operators $Y(f,\tau)$ of $V$ and such that 
\[
G(\tau)1_V=\mathcal{R}_V(\tau)g,
\]
for some $g\in V$. Then 
\[
G(\tau)=Y(g,\tau).
\]
\end{thm}

The proof is the same as in the usual case, using the vacuum axioms,
skew symmetry and rationality, cf., \cite{MR99f:17033}, Thm. 4.4.

Now we will express the vertex operators of an $H_T$-vertex algebra in
terms of normal ordered products, see Definition \ref{def:normalorderedprod}.

\begin{lem}\label{lem:Yhomvertexproducts}
  For $f,g\in V$ and $F\in\hat\K$ we have
  \[
  Y(f_{\{F\}}g,\tau)=f(\tau)_{\{F\}}g(\tau).
  \]
More generally, if $f^1,f^2\dots,f^k\in V$ and
$F_1,F_2,\dots,F_k\in\hat\K$,
then
\[
Y(f^{1}_{\{F_1\}}f^{2}_{\{F_2\}}\dots
  f^{k}_{\{F_k\}}1_V,\tau)=f^{1}(\tau)_{\{F_1\}}f^{2}(\tau)_{\{F_2\}}\dots  f^{k}(\tau)_{\{F_k\}}1_V.
\]
\end{lem}

We define the normal ordered product of more factors as usually from
right to left:
\begin{multline*}
  f^{1}(\tau)_{\{F_1\}}f^{2}(\tau)_{\{F_2\}}\dots
  f^{k}(\tau)_{\{F_k\}}g(\tau)= \bigl(f^{1}(\tau)_{\{F_1\}}
  (f^{2}(\tau)_{\{F_2\}}\dots\\\dots (f^{k}(\tau)_{\{F_k\}}g(\tau))\dots)\bigr)
\end{multline*}
\begin{proof}
  We have, by the vacuum axioms, Lemma \ref{lem:firstprop} and Lemma
  \ref{lem:holovertexexponential} 
\[
Y(f_{\{F\}}g,\tau)1_V=\mathcal{R}_V(\tau)f_{\{F\}}g=f(\tau)_{\{F\}}g(\tau))1_V,
\]
so that the first part of the Lemma follows from Goddard's
uniqueness. The second part follows by induction.
\end{proof}

\begin{lem}(Borcherds OPE)
  Let $V$ be an $H_T$-vertex algebra, $f,g\in V$. Then
  \begin{align*}
    [Y(f,\tau_1),Y(g,\tau_2) 
    &=\sum_{n,k}f_{\{e_{n,k}^*\}}g(\tau_2)
    (e_{n,k})_2\delta(\tau_1,\tau_2)\\
    &=\sum_{\ell\ge0}f_{\{\tau(\ell)\}}g(\tau_2)
    (\Delta[\ell])_2\delta(\tau_1,\tau_2).    
  \end{align*}
\end{lem}

Note that the first sum on the RHS in the Lemma will be finite, the
second infinite, in general.

\begin{proof}
  The first equality follows from the mutual rationality of vertex operators by
  combining \eqref{eq:expansionrationalcommutator},
  \eqref{eq:finitecomponentrationalcommutator} and the previous Lemma
  \ref{lem:Yhomvertexproducts}. The second equality follows from
  \eqref{eq:infiniteexpansioncommutator} combined with Lemma \ref{lem:Yhomvertexproducts}.
\end{proof}
\section{Alternative Axiomatization}
\label{sec:AltAxiom}

We took as one of the axioms of an $H_T$-vertex algebra the
$H_T$-covariance axiom
\[
Y(h_Vf,\tau)=h_K.Y(f,\tau),
\]
and derived in Lemma \ref{lem:adcovariancevertxop} the
$\ad$-covariance
\begin{equation}
Y(h_Vf,\tau)=\ad_h^VY(f,\tau).\label{eq:adhcovariance}
\end{equation}
This gives a set of axioms similar to those for an $H_T$-vertex
Poisson algebra, see Definition \ref{defn:VertexPoisson}, and of an
$H_T$-conformal algebra, see Definition \ref{def:HTconformal}.
However, these axioms are not very efficient, in particular we assumed
the skew-symmetry \eqref{eq:Skewsymmetry} property. In this section we
give a different set of axioms, including the $\ad$-covariance
\eqref{eq:adhcovariance}, which are easier to check and are analogous
to those in \cite{MR99f:17033}. In particular we don't assume
skew-symmetry as an axiom, but derive it.

\begin{prop}
  Let $V$ be an $H_T$-module with vacuum vector $1_V$ and a
  state-field correspondence satisfying the vacuum axioms
  \eqref{eq:Vertexopvacuum}, \eqref{eq:constantterm}, and such that
  \begin{Myenumerate}
  \item (Compatibility of State-Field Correspondence with
  $H_T$-action) For all $f\in V$ and $h\in H_T$ we have
\[
h_V.f=f_{[h]}1_V,
\]
where the right handside is defined by \eqref{eq:DefhsubV}.
\item ($\ad_h$-covariance) For all $f\in V$ and $h\in H_T$
\[
\ad_h^V Y(f,\tau)=h_K.Y(f,\tau).
\]
\item (Mutual Rationality) The fields $Y(f,\tau_1)$ and
  $Y(g,\tau_2)$ are for all $f,g\in V$ mutually rational.
  \end{Myenumerate}
Then $V$ is an $H_T$-vertex algebra as defined in Definition \ref{Def:HTvertexalgebra}.
\end{prop}

\begin{proof}
  We need to prove skew-symmetry \eqref{eq:Skewsymmetry} and
  $H_T$-covariance \eqref{eq:HTcovariance}.

First we note that
\begin{equation}
  \label{eq:vacuumexponential}
  Y(f,\tau)1_V=\mathcal{R}_V(\tau)f,
\end{equation}
because the proof of Lemma \ref{lem:firstprop} still works. We also
have by $\ad_h$-co\-vari\-ance \eqref{eq:adhcovariance} and Lemma \ref{lem:Adjointactionexponentials}
\begin{equation}
  \label{eq:KVexponentials}
  \mathcal{R}_{\K}^S(\tau_1)Y(g,\tau_2)=\mathcal{R}_V^S(\tau_1)Y(g,\tau_2)\mathcal{R}_V(\tau_1).
\end{equation}
Finally, for any \emph{holomorphic} distribution $\mathcal{D}_\Hol$ we
have
\[
\mathcal{R}^S(\tau_1)\mathcal{D}_\Hol(\tau_2)|_{\tau_1=0}=\mathcal{D}_\Hol^S(\tau_1), 
\]
where $\mathcal{D}^S$ is the antipodal distribution defined by
\eqref{eq:antipodalD}. 

Now by mutual rationality we have, for some $F\in\Czpol$
\begin{align*}
  F_{12}Y(f,\tau_1)Y(g,\tau_2)1_V 
  &=  F_{12}Y(g,\tau_2)Y(f,\tau_1)1_V&&\\
  &=  F_{12}Y(g,\tau_2)\mathcal{R}_V(\tau_1)f&&\text{by
  \eqref{eq:vacuumexponential}}\\
  &=  F_{12}\mathcal{R}_{\K}^S(\tau_1)\mathcal{R}_V(\tau_1)Y(g,\tau_2)f&&\text{by \eqref{eq:KVexponentials}}.
\end{align*}
Now the left hand side is manifestly holomorphic in $\tau_2$, whereas
in the right hand side $F_{12}\mathcal{R}_{\K}^S(\tau_1)Y(g,\tau_2)f$
can be assumed to holomorphic by suitable choice of $F$. So we can put
$\tau_2=0$ and cancel $F_{12}|_{\tau_2=0}$ so get skew-symmetry
\eqref{eq:Skewsymmetry}. 

Next note that we have
\[
f(\tau)_{\{F\}}g(\tau)1_V=\mathcal{R}_V(\tau)f_{\{F\}}g,
\] as the proof of \ref{lem:tracefmStwod} just uses
$\ad_h$-covariance of fields. Now we use Goddard's uniqueness in our
situation: if $V$ has a state-field correspondence $Y$ as in the
statement of the theorem and we have a field $G(\tau)$ mutually
rational to all $Y(f,\tau)$ such that for some $g\in V$
\[
Y(g,\tau)1_V=G(\tau)1_V,
\]
then in fact $G(\tau)=Y(g,\tau)$. We use this to conclude that for
all$f,g\in V$
\begin{equation}
  \label{eq:Fproducthomomorphism}
Y(f_{\{F\}}g,\tau)=f(\tau)_{\{F\}}g(\tau).  
\end{equation}
In particular, for $F=S(h)\frac1\tau$, $h\in H_T$ and $g=1_V$ we find
on the one hand, by compatibility of state-field correspondence with
the $H_T$-action
\[
Y(f_{\{S(h)\frac1\tau\}}1_V,\tau)=Y(f_{[h]}1_V,\tau),
\]
and on the other hand, by~\eqref{eq:Fproducthomomorphism}
\[
Y(f_{\{S(h)\frac1\tau\}}1_V,\tau)=:h_K\langle f(\tau)\rangle
Y(1_V,\tau):,
\]
so that $H_T$-covariance \eqref{eq:HTcovariance} follows from the
vacuum axioms.
\end{proof}

\section{Existence}
\label{sec:existence}
\begin{prop}\cite{MR1334399}\label{prop:Existence}
  Let $V$ be an $H_T$-module with distinguished vector $1_V$ and a
  collection of fields $f^\alpha(\tau)$ such that
  \begin{Myenumerate}
  \item $\ad_h^V f^\alpha(\tau)=h_{\K}.f(\tau)$, $h\in H_T$.
  \item $h 1_V=\epsilon(h)1_V$, $h\in H_T$.
  \item $f^\alpha(\tau)1_V=f^\alpha_\Hol(\tau)1_V$, and, defining
    $f^\alpha\in V$ by $f^\alpha(\tau)1_V|_{\tau=0}=f^\alpha$, the map
    $f^\alpha(\tau)\mapsto f^\alpha$ is injective.
  \item All fields $f^\alpha(\tau)$ are mutually rational.
  \item $V$ is spanned by
  $f^{\alpha_1}_{\{F_1\}}f^{\alpha_2}_{\{F_2\}}\dots
  f^{\alpha_k}_{\{F_k\}}1_V$, for $F_1,F_2,\dots, F_K\in\hat\K$.
  \end{Myenumerate}
Then $V$ is a vertex algebra, with state field correspondence given by
$F$-normal ordered products:
\[
Y(f^{\alpha_1}_{\{F_1\}}f^{\alpha_2}_{\{F_2\}}\dots
  f^{\alpha_k}_{\{F_k\}}1_V,\tau)=f^{\alpha_1}(\tau)_{\{F_1\}}f^{\alpha_2}(\tau)_{\{F_2\}}\dots
  f^{\alpha_k}(\tau)_{\{F_k\}}1_V.
\]
\end{prop}

The proof in \cite{MR99f:17033}, Thm. 4.5 works in our situation,
mutatis mutandis, using the alternative axioms of section \ref{sec:AltAxiom}.

\section{Affine $H_T$-vertex algebras}
\label{sec:AffHTVertexAlg}
Let $\mathfrak g$ be a Lie algebra. Recall from Example
\ref{ex:affineLiealgebra} the notion of the conformal affinization
$\mathcal{L}\mathfrak g$ associated to the conformal algebra $\Cg$,
see Example \ref{ex:AffineHTconformalalg}. Define a decomposition
$\Lg=\Lgcr\oplus\Lgann$ in creation and annihilation Lie subalgebras (and $\hat
H_T$-submodules), where
\begin{align*}
\Lgcr&=\mathrm{Span}(X_{\{p\}}; X\in\mathfrak g,p\in \Ksing),\\
\Lgann&=\mathrm{Span}(X_{\{p\}}; X\in\mathfrak g,p\in \Czpol).
\end{align*}
Define a 1-dimensional trivial $\Lgann$ module $\mathbb C_0=\mathbb C m_0$ and let
$V(\mathfrak g)$ be the induced $\Lg$-module:
\[
V(\mathfrak g)=\mathcal{U}(\Lg)\otimes_{\mathcal{U}(\Lgann)}\mathbb
C_0\simeq \mathcal{U}(\Lgcr)v_0,\quad v_0=1\otimes m_0.
\]
Introduce an action of $H_T$ on $V(\mathfrak g)$ by putting
$hv_0=\epsilon(h)v_0$ and extending the $H_T$ action
(\ref{eq:defHactionLM}) on $\Lg$ to $\mathcal{U}(\Lg)$ by the Leibniz
rule. Define for $\{X^\alpha\}$ a basis for $\mathfrak g$ currents
$X^\alpha(\tau)$ as in (\ref{eq:DefCurrentsLC}). These are fields on
$V(\mathfrak g)$ that are mutually rational, and the action of the
coefficients of these fields on the vacuum $v_0$ obviously produces a
spanning set for $V(\mathfrak g)$. Also note that, if
$X^\alpha(\tau)v_0|_{\tau=0}=X^\alpha_0$, the the map
$X^\alpha(\tau)\mapsto X^\alpha_0$ is injective.
Furthermore, we have $\ad_h$-covariance of these fields.

\begin{lem}
  For all $h\in H_T$ and all $X^\alpha$ we have $\ad_h^{V(\mathfrak
    g)} X^\alpha(\tau)=h_{\K} X^\alpha(\tau)$.
\end{lem}

\begin{proof}
  For all $v\in V=V(\mathfrak g)$ we have
  \begin{align*}
    \ad_h^{V} X^\alpha(\tau_1)&=\sum h^\prime_V\circ X^\alpha(\tau_1)\circ S(h^{\prime\prime})_Vv\\
    &=\sum \Tr_{0}\left(X^\alpha\otimes
    S(h^\prime)_{\tau_0}\delta(\tau_0,\tau_1)\right)
    h^{\prime\prime}_VS(h^{\prime\prime\prime})_Vv\\
    &=\sum \Tr_{\tau_0}\left(X^\alpha\otimes
    S(h^\prime)_{0}\delta(\tau_0,\tau_1)\right)
    \epsilon(h^{\prime\prime})1_Vv\\
    &=\sum \Tr_{\tau_0} \left(X^\alpha\otimes S(h)_{0}\delta(\tau_0,\tau_1)\right)v\\
    &=\sum \Tr_{\tau_0}(X^\alpha\otimes h_{1}\delta(\tau_0,\tau_1))v\\
    &=h_{\K} X^\alpha(\tau_1)v
  \end{align*}
\end{proof}
Then is follows from Proposition \ref{prop:Existence} that
$V(\mathfrak g)$ is an $H_T$-vertex algebra, called the affine
$H_T$-vertex algebra of $\mathfrak g$.

\section{Toda Vertex Algebra}
\label{sec:TodaVertexAlg}
Recall from Example \ref{ex:TodaLiealg} the Lie algebra $\LToda$
associated to the Toda conformal algebra of Example
\ref{ex:CToda}. Define as before a decomposition
$\LToda=\LTodacr\oplus\LTodaann$ in creation and annihilation Lie subalgebras (and $\hat
H_T$-submodules), where
\begin{align*}
{\LTodacr}&=\mathrm{Span} (B_{\{p\}}, C_{\{p\}}; p\in \Ksing),\\
\LTodaann&=\mathrm{Span}(B_{\{p\}}, C_{\{p\}}; p\in \Czpol).
\end{align*}
Define a 1-dimensional trivial $\LTodaann$ module $\mathbb C_0=\mathbb
C m_0$ and let
$\VToda$ be the induced $\LToda$-module:
\[
\VToda=\mathcal{U}(\LToda)\otimes_{\mathcal{U}(\LTodaann)}\mathbb
C_0\simeq \mathcal{U}(\LTodacr)v_0,\quad v_0=1\otimes m_0.
\]
In exactly the same as for the affine vertex algebra we see that
$\VToda$ is an $H_T$-vertex algebra, called the  Toda vertex
algebra. It is a quantization of the Toda vertex Poisson algebra
discussed in Chapter \ref{chap:HamstrucInfToda}.

Note that $\VToda$ is non commutative, for the multiplication
$f\otimes g\mapsto f_{\{1\}}g$, because the subalgebra $\LTodacr$ is
non-Abelian, see (\ref{eq:BCcommutator}).

\begin{thebibliography}{FKRW95}

\bibitem[Abe80]{MR83a:16010}
Eiichi Abe, \emph{Hopf algebras}, Cambridge Tracts in Mathematics, vol.~74,
  Cambridge University Press, Cambridge, 1980, Translated from the Japanese by
  Hisae Kinoshita and Hiroko Tanaka. \MR{83a:16010}

\bibitem[Bor98]{MR1653021}
Richard~E. Borcherds, \emph{Vertex algebras}, Topological field theory,
  primitive forms and related topics (Kyoto, 1996), Birkh\"auser Boston,
  Boston, MA, 1998, pp.~35--77. \MR{1 653 021}

\bibitem[CP95]{MR96h:17014}
Vyjayanthi Chari and Andrew Pressley, \emph{A guide to quantum groups},
  Cambridge University Press, Cambridge, 1995, Corrected reprint of the 1994
  original. \MR{96h:17014}

\bibitem[Dic03]{MR1964513}
L.~A. Dickey, \emph{Soliton equations and {H}amiltonian systems}, second ed.,
  Advanced Series in Mathematical Physics, vol.~26, World Scientific Publishing
  Co. Inc., River Edge, NJ, 2003. \MR{MR1964513 (2004c:37160)}

\bibitem[FKRW95]{MR1334399}
Edward Frenkel, Victor Kac, Andrey Radul, and Weiqiang Wang, \emph{{$\mathcal
  W\sb {1+\infty}$} and {$\mathcal W(\mathfrak g\mathfrak l\sb N)$} with
  central charge {$N$}}, Comm. Math. Phys. \textbf{170} (1995), no.~2,
  337--357. \MR{MR1334399 (96i:17024)}

\bibitem[FLM88]{MR90h:17026}
Igor Frenkel, James Lepowsky, and Arne Meurman, \emph{Vertex operator algebras
  and the {M}onster}, Academic Press Inc., Boston, MA, 1988. \MR{90h:17026}

\bibitem[Ike94]{MR1304086}
Kaoru Ikeda, \emph{The {P}oisson structure on the coordinate ring of discrete
  {L}ax operator and {T}oda lattice equation}, Adv. in Appl. Math. \textbf{15}
  (1994), no.~4, 379--389. \MR{MR1304086 (96g:58076)}

\bibitem[Kac98]{MR99f:17033}
Victor Kac, \emph{Vertex algebras for beginners}, second ed., American
  Mathematical Society, Providence, RI, 1998. \MR{99f:17033}

\bibitem[Kup85]{MR86m:58070}
B.~A. Kuperschmidt, \emph{Discrete {L}ax equations and differential-difference
  calculus}, Ast\'erisque (1985), no.~123, 212. \MR{86m:58070}

\bibitem[Li]{math.QA/0209310}
Haisheng Li, \emph{{Vertex algebras and vertex poisson algebras}},
  arXiv:math.QA/0209310.

\bibitem[Sny]{math.QA/9904104}
Craig~T. Snydal, \emph{{Equivalence of Borcherds G-Vertex Algebras and
  Axiomatic Vertex Algebras}}, arXiv:math.QA/9904104.

\bibitem[Tod89]{MR971987}
Morikazu Toda, \emph{Theory of nonlinear lattices}, second ed., Springer Series
  in Solid-State Sciences, vol.~20, Springer-Verlag, Berlin, 1989. \MR{MR971987
  (89h:58082)}

\end{thebibliography}
\bibliographystyle{amsalpha}
\def\cprime{$'$}
\providecommand{\bysame}{\leavevmode\hbox to3em{\hrulefill}\thinspace}
\providecommand{\MR}{\relax\ifhmode\unskip\space\fi MR }
\providecommand{\MRhref}[2]{%
  \href{http://www.ams.org/mathscinet-getitem?mr=#1}{#2}
}
\providecommand{\href}[2]{#2}

\end{document}